Universidade Federal de Minas Gerais

Programa de Pós-Graduação em Engenharia Elétrica

# Maximum *a Posteriori* Joint State Path and Parameter Estimation in Stochastic Differential Equations

Dimas Abreu Dutra

**Universidade Federal de Minas Gerais**

**Escola de Engenharia**

**Programa de Pós-Graduação em Engenharia Elétrica**

# MAXIMUM A POSTERIORI JOINT STATE PATH AND PARAMETER ESTIMATION IN STOCHASTIC DIFFERENTIAL EQUATIONS

Dimas Abreu Dutra

Tese de Doutorado submetida à Banca Examinadora designada pelo Colegiado do Programa de Pós-Graduação em Engenharia Elétrica da Escola de Engenharia da Universidade Federal de Minas Gerais, como requisito para obtenção do Título de Doutor em Engenharia Elétrica.

Orientador: Prof. Luis Antonio Aguirre

## "Maximum A Posteriori Joint State Path and Parameter Estimation in Stochastic Differential Equations"

**Dimas Abreu Dutra**

**Tese de Doutorado submetida à Banca Examinadora designada pelo Colegiado do Programa de Pós-Graduação em Engenharia Elétrica da Escola de Engenharia da Universidade Federal de Minas Gerais, como requisito para obtenção do grau de Doutor em Engenharia Elétrica.**

**Aprovada em 13 de agosto de 2014.**

**Por:**

**Prof. Dr. Luis Antonio Aguirre**
**DELT (UFMG) - Orientador**

**Prof. Dr. Bruno Otávio Soares Teixeira**
**DELT (UFMG)**

**Prof. Dr. Eduardo Mazoni Andrade Marçal Mendes**
**DELT (UFMG)**

**Prof. Dr. Alexandre Rodrigues Mesquita**
**DELT (UFMG)**

**Prof. Dr. Hani Camille Yehia**
**DELT (UFMG)**

**Prof. Dr. Oswaldo Luiz do Valle Costa**
**Eletrônica (USP-Poli)**

**Prof. Dr. Marcelo Dutra Fragoso**
**LNCC (LNCC)**

Having all the answers just means you've been asking boring questions.

Emily Horn and Joey Comeau
A Softer World #869, *The freedom of uncertainty*

# Abstract

A wide variety of phenomena of engineering and scientific interest are of a continuous-time nature and can be modeled by stochastic differential equations (SDEs), which represent the evolution of the uncertainty in the states of a system. For systems of this class, some parameters of the SDE might be unknown and the measured data often includes noise, so state and parameter estimators are needed to perform inference and further analysis using the system state path. One such application is the flight testing of aircraft, in which flight path reconstruction or some other data smoothing technique is used before proceeding to the aerodynamic analysis or system identification.

The distributions of SDEs which are nonlinear or subject to non-Gaussian measurement noise do not admit tractable analytic expressions, so state and parameter estimators for these systems are often approximations based on heuristics, such as the extended and unscented Kalman smoothers, or the prediction error method using nonlinear Kalman filters. However, the Onsager–Machlup functional can be used to obtain *fictitious* densities for the parameters and state-paths of SDEs with analytic expressions.

In this thesis, we provide a unified theoretical framework for maximum *a posteriori* (MAP) estimation of general random variables, possibly infinite-dimensional, and show how the Onsager–Machlup functional can be used to construct the joint MAP state-path and parameter estimator for SDEs. We also prove that the minimum energy estimator, which is often thought to be the MAP state-path estimator, actually gives the state paths associated to the MAP *noise* paths. Furthermore, we prove that the discretized MAP state-path and parameter estimators, which have emerged recently as powerful alternatives to nonlinear Kalman smoothers, converge hypographically as the discretization step vanishes. Their hypographical limit, however, is the MAP estimator for SDEs when the trapezoidal discretization is used and the minimum energy estimator when the Euler discretization is used, associating different interpretations to each discretized estimate.

Example applications of the proposed estimators are also shown, with both simulated and experimental data. The MAP and minimum energy estimators are compared with each other and with other popular alternatives.

# Resumo

Uma grande variedade de fenômenos de interesse para engenharia e ciência são a tempo contínuo por natureza e podem ser modelados por equações diferenciais estocásticas (EDEs), que representam a evolução da incerteza nos estados do sistema. Para sistemas dessa classe, alguns parâmetros da EDE podem ser desconhecidos e os dados coletados frequentemente incluem ruídos, de modo que estimatores de esstados e parâmetros são necessários para realizar inferência e análises adicionais usando a trajetória dos estados do sistema. Uma dessas aplicações é em ensaios em voo de aeronaves, para os quais reconstrução de trajetória de voo ou outras técnicas de suavização são utilizadas antes de se proceder para análise aerodinâmica ou identificação de sistemas.

As distribuições de EDEs não lineares ou sujeitas a ruído de medição não Gaussiano não admitem expressões analíticas utilizáveis, o que leva a estimadores de estados e parâmetros para esses sistemas a basearem-se em heurísticas como os suavizadores de Kalman estendido e *unscented*, ou o método de predição de erro utilizando filtros de Kalman não lineares. No entanto, o funcional de Onsager–Machlup pode ser utilizado para obter densidades fictícias conjuntas para trajetórias de estado e parâmetros de EDEs com expressões analíticas.

Nesta tese, um arcabouço teórico unificado é desenvolvido para estimação máxima *a posteriori* (MAP) de variáveis aleatórias genéricas, possivelmente infinito-dimensionais, e é mostrado como o funcional de Onsager–Machlup pode ser utilizado para a construção do estimador MAP conjunto de trajetórias de estado e parâmetros de EDEs. Também é provado que o estimador de mínima energia, comumente confundido com com o estimador de MAP, obtém as trajetórias de estado associadas às trajetórias de *ruído* MAP. Além disso, é provado que os estimadores conjuntos de trajetória de estados e parâmetros MAP discretizados, que emergiram recentemente como alternativas poderosas para os estimadores de Kalman não lineares, convergem hipograficamente à medida que o passo de discretização diminue. O seu limite hipográfico, no entanto, é o estimador MAP para EDEs quando a discretização trapezoidal é utilizada e o estimador de mínima energia quando a discretização de Euler é utilizada, associando interpretações diferentes a cada estimativa discretizada.

Exemplos de aplicações dos estimadores propostos são apresentadas com dados simulados e experimentais, nas quais os estimadores MAP e de mínima energia são comparados entre si e com alternativas mais bem sedimentadas.

# Contents

# Notation

> A foolish consistency is the hobgoblin of little minds.
>
> Ralph Waldo Emerson, *Self-Reliance*

Here is a list of the mathematical typography and notation used throughout this thesis. Most are referenced on first use and are standard in the literature, still they are collected here for easy reference. To avoid a notational overload, many representations have a simpler form omitting parameters which can be clearly deduced from the context.

## Typographical conventions

We begin with some typographical conventions used to represent different mathematical objects. Note that these conventions are sometimes broken when they become cumbersome or deviate from the standard convention of the literature.

*Identifing subscripts* are typeset upright, e.g., $x_{\mathrm{a}}, x_{\mathrm{b}}$.

*Matrices* are typeset in uppercase bold, e.g., $\boldsymbol{A}, \boldsymbol{B}, \boldsymbol{\Gamma}, \boldsymbol{\Phi}$.

*Random variables* are typeset in uppercase while values they might take are represented in lowercase, e.g., if $X\colon \Omega \to \mathbb{R}$ is a $\mathbb{R}$-valued random variable over the probability space $(\Omega, \mathcal{E}, P)$, then $x \in \mathbb{R}$ can be used to represent specific values it can take. The dependence on the outcome $\omega$ will be omitted when unambiguous.

*Sets* are typeset in uppercase blackboard bold, e.g., $\mathbb{A}, \mathbb{B}$.

*k-chains* are typeset in lowercase blackboard bold, e.g., $\mathbb{a}, \mathbb{b}$.

*Time indices* of stochastic processes are typeset as subscripts when unambiguous, e.g., if $X\colon \mathbb{R} \times \Omega \to \mathbb{R}$ is an $\mathbb{R}$-valued stochastic process over the probability space $(\Omega, \mathcal{E}, P)$, then $X(t, \omega)$ can be written as $X_t(\omega)$ or $X_t$.

*Topological spaces* are typeset in uppercase calligraphic, e.g., $\mathcal{A}, \mathcal{B}$.

## General symbols

$\mathbb{N} := \{1, 2, 3, \dots\}$ is the set of natural numbers (strictly positive integers).

$\mathbb{R}$ is the set of real numbers.

$\mathbb{R}_{\geq 0} := \{x \in \mathbb{R} \,|\, x \geq 0\}$ is the set of nonnegative real numbers.

$\mathbb{R}_{>0} := \{x \in \mathbb{R} \,|\, x > 0\}$ is the set of strictly positive real numbers.

$\overline{\mathbb{R}} := \mathbb{R} \cup \{-\infty, \infty\}$ is the extended real number line.

## Linear algebra

$\boldsymbol{A}^{-1}$ is the inverse of the matrix $\boldsymbol{A}$.

$\boldsymbol{A}^{\top}$ is the transpose of the matrix $\boldsymbol{A}$.

$\boldsymbol{I}_n$ is an $n \times n$ identity matrix. The subscript indicating the size might be dropped if it can be deduced from the context.

$\boldsymbol{A}^{(ij)}$ is the element in the $i$th row and $j$th column of the matrix $\boldsymbol{A}$.

$a^{(i)}$ is the $i$th element of the vector $a$.

## Probability and measure theory

$(\Omega, \mathcal{E}, P)$ is a standard probability space (see Ikeda and Watanabe, 1981, Defn. 1.3.3) on which all random variables are defined.

$\omega \in \Omega$ is the random outcome.

$\mathcal{B}_{\mathcal{X}}$ is the Borel $\sigma$-algebra of the topological space $\mathcal{X}$, i.e., the $\sigma$-algebra of subsets of $\mathcal{X}$ generated by the topology of $\mathcal{X}$.

$\{\mathcal{E}_t\}_{t \geq 0}$ is a filtration on the probability space $(\Omega, \mathcal{E}, P)$, i.e., $\mathcal{E}_t \subset \mathcal{E}_s \subset \mathcal{E}$ for all $t, s \in [0, \infty)$ such that $t \leq s$.

$[\![A, B]\!]$ is the process of quadratic covariation between $A$ and $B$.

$\mathrm{supp}(\mu)$ is the support of a measure $\mu$ over a measurable space $(\mathcal{X}, \mathcal{B}_{\mathcal{X}})$, i.e., the set $\mathbb{A} \in \mathcal{B}_{\mathcal{X}}$ of all points whose every open neighbourhood has strictly positive $\mu$-measure (cf. Ikeda and Watanabe, 1981, Sec. 6.8).

$I_{\mathbb{A}}(x)$ is the indicator function of the set $\mathbb{A}$, i.e., $I_{\mathbb{A}}(x) := 1$ if $x \in \mathbb{A}$, 0 if $x \notin \mathbb{A}$.

In addition, we will use the piece of jargon "for $P$-almost all $\omega \in \mathbb{E}$" to say that a property holds for all but a zero-measure subset of an event $E \in \mathcal{E}$, i.e., when the property holds for all $\omega \in \mathbb{E} \backslash \mathbb{N}$, where $\mathbb{N} \in \mathcal{E}$ is a $P$-null event $P(\mathbb{N}) = 0$.

## Analysis

$|x|$ is the Euclidean norm of $x \in \mathbb{R}^d$.

$\|x\|$ is a norm of $x$. If unambiguous, the norm shall be inferred from the space of $x$.

$\|\!|f|\!\|:=\max_{x\in\mathcal{X}}\|f(x)\|_{\mathcal{Y}}$ is the supremum norm of a function $f\colon\mathcal{X}\to\mathcal{Y}$, also known as the infinity or uniform norm.

$\langle a,b\rangle$ is the inner product between $a$ and $b$. If $a,b\in\mathbb{R}^n$, the Euclidean inner product is implied $\langle a,b\rangle_2:=a^\mathsf{T}b$, unless otherwise specified by a subscript.

$L_n^p(\mathcal{X},\mathcal{F},\mu)$ is the Banach space of $f\colon\mathcal{X}\to\mathbb{R}^n$ functions endowed with the norm $\|f\|_{L_n^p}:=\left(\int_{\mathcal{X}}\|f(x)\|^p\,\mathrm{d}\mu(x)\right)^{\frac{1}{p}}$, where $(\mathcal{X},\mathcal{F},\mu)$ is a measure space. For $\mathcal{X}\subset\mathbb{R}^m$ the Lebesgue $\sigma$-algebra and measure are implied and the notation may be shortend to $L_n^p(\mathcal{X})$. Additionally, the subscript may be ommited for $n=1$, i.e., $L^p:=L_1^p$.

$L_n^2(\mathcal{X},\mathcal{F},\mu)$ is the Hilbert space of $f\colon\mathcal{X}\to\mathbb{R}^n$ functions endowed with the inner product $\langle f,g\rangle_{L_n^2}:=\int_{\mathcal{X}}f(x)^\mathsf{T}g(x)\,\mathrm{d}\mu(x)$, where $(\mathcal{X},\mathcal{F},\mu)$ is a measure space. For $\mathcal{X}\subset\mathbb{R}^m$ the Lebesgue $\sigma$-algebra and measure are implied and the notation may be shortend to $L_n^2(\mathcal{X})$. Additionally, the superscript may be ommited for $n=1$, i.e., $L^2:=L_1^2$.

$\mathcal{W}_n^2([a,b])$ is the Hilbert space of absolutely continuous $g\colon[a,b]\to\mathbb{R}^n$ functions with square-integrable weak derivatives $\dot{g}\in L_n^2([a,b])$, endowed with the inner product $\langle f,g\rangle_{\mathcal{W}_n^2}:=f(a)^\mathsf{T}g(a)+\int_a^b\dot{f}(t)^\mathsf{T}\dot{g}(t)\,\mathrm{d}t$, which is the direct sum of $n$ copies of the Sobolev space $\mathcal{W}^{2,1}([a,b])$.

$\mathcal{C}(\mathcal{X},\mathcal{Y})$ is the space of continuous functions $f\colon\mathcal{X}\to\mathcal{Y}$ between the topological spaces $\mathcal{X}$ and $\mathcal{Y}$. The domain and codomain can be ommited if they can be inferred from the context. If $\mathcal{X}$ is a compact Hausdorff space, like a closed interval of $\mathbb{R}$, and $\mathcal{Y}=\mathbb{R}^n$, then $\mathcal{C}$ is assumed to be endowed with the supremum norm $\|\!|\cdot|\!\|$, unless otherwise noted.

$\mathrm{PL}(\mathcal{P},\mathcal{Y})$ is the space of piecewise linear functions, with breaks over the partition $\mathcal{P}$, from $[\min(\mathcal{P}),\max(\mathcal{P})]$ to $\mathcal{Y}$.

$\partial\mathbb{A}$, $\partial\mathbb{c}$ is the boundary of a set $\mathbb{A}$ or the boundary of a $k-chain$ $\mathbb{c}$.

$\bar{\mathbb{A}}$ is the closure of a set $\mathbb{A}$.

$\operatorname{int}\mathbb{A}$ is the interior of a set $\mathbb{A}$.

$\int A\circ\mathrm{d}B$ is the Stratonovich integral of the process $A$ with respect to the process $B$.

## Miscellaneous

$\nabla_{\mathrm{x}}f(a):=\left[\frac{\partial f^{(i)}}{\partial x^{(j)}}(a)\right]^{(ij)}$ is the Jacobian matrix of the function $f$ with respect to the input $x$, evaluated at the point $a$.

$\operatorname{div}_{\mathrm{x}}f(a):=\sum_i\frac{\partial f(i)}{\partial x^{(i)}}(a)$ is the divergence of the function $f$ with respect to the input $x$, evaluated at the point $a$, i.e., the trace $\operatorname{tr}(\nabla_{\mathrm{x}}f(a))$ of its Jacobian matrix with respect to $x$.

# List of Acronyms

*BVP* boundary value problem,

*COIN-OR* computational infrastructure for operations research,

*EKF* extended Kalman filter,

*EKS* extended Kalman smoother,

*IPOPT* interior point optimizer,

*JMAPSPPE* joint maximum *a posteriori* state path and parameter estimator,

*IAE* integrated absolute error,

*ISE* integrated square error,

*MAP* maximum *a posteriori*,

*MEE* minimum energy estimator,

*MMSE* minimum mean square error,

*NLP* nonlinear program,

*ODE* ordinary differential equation,

*OEM* output error method,

*PEM* prediction error method,

*SDE* stochastic differential equation,

*UKF* unscented Kalman filter,

*UKS* unscented Kalman smoother,

*UFMG* Universidade Federal de Minas Gerais.

# Chapter 1

# Introduction

The subject of this thesis is joint maximum *a posteriori* state path and parameter estimation in systems described by stochastic differential equations, and its main contribution is the introduction of two new estimators. Besides their obvious use in joint state path parameter estimation, the proposed estimators can also be employed in system identification and in smoothing, i.e., in applications where both the parameters of the system and its states are unknown but either only the parameters or the states are sought after. Consequently, this thesis is concerned not only with the intersection of parameter and state estimation, but also with each of these fields of study independently.

We begin this chapter with an overview of smoothing and dynamical system parameter estimation, presented in Section 1.1, covering both their historical developments and the state of the art. We then proceed to present the motivation, objectives and contributions of this thesis in Section 1.2.

## 1.1 A brief survey of the related literature

In this section we present a brief survey of the literature related to the topics investigated herein. We begin by presenting important definitions and terminology which will be used throughout this thesis.

### 1.1.1 Classes of models and state estimation

To best estimate quantities associated with dynamical systems from measurements spread out over time, the measured values should be used in conjunction with knowledge of the system dynamics. This is especially important when both the measurements and the dynamics are uncertain, and described by a stochastic model instead of deterministic functions. Taking into account that for most real-world systems all models are approximations and all measurements have finite precision, models that take into account these shortcomings are more realistic than those that do not.

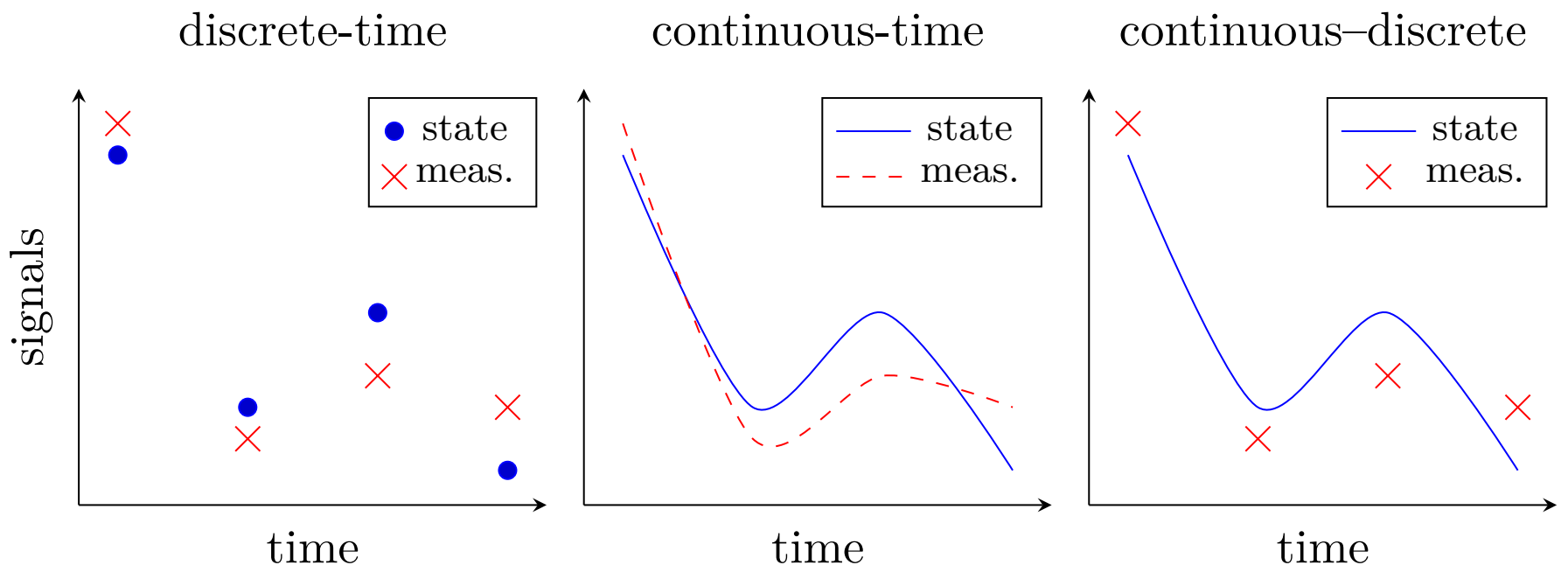

Figure 1.1: Graphical representation of the state and measurements for the three model classes.

The uncertainty in stochastic dynamical models can be due to unknown external disturbances acting on the system or can be simply a way for the modeler to express his lack of confidence and repeatability on the outcomes of the system. Examples of random external disturbances include electromagnetic interference in a circuit, thermal noise in a conductor, turbulence in an airplane and solar wind in a satellite, to name a few. Whether the behavior of the system—or of the corresponding disturbances—is truly random does not matter, stochastic models for dynamical systems are an important tool to make inference taking into account the uncertainties involved in the dynamics and data acquisition. These tools can be used with both the subjectivist and objectivist Bayesian points of view (for more information on these interpretations and their dichotomy see Press, 2003, Chap. 1).

Stochastic models for dynamical systems can be classified into three major classes, according to how the dynamics and measurements are represented (cf. Jazwinski, 1970, p. 144):

*discrete-time* models are those in which both the measurements and the underlying system dynamics are represented in discrete-time, over a possibly infinite but countable set of time points;

*continuous-time* models are those in which both the measurements and the underlying system dynamics are represented in continuous-time, over an uncountable set of time points;

*continuous–discrete* models, also known as sampled-data or hybrid models, are those in which the measurements are taken in discrete-time but the underlying system dynamics is represented in continuous-time. The measurements times are then a countable subset of the time interval over which the dynamics is represented.

These classes are illustrated graphically in Figure 1.1. Continuous–discrete models are particularly important as various phenomena of interest to science and engineering are of a continuous-time nature, yet the measurements available for inference are sampled at discrete time instants.

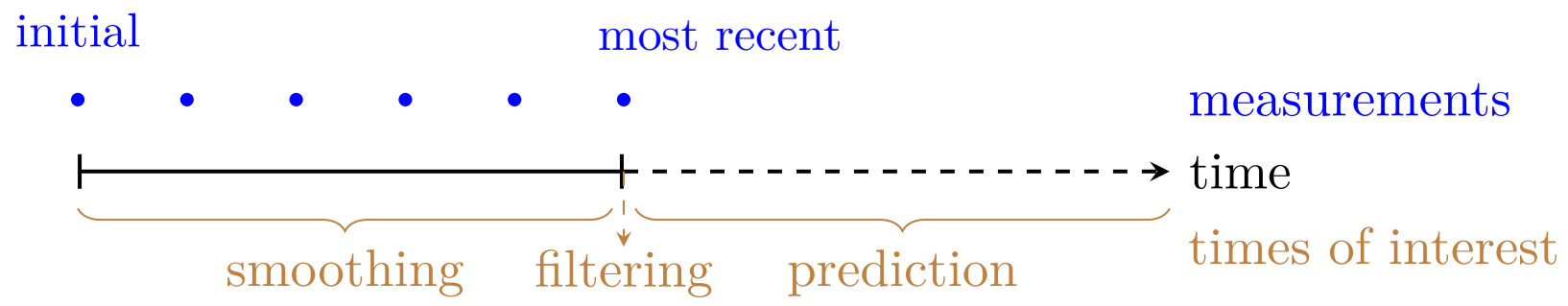

Figure 1.2: Graphical representation of the classes of state estimation.

The dynamics of discrete-time systems are usually represented by stochastic difference equations. Likewise, the dynamics of continuous-time and continuous–discrete systems are usually represented by stochastic differential equations (SDEs). Models represented by SDEs are used to make inference in a wide range of applications, including but not limited to radar tracking (Arasaratnam et al., 2010), aircraft flight path reconstruction (Mulder et al., 1999) and investment finance and option pricing (Black and Scholes, 1973); see the reviews by Kloeden and Platen (1992, Chap. 7), and Nielsen et al. (2000) for a comprehensive list.

In a very influential article, Kalman (1960) defined three related classes of estimation problems for signals under noise, given measurements spread out over time. His definitions are now standard terminology, and consist of the following classes:

*filtering,* where the goal is to estimate the signal at the time of the latest measurement;

*prediction,* also known as extrapolation or forecasting, where the goal is to estimate the signal at a time beyond that of the latest measurement;

*smoothing,* also known as interpolation, where the goal is to estimate the signal over times up to that of the latest measurement.

The filtering and prediction problem usually arise in online (real-time) applications, such as process control and supervision. The smoothing problem, on the other hand, arises in offline (*post-mortem*, after the fact) applications, or online applications where a delay between measurement and estimation is tolerable. In smoothing, future data (with respect to the estimate) is used to obtain better estimates (Jazwinski, 1970, p. 143). These classes are represented graphically in Figure 1.2.

Three additional subclasses of the smoothing problem were defined by Meditch (1967) and have also been adopted into the jargon:

*fixed-interval* smoothing, where the goal is to estimate the whole signal in between the earliest and latest measurement;

*fixed-point* smoothing, where the goal is to estimate the signal at a specific time-point before the latest measurement;

*fixed-lag* smoothing, where the goal is to estimate the signal at a specific time-distance from the latest measurement.

The difference between fixed-point and fixed-lag smoothing is only relevant in online applications.

In the context of Bayesian statistics, the smoothing estimates are usually chosen as the posterior mean or the posterior mode. The posterior mean is known as the minimum mean square error (MMSE), as it minimizes the expected square error with respect to the true states, over all outcomes for which the measured values would have been observed. The posterior mode is known as the maximum *a posteriori* (MAP) estimate, as it maximizes the posterior density. It is interpreted as the most probable outcome, given the measurements. Fixed-interval MAP smoothing can, furthermore, be divided into two classes:

*MAP state-path* smoothing, where the estimates are the joint posterior mode of the states along all time instants in the interval, given all the measurements available;

*MAP single instant* smoothing, where the estimate for each time instant is the marginal posterior mode of the state at the instant, given all the measurements available.

Godsill et al. (2001) refer to these classes as the joint and marginal MAP estimation. The MAP state-path smoothing is also referred to as MAP state-trajectory smoothing or, in the discrete-time case, MAP state-sequence smoothing. As there is not much literature on this class of estimators, the terminology is not well cemented.

### 1.1.2 Smoothing

The modern theory of smoothing has its origins in the works of Wiener (1964)[1] and Kolmogorov (1941)[2], who developed the solution to stationary linear systems subject to additive Gaussian noise. Wiener solved the continuous-time problem while Kolmogorov solved the discrete-time one. Their approach used the impulse response representation of systems to derive the results and represent the solutions. Wiener's work seems to have pioneered the use of statistics and probability theory to both formulate and solve the problem of filtering and smoothing of signals subject to noise.

While many built uppon the Wiener–Kolmogorov filtering theory—refer to Meditch (1973) for a comprehensive survey—the main breakthrough came in the seminal work of Kalman (1960), who solved the problem of non-stationary estimation in linear discrete-time and continuous–discrete systems subject to additive Gaussian noise. The Kalman filter differed significantly from the Wiener filter by the use of the state-space formulation in the time domain to derive the results and represent the solution. Shortly after, Kalman and Bucy (1961) proposed the Kalman–Bucy filter, the analog of the Kalman filter for continuous-time systems with continuous-time measurements. It should be noted that for linear–Gaussian systems there exist *exact* formulas for the discretization, so there is little difference between the discrete-time and continuous–discrete formulations.

[1] Originally circulated in 1942 as a classified memorandum, connected to the war effort.
[2] In Russian, for an English translation see Kolmogorov (1992).

In his initial paper, Kalman (1960) mentions that smoothing falls within the theory developed, but he does not address the smoothing problem directly. The Kalman filter was, however, readily extended to solve the smoothing problem by Bryson and Frazier (1963) for the continuous-time case, and by Cox (1964) and Rauch et al. (1965)[3] for the discrete-time case. The above mentioned smoothers work by combining the results of a standard forward Kalman filter with a backward smoothing pass, which led to them being classified as *forward–backward* smoothers. Another approach, introduced by Mayne (1966), consists of combining the results of two Kalman filters, one running forward in time and another one running backward, which came to be called *two-filter* smoothing.

The two-filter and forward backward smoothers can also be adapted to general nonlinear systems and to non-Gaussian systems. However, in these genearal cases, the smoothing distributions do not admit known closed-form expressions, except in some very specific cases. For these systems, consequently, practical smoothers must rely on approximations. One of the first of these approximations was the linearization of the system dynamics, which led to the development of the extended Kalman filter and smoothers (for a historical perspective on its development see Schmidt, 1981).

Nonlinear Kalman smoothers like the extended Kalman smoother follow roughly the same steps as their linear counterparts, summarizing the relevant distributions by their mean and variance. This implies that their underlying assumption is that the densities involved in the smoothing problem are approximately Gaussian. Consequently, the two-filter Rauch–Tung–Striebel or Bryson–Frazier formulas can be used if appropriate linear–Gaussian approximations are made. Besides the extended Kalman approach of linearization by truncating the Taylor series of the model functions, other approaches to obtain linear–Gaussian approximations include the unscented transform (Julier and Uhlmann, 1997), the statistical linear regression (Lefebvre et al., 2002) or Monte Carlo methods (Kotecha and Djuric, 2003), to name a few. The unscented Kalman filter of Julier and Uhlmann (1997), in particular, was formalized for smoothing with the two-filter formula by Wan and van der Merwe (2001, Sec. 7.3.2) and with the Rauch–Tung–Striebel forward–backward correction by Särkkä (2006, 2008).

[3] An early version of this paper was published in 1963 as a company report.

More general smoothers for nonlinear systems are based on non-Gaussian approximations of the smoothing distributions. The two-filter smoother of Kitagawa (1994), for example, uses Gaussian mixtures to represent the distributions and the Gaussian sum filters of Sorenson and Alspach (1971) to perform the estimation. Sequential Monte Carlo two-filter (Kitagawa, 1996; Klaas et al., 2006) and forward–backward smoothers (Doucet et al., 2000; Godsill et al., 2004) represent the distributions by weighted samples using the sequential importance resampling method of Gordon et al. (1993). Sequential Monte Carlo methods can also be used for MAP state-path smoothing by employing a forward filter in conjunction with the Viterbi algorithm to select the most probable particle sequence, as proposed by Godsill et al. (2001).

MAP state-path smoothers, nevertheless, do not need to be formulated or implemented in the sequential Bayesian estimator framework. For a wide variety of systems, the posterior state-path probability density admits tractable closed-form expressions[4] which can be easily evaluated on digital computers. Therefore, by the use of nonlinear optimization tools, the MAP state-path can be obtained without resorting to Bayesian filters.

Since the early developments in Kalman smoothing, MAP state-path estimation by means of nonlinear optimization was considered as a generalization of linear smoothing to nonlinear systems subject to Gaussian noise (Bryson and Frazier, 1963; Cox, 1964; Meditch, 1970). However, the limitation of the computers available at the time meant that the estimates could not be obtained by solving a single large-scale nonlinear program. Therefore, the division of the large problem into several small problems, as performed by the extended Kalman smoother, was advantageous. Furthermore, for discrete-time systems, the extended Kalman smoother was shown (Bell, 1994; Johnston and Krishnamurthy, 2001) to be equivalent to Gauss–Newton steps in the maximization of the posterior state-path density. Because of those similarities, the terminology is still not standard, with many different methods falling under the umbrella of Kalman smoothers. In this thesis, we will refer to as *nonlinear Kalman smoothers* those methods which rely on two-filter or forward–backward formulas, and as *MAP state-path smoothers* those methods which rely on nonlinear optimization to obtain their estimates.

As the tools of nonlinear optimization matured, optimization-based MAP state-path smoothers emerged again as powerful alternatives to nonlinear Kalman smoothers (Aravkin et al., 2011, 2012a,b,c, 2013; Bell et al., 2009; Dutra et al., 2012, 2014; Farahmand et al., 2011; Monin, 2013). Furthermore, unlike nonlinear Kalman smoothers, these smoothers are not limited to sys-tems with approximatelly Gaussian distributions, which lends them a wider applicability. It also opens the possibility to model the systems with non-Gaussian noise. Heavy-tailed distributions are often used to add robustness to the estimator against outliers, even when the measurements actually come from a different distribution. Examples of robust MAP state-path smoothing with heavy-tailed distributions include the work of Aravkin et al. (2011) and Farahmand et al. (2011) with the $\ell_1$-Laplace, and the work of Aravkin et al. (2012a,c) with piecewise linear-quadratic log densities and Student's $t$. Distributions with limited support, such as the truncated normal, can also be used to model constraints (Bell et al., 2009).

[4]Up to a constant factor which does not influence the location of maxima.

One theoretical issue with MAP state-path smoothing that persisted up to recently was its correct definition, interpretation, and formulation for continuous-time or continuous–discrete systems, and the applicability of discrete-time methods to discretized continuous systems. The definition is not trivial because the state-paths of these systems lie in infinite-dimensional separable Banach spaces, for which there is no analog of the Lebesgue measure, i.e., a locally finite, strictly positive and translation-invariant measure. Consequently, the measure-theoretic definition of density as the Radon–Nikodym derivative of the probability measure is not applicable for the purposes of obtaining the posterior mode (Dürr and Bach, 1978, p. 155).

Nevertheless, systems described by stochastic differential equations do possess asymptotic probability of small balls, given by the Onsager–Machlup functional. This functional, first derived by Stratonovich (1971) for nonlinear SDEs, can be thought of as an analog of the probability density for the purpose of defining the mode. In the statistics and physics literature, a consensus arose that the maxima of the Onsager–Machlup functional represented the most probable paths of diffusion processes (Dürr and Bach, 1978; Ito, 1978; Takahashi and Watanabe, 1981; Ikeda and Watanabe, 1981, Sec. VI.9; Hara and Takahashi, 1996) and that it should be used for construction of MAP state-path estimators (Zeitouni and Dembo, 1987).

Despite that, in the automatic control and signal processing literature, the Onsager–Machlup functional seemed to be virtually unknown until the late 20th century. One of the few uses of the Onsager–Machlup functional for MAP state-path estimation in the engineering literature was the work of Aihara and Bagchi (1999a,b), which did not enjoy much popularity.[5] In fact, as is detailed in the next paragraph, many authors incorrectly obtained the energy functional as the merit function for MAP state-path estimation. The energy functional differs from the Onsager–Machlup functional by the lack of a correction term which, when ommitted, can lead to different estimates.

Cox (1963, Chap. III), for example, obtained the energy functional by applying the probability density functional of Parzen (1963) to the driving noise of the SDE, without taking into account that the nonlinear drift might cause the mode of the state path to differ from the mode of the noise path, as proved to be the case in Chapter 2 of this thesis and in one of its derivative works (Dutra et al., 2014). Mortensen (1968) obtained the energy functional by Feynman path integration in function spaces, but apparently used an incorrect integrand, as the Onsager–Machlup function is the correct Lagrangian in the path integral formulation of diffusion processes (Graham, 1977).[6] Jazwinski (1970, p. 155), on the other hand, obtained the energy by considering the limit of the mode of the Euler-discretized systems. However, as argued by Horsthemke and Bach (1975), the Euler discretization cannot be used for such purposes due to insufficent approximation order.

[5] Only one citation from other authors, up to 2013, was found for both these articles in the SciVerse Scopus, Web of Science, and Google Scholar databases.

The recent renewed interest in MAP state-path estimators, afforded by the maturity of the necessary optimization tools, focused mostly on discrete-time systems. When applied to continuous-discrete systems described by SDEs, discretization schemes were used, as done by Bell et al. (2009) and Aravkin et al. (2011, 2012c). For nonlinear systems, SDE discretization schemes are approximations which improve as the discretization step decreases. Applications in Kalman filtering that use such discretization schemes, for example, usually require fine discretizations to produce meaningful results (Arasaratnam et al., 2010; Särkkä and Solin, 2012). An important requirement for the adoption of discrete-time MAP state-path smoothing for discretized continuous-time systems is then understanding if and how the discretized MAP state path relates to the MAP state path of the original system.

In this thesis we proved that, under some regularity conditions, discretized MAP state-path estimation converges hypographically to variational estimation problems, as the discretization step vanishes. However, the limit depends on the discretization scheme used. For the Euler discretization scheme, the most popular, the discretization converges to the minimum energy estimation, which is proved to correspond to the state path associated with the MAP *noise path*. When the trapezoidal scheme is used, on the other hand, the estimation converges to the MAP *state-path* estimation using the Onsager–Machlup functional. Our work also offers a formal definition of mode and MAP estimation of general random variables in possibly infinite-dimensional spaces, under the framework of Bayesian decision theory. Finally, it highlights the links between continuous–discrete MAP state-path estimation and optimal control, opening up an even wider range of tools for this class of estimation problems. These results were published in a derivative work of this thesis (Dutra et al., 2014).

Having layed down a solid theoretical foundation for MAP state-path estimation in systems described by SDEs, the field is now ripe for its use in non-Gaussian applications, such as robust smoothing, and fields typically dominated by nonlinear Kalman smoothing, such as aircraft flight path reconstruction (Mulder et al., 1999; Jategaonkar, 2006, Chap. 10; Teixeira et al., 2011).

[6]The exact integrand used could not be found, as the derivation was performed in "mimeographed lecture notes" which could not be recovered.

### 1.1.3 System identification

System identification is the inverse problem of obtaining a description for the behavior of a system based on observations of its behavior. For dynamical systems, the system is represented by a suitable mathematical model and the observations are measured input–output signals of its operation. The process of dynamical-system identification is often structured as composed of four main tasks:

*data collection,* the process of planning and executing tests to obtain representative observations of the system dynamics;

*model structure determination,* the process of choosing a parametrized structure and mathematical representation for the model, using both *a priori* information and some of the collected data;

*parameter estimation,* the process of obtaining parameter values from the collected data and prior knowledge;

*model validation,* the process of evaluating if the chosen model with its estimated parameters is an adequate description of the system.

Each of these tasks is dependent upon the preceding ones, but the overall process can be iterative and involve multiple repetitions until a suitable combination of data, model, and parameters is found.

With respect to the model structure, three main classes are usually considered:

*white-box* models, also known as phenomenological models, constructed from first principles and knowledge of the system internals;

*black-box* models, also known as behavioral models, of a general structure chosen to represent the input–output relationship without regard to the system internals;

*grey-box* models, those which combine aspects of both white- and black-box models, featuring components whose dynamics are derived from first principles and others whose dynamics only represent some cause–effect relationship.

Grey- and white-box models have the advantage that their parameters can have physical meaning and can be related to other sources. Often, discrete-time models are used in black-box modelling, and continuous-time and continuous–discrete models are used in grey- and white-box modelling, as phenomenological descriptions of the dynamics of most systems are formulated in continuous-time. Additionally, the representation of the state dynamics by stochastic differentential equations is a valuable tool for grey-box modelling, as it can put together a simplified phenomenological representation with a source of uncertainty that accounts for the simplifications performed (cf. Kristensen et al., 2004a,b).

### 1.1.4 Parameter estimation in stochastic differential equations

In the process of system identification, parameter estimation can be thought of as procedures for extracting the information from the data. For systems for which the states are directly observable, without measurement noise, the field is well-consolidated and there exist a plethora of estimation techniques; see Nielsen et al. (2000, Secs. 3–4) and Bishwal (2008) for reviews. However, while the assumption of no measurement noise might be realistic for stock market and finantial modelling, it is overly optimistic for most engineering applications. When measurement noise is considered, the estimation problem becomes more complex and the field is still under active development.

In the remainder of this section, we will focus only on methods and techniques which take measurement noise into account. As the states are not directly measured, this class of parameter estimators is closely related to state estimators. Intuitively thinking, to model the evolution of the states one first needs estimates of their values.

For Markovian systems, the likelihood function can be decomposed into the product of the one-step-ahead predicted output densities, given all previous measurements, evaluated at the measured values. Consequently, for linear systems subject to Gaussian noise, the likelihood function can be obtained by employing a Kalman filter to obtain the predicted output distributions. This approach was pioneered by Mehra (1971) and came to be known as the prediction-error method or filter-error method. Among its first applications was the parameter estimation of aircraft models from flight-test data collected under turbulence (Mehra et al., 1974). These methods can be used for both maximum likelihood and maximum *a posteriori* parameter estimation; see the survey by Åström (1980) for a review of its early developments.

Just like the Kalman filter was adapted to nonlinear systems by linearization of the model functions, the prediction-error method was readily adapted to nonlinear systems by employing the extended Kalman filter to obtain the predicted output densities (Mehra et al., 1974). Like its underlying nonlinear Kalman filter, this approach is based on the the implicit assumption that the the state and output distributions involved are all approximately conditionally Gaussian, given the parameters. Similarly, any other nonlinear filtering approach, such as the unscented Kalman filter (Chow et al., 2007) or particle filters (Andrieu et al., 2004; Doucet and Tadić, 2003), can be used in prediction error methods. The limitations and difficulties associated with these methods are, likewise, those of their underlying nonlinear filters. Furthermore, the performance of the filters is usually significantly degraded when the parameters are far from the optimum, which can make the estimator sensitive to the initial guess.

Another popular approach for parameter estimation worth mentioning consists of treating parameters as augmented states and using standard state estimation techniques (see Jazwinski, 1970, Sec. 8.4; or Ljung, 1979, and the references therein). These methods have the serious disadvantage, however, that if no process noise is assumed to act on the augmented states corresponding to the parameters, their estimated variance will converge to zero. Given all the approximations involved, that is overly optimistic and, furthermore, can cause the estimators to fail. Consequently, a small artificial noise is assumed to act on the dynamics of, but this workaround also introduces its own drawbacks. The choice of the "parameter noise" variance is difficult and its effect on the original problem is not always clear. Also, the final output of the method is not a single estimate for the parameters but a time-varying one. It is not obviuos how to choose one instant or combine the estimate a different instants to obtain the best value.

Up to the early 21st century, not much new development was done in parameter estimation in systems described by SDEs subject to measurement noise, as argued by Kristensen et al. (2004b, p. 225). A review by Nielsen et al. (2000) pointed to the decades-old prediction error method with the extended Kalman filter as the "most general and useful approach" at the time for parameter estimation in this class of systems.

A new development came in the work of Varziri et al. (2008b), which they termed the approximate maximum likelihood estimator. As explained in (Varziri et al., 2008a, p. 830), their estimator can be interpreted as a joint MAP state-path and parameter estimator with a non-informative parameter prior. If the state-path and parameter posterior is then approximately jointly symmetric, the estimated parameter should be close to the true maximum likelihood estimate. However, it should be noted that they used the energy functional, as derived by Jazwinski (1970, p. 155), instead of the Onsager–Machlup functional to construct their merit function. This means that their estimates have a different interpretation than intended, as is proved in Chapter 2 of this thesis. In Chapter 3 of this thesis we prove, furthermore that their derivation would be the one intended if the trapezoidal discretization scheme had been used instead of the Euler scheme.

One limitation of Varziri's approximate maximum likelihood estimator was that the measurement and process noise variances had to be known *a priori*. That limitation was overcome first by using iterative procedures based on heuristics (Karimi and McAuley, 2014b; Varziri et al., 2008c). However, a more rigorous and promising approach is using the Laplace approximation to marginalize the states in the joint state-path and parameter density to obtain an approximation of the parameter likelihood function (Karimi and McAuley, 2013, 2014a; Varziri et al., 2008a). In this context, the Laplace approximation is a technique to marginalize the state-path out of the posterior joint state-path and paremeter density, obtaining the posterior parameter density. It consists of making a Gaussian approximation of the state-path at its mode from the Hessian of the log-density. It should be noted that the assumption that the posterior state-path smoothing distribution is approximately Gaussian is usually less strict than the assumption that the one-step-ahead predicted output is approximately Gaussian.

Once again, these estimators used the energy functional instead of the Onsager–Machlup functional, implying that they are, in fact, applying the Laplace approximation to marginalize the noise path out of the joint noise-path and parameter distribution. The noise-path is usually less observable than the state-path, implying that the Hessian of its log-posterior usually has a smaller norm and the Laplace approximation is coarser. We also remark that, although easily adapted to non-Gaussian measurement and prior distributions, the literature on these estimators is restricted to the case where these distributions are all Gaussian (Karimi and McAuley, 2013, 2014a,b; Varziri et al., 2008a,b; Varziri et al., 2008c).

We can then see that joint MAP state-path and parameter estimators figure as an important tool for further development in parameter estimators for systems described by stochastic differential equations. However, theoretical difficulties with the correct definition of MAP state-path and the effects of discretizations must be clearly resolved for these methods to gain a wider adoption. Once these issues are resolved, these estimators can be used in applications dominated by nonlinear-Kalman prediction-error methods and non-Gaussian applications. In particular, just like with MAP state-path estimation, heavy-tailed distributions can be used to perform robust system identification in the presence of outliers.

## 1.2 Purpose and contributions of this thesis

In this section, we lay down the purpose and contributions of this thesis to its related fields of study. We begin with its theoretical and technological motivations and proceed to the objectives of this work and its contributions to the literature.

### 1.2.1 Motivations

The motivations of this thesis are twofold: a technological driver based on demand from applications; and the need to better understand the theory behind maximum *a posteriori* state-path and joint state-path and parameter estimation in systems described by stochastic differential equations. Inportant theoretical considerations that need a better understanding are the probabilistic interpretations of the minimum energy and MAP joint state-path and parameter estimators and their relationship to the discretized estimators.

The *Universidade Federal de Minas Gerais*[7] (UFMG) is one of the leading centers of aeronautical technology research and education in Brazil. Its activities encompass all aspects of aeronautical engineering, including desing, construction and operation of airplanes and aeronautical systems. Among its current projects are a light airplane with assisted piloting, a racing airplane to beat world speed records, and hand-lauched autonomous unmanned aerial vehicles.

Flight testing is an integral part of design and operation of aircraft and aeronautical systems. Its purposes include analysis of aircraft performance and behaviour; identification of dynamical models for control and simulation; and testing of systems and algorithms. Flight test data is subject to noise and sensor imperfections; it is typically preprocessed before any analysis is done. The process of recovering the flight path from test data is known as flight-path reconstruction (Mulder et al., 1999; Jategaonkar, 2006, Chap. 10; Klein and Morelli, 2006, Chap. 10) and is an essential step for flight-test data analysis. In adition, flight vehicle system identification is an indispensable tool in aeronautics with great engineering utility; see the editorial by Jategaonkar (2005) from the special issues of the Journal of Aircraft on this topic and some of the reviews therein (Jategaonkar et al., 2004; Morelli and Klein, 2005; Wang and Iliff, 2004) for more information.

Small hand-launched unmanned aerial vehicles, like the ones being designed and operated at UFMG, present additional challenges to flight-path reconstruction and system identification. Their small payload limits the weight of the flight-test instrumentation that can be carried onboard, restricting the instrumentation to lightweight sensors with more noise, bias and imperfections. Their small weight also makes them more sensitive to turbulence, which needs to be accounted for. The technological motivation of this thesis is then to develop state and parameter estimators for flight-path reconstruction and system identification which are appropriate for use in small hand-launched unmanned aircraft.

Related to our technological motivation, we have that the tools of MAP state-path estimation, which show great promise for our intended application, still has an unclear definition for systems described by stochastic differential equations. Furthermore, the relationship between discretized estimators and the continuous-time underlying problem needs to be understood with a rigorous mathematical analysis. Consequently, our theoretical motivations are to provide a firm theoretical basis for MAP state-path and parameter estimation in continuous-time and continuous–discrete systems.

[7]Federal University of Minas Gerais

### 1.2.2 Objectives and contributions

Having stated the motivations behind this work, the objectives of this thesis are then

- to provide a rigorous definition of mode and MAP estimation for random variables in possibly infinite-dimensional spaces which coincides with the usual definition for continuous and discrete random variables;
- to derive the joint MAP state-path and parameter estimator for systems described by SDEs, together with conditions and assumptions for its applicability;
- to obtain a probabilistic interpretation for the joint minimum energy state-path and parameter estimator, together with conditions and assumptions for its applicability;
- to relate discretized joint MAP state-path and parameter estimators to their continuous-time counterparts;
- to demonstrate the viability of the derived estimators in example applications with both simulated and experimental data.

We believe these are important steps for further advancement of the field and the use of these techniques in novel and challenging applications.

The specific contributions and novelties of this thesis are

- formalization of the concept of fictitious densities, in Definition 2.4;
- formalization of the concept of mode and MAP estimator for random variables in metric spaces, using the concept of fictitious densities, in Proposition 2.6 and Definitions 2.5 and 2.7;
- derivation of the Onsager–Machlup functional for systems with unknown parameters and possibly singular diffusion matrix, in Theorem 2.9;
- probabilistic interpretation of the energy functional as the fictitious density of the associated noise path, in Theorem 2.26, extending the results published in (Dutra et al., 2014) to systems with unknown parameters and possibly singular diffusion matrices;
- derivation of the hypographical limits of the Euler- and trapezoidally-discretized joint state path and parameter densities, in Theorems 3.7 and 3.14, once more extending the results published in (Dutra et al., 2014) to systems with unknown parameters and possibly singular diffusion matrices.

## 1.3 Outline of this thesis

The remainder of this thesis is organized as follows. In Chapter 2 we formally define the MAP estimator for general random variables and derive it for joint state-path and parameter estimation in systems described by SDEs. In addition, we define the minimum energy estimator for state-path and parameter estimation in SDEs and prove that it is associated with the MAP noise-paths and parameters. In Chapter 3, we briefly present the concepts of hypographical convergence and show how it can be applied to the discretized MAP estimation. We then derive the Euler- and trapezoidally-discretized joint MAP state-path and parameter estimators and obtain their hypographical limits. The proposed estimators are then demonstrated in Chapter 4 with both simulated and experimental data and the conclusions and future directions for extension of this work presented in Chapter 5.

# Chapter 2

# MAP estimation in SDEs

Don't panic!

Douglas Adams, *The Hitchhiker's Guide to the Galaxy*

This chapter contains the main theoretical contributions of this thesis. We begin with the definition and interpretations of maximum *a posteriori* (MAP) estimation in Section 2.1, and then proceed to apply the concepts developed for the derivation of the joint MAP state-path and parameter estimator for systems described by stochastic differential equations (SDEs) in Section 2.2. Then, in Section 2.3, the joint MAP noise-path and parameter estimator is derived and shown to be the minimum energy estimator obtained by omitting the drift divergence in the Onsager–Machlup functional.

## 2.1 Foundations of MAP estimation

In this section, we present the theoretical foundations of maximum *a posteriori* (MAP) estimation. We begin with a brief presentation of Bayesian point estimation and show how some parameters of the posterior distribution, such as the mean and mode, are the Bayesian estimates associated with popular loss functions. We then define the mode of discrete and continuous random variables over $\mathbb{R}^n$ and extend this definition to random variables over general metric spaces, like paths of stochastic processes. The MAP estimator is then defined as the posterior mode and interpreted in the context of Bayesian decision theory.

### 2.1.1 General Bayesian estimation

In Bayesian statistics, the posterior distribution of a random quantity of interest, given an observed event, is the aggregate of all the information available on the quantity (Migon and Gamerman, 1999, p. 79). This information is represented in the form of a probabilistic description, and includes the prior, what was known before the observation, and the likelihood, what is added by the observation. This information can be used to make optimal decisions in the face of uncertainty, with the application of Bayesian decision theory.

In Bayesian decision theory, a loss function is chosen to represent the undesirability of each choice and random outcome. The Bayesian choice is then that which minimizes the expected posterior loss, over all possible decisions (Robert, 2001). Alternatively, the problem can be modeled in terms of the gain associated with each choice and random outcome, in which case the Bayesian choice is the one which maximizes the expected posterior gain. Gain functions are also refered to as utility functions.

One application of Bayesian decision theory is Bayesian point estimation. Although any summarization of the posterior distribution leads to some loss of information, it is often necessary to reduce the distribution of a random variable to a single value for some reason, among them reporting, communication or further analysis which warrants a single, *most representative*, value. A loss function is then used to quantify the suitability of each estimate, for each possible outcome of the random variable of interest, and the point estimate is the choice which minimizes the expected posterior loss. This concept is detailed below.

Let $X$ be an $\mathcal{X}$-valued random variable of interest and $Y$ be a $\mathcal{Y}$-valued observed random variable. An estimator is a deterministic rule $\delta\colon \mathcal{Y} \to \mathcal{X}$ to obtain estimates for $X$ given observed values of $Y$. The performance of estimators is compared using measurable loss functions $L\colon \mathcal{X} \times \mathcal{X} \to \mathbb{R}_{\geq 0}$. For each outcome $X = x$, $L(\hat{x}, x)$ evaluates the penalty (or error) associated with estimate $\hat{x}$. The integrated risk (or Bayes risk) associated with each estimator is then defined by

$$r(\delta) := \mathrm{E}\big[L(\delta(Y), X)\big]\,,$$

where $\mathrm{E}[\cdot]$ denotes the expectation operator. The risk is the mean loss over all possible values of the variable of interest and the observation. Having chosen an appropriate loss function for the problem at hand, the Bayesian estimators are defined as follows.

**Definition 2.1** (Bayesian estimator). *A Bayesian estimator $\delta^{\mathrm{b}}: \mathcal{Y} \to \mathcal{X}$ for the random variable $X$, associated with the loss function $L$, given the observed value $y \in \mathcal{Y}$ of $Y$, is any which minimizes the expected posterior loss, i.e.,*

$$\mathrm{E}\big[L(\delta^{\mathrm{b}}(y), X) \,\big|\, Y = y\big] = \min_{x \in \mathcal{X}} \mathrm{E}\big[L(x, X) \,|\, Y = y\big]\,.$$

Note that the minimum—and hence the Bayesian estimator—is not guaranteed to exist and might not be unique. However, when it exists, no other estimator is better according to the criterion defined by $L$, the probabilistic model in $P$ and the observation $y$. Multiplicity of minima means that there exist several choices which are indistinguishable with respect to that criterion. In addition, if a Bayesian estimator exists for $P_Y$-almost everywhere $y$, then it attains the lowest possible integrated risk (Robert, 2001, Thm. 2.3.2).

One way to interpret these quantities is making an analogy to a game of chance. The estimate $\hat{x}$ is a bet on the value of $X$, given that the player knows $Y = y$ about the state of the game. The estimator represents a betting strategy and the loss function the financial losses for each outcome and bet, that is, the payoff. The integrated risk is then the average loss expected for a given betting strategy. Thus, an optimal betting strategy would minimize the integrated risk, leading to the minimum financial losses in the long run.[1]

There exist many canonical loss functions whose Bayesian estimators are certain statistics of the posterior distribution. Take, for example, quadratic error losses:

$$L_2(\hat{x}, x) := h(x - \hat{x}, x - \hat{x}),$$

where $h\colon \mathcal{X} \times \mathcal{X} \to \mathbb{R}_{\geq 0}$ is a coercive and bounded bilinear functional. For $\mathcal{X} = \mathbb{R}^n$ then $h$ is a positive-definite quadratic form, i.e.,

$$L_2(\hat{x}, x) = (x - \hat{x})^\mathsf{T} Q (x - \hat{x})$$

for any positive definite $Q \in \mathbb{R}^{n \times n}$. The Bayesian estimator associated with this family of loss functions is the posterior mean of $X$. This result was proved for $\mathbb{R}$-valued random variables by Gauss and Legendre (see Robert, 2001, Sec. 2.5.1) but also holds, under some regularity assumptions, when $\mathcal{X}$ is a more general reflexive Banach space.

A similar result holds for the absolute error loss,

$$L_1(\hat{x}, x) := \|\hat{x} - x\| .$$

When $X$ is an absolutely integrable $\mathbb{R}$-valued random variable, then the Bayesian estimator associated with the $L_1$ loss is is the posterior median, a result initially proved by Laplace (see Robert, 2001, Prop. 2.5.5). When $X$ is $\mathbb{R}^n$-valued and $\mathrm{E}[|X| \mid Y = y] < \infty$, then the posterior spatial median (also known as the $L_1$-median) is the Bayesian estimator associated with the $L_1$ loss. This concept also generalizes well to infinite-dimensional–valued $X$; the spatial median is also the Bayesian estimator associated with the $L_1$ loss when $\mathcal{X}$ is a reflexive Banach space (Averous and Meste, 1997).

The modes, which are known as maximum *a posteriori* estimates, are only *degenerate* Bayesian estimates, however, for random variables over general spaces. Nevertheless, they can be interpreted in the framework of Bayesian decision theory and Bayesian point estimation. Before showing how MAP estimation fits into this framework, we present in the following subsection the formal definition of modes for random variables over possibly infinite-dimensional spaces.

[1]The requirement of $L \geq 0$ would not attract many gamblers to the game, however, since it means that no bet and outcome combination yields profit. In this case the best strategy would actually be not entering the game. Perhaps this is a critique of gambling by the Bayesian decision theorists.

### 2.1.2 Modes of random variables

The modes of a random variable are population parameters that correspond to the region of the variable's sample space where the probability is most concentrated. For discrete random variables, i.e., when there is a countable subset of the variable's sample space with probability one, then the modes are its most probable outcomes. For continuous random variables over $\mathbb{R}^n$, i.e., those whose probability measure admits a density with respect to the Lebesgue measure, all individual outcomes of the random variable have probability zero, so the modes are defined as the densest outcomes (Prokhorov, 2002). These definitions are formalized below.

**Definition 2.2** (modes of discrete random variables)**.** *Let $X$ be an $\mathcal{X}$-valued random variable such that there exists a countable set $\mathbb{A} \subset \mathcal{X}$ with $P_X(\mathbb{A}) = 1$. The* modes *of $X$ are the points $\hat{x} \in \mathcal{X}$ satisfying*

$$P(X = \hat{x}) = \max_{x \in \mathcal{X}} P(X = x),$$

*at least one which exists.*

**Definition 2.3** (modes of continuous random variables)**.** *Let $X$ be an $\mathbb{R}^n$-valued random variable that admits a continuous density $f\colon \mathbb{R}^n \to \mathbb{R}_{\geq 0}$ with respect to the Lebesgue measure. The* modes *of $X$, if they exist, are the points $\hat{x} \in \mathbb{R}^n$ satisfying*

$$f(\hat{x}) = \max_{x \in \mathbb{R}^n} f(x).$$

We note that Definition 2.3 is restricted to random variables with continuous densities because otherwise the definition is ambiguous. Any function which is equal Lebesgue-almost everywhere to a probability density function is also a density for the same random variable. If the random variable admits a continuous density, then the continuous one is considered as the best density, as it is the one with the largest Lebesgue set and coincides everywhere with the Lebesgue derivative of the induced probability measure (Stein and Rami, 2005, p. 104). If the random variable only admits discontinuous densities, however, then it is not clear which density is the best and should be used for calculating the mode.

When generalizing the above definitions for paths of diffusion processes, probability densities (in the measure-theoretic sense) cannot be used. That is because for infinite-dimensional separable Banach spaces, like the space of paths of $\mathbb{R}^n$-valued diffusions, there is no analog of the Lebesgue measure, i.e., a translation-invariant, locally finite and strictly positive measure. Any translation-invariant measure on an infinite-dimensional separable Banach space would assign either infinite or zero measure to *all* open sets. Thus any measure with respect to which a density (Radon–Nikodym derivative) is taken would weight *equally-sized*[2] neighbourhoods differently. Similar problems arise in more general metric spaces, like the space of paths of diffusion processes over manifolds (cf. Takahashi and Watanabe, 1981).

This limitation is overcome by using a functional that quantifies the concentration of probability in the neighborhood of paths but, unlike the density, cannot be used to recover the probability of events via integration. The concentration is quantified by the asymptotic probability of metric $\epsilon$-balls, as $\epsilon$ vanishes. The normalization factor of this asymptotic probability weights balls of equal radius equally. Stratonovich (1971), who first proposed these ideas, called it the probability density functional of paths of diffusion processes. We believe, however, that the nomenclature of Takahashi and Watanabe (1981) is more appropriate: the functional is "an *ideal* density with respect to a *fictitious* uniform measure"[3]. Zeitouni (1989) shortened this nomenclature to *fictitious density*, which is the term that is used in this thesis. This concept is formalized below.

**Definition 2.4** (fictitious density)**.** *Let $X$ be an $\mathcal{X}$-valued random variable, where $(\mathcal{X}, d)$ is a metric space. The function $f\colon \mathcal{X} \to \mathbb{R}_{\geq 0}$ is an ideal density with respect to a fictitious uniform measure, or a* fictitious density *for short, if there exists $g\colon \mathbb{R}_{>0} \to \mathbb{R}_{>0}$ such that*

$$\lim_{\epsilon \downarrow 0} \frac{P(d(X,x) < \epsilon)}{g(\epsilon)} = f(x) \qquad \text{for all } x \in \mathcal{X} \tag{2.1}$$

*and $f(x) > 0$ for at least one $x \in \mathcal{X}$.*

Both probability mass functions and continuous probability density functions are fictitious densities, according to Definition 2.4. As the modes of discrete and continuous random variables are the location of the maxima of the fictitious density, it is natural to define the mode of general random variables that admit a fictitious density as the location of the maxima of the fictitious densities.

**Definition 2.5** (mode in metric spaces)**.** *Let $X$ be an $\mathcal{X}$-valued random variable, where $(\mathcal{X}, d)$ is a metric space. If $X$ admits a fictitious density $f$, its mode, if it exists, is any $\hat{x} \in \mathcal{X}$ that satisfies*

$$f(\hat{x}) = \max_{x \in \mathcal{X}} f(x).$$

[2]Such as metric balls of the same radius or translations of the same set.

[3]Emphasis in the original.

An alternative way to define the mode, which can help with the interpretation of the its statistical meaning, is the definition proposed in one of the derivative works of this thesis (Dutra et al., 2014, Defn. 1) and presented below.

**Proposition 2.6** (alternative definition of mode)**.** *Let $X$ be an $\mathcal{X}$-valued random variable, where $(\mathcal{X}, d)$ is a metric space. Any $\hat{x} \in \mathcal{X}$ is a mode of $X$, according to Definition 2.5, if and only if*

$$\lim_{\epsilon \downarrow 0} \frac{P\big(d(X,x) < \epsilon\big)}{P\big(d(X,\hat{x}) < \epsilon\big)} \leq 1 \qquad \text{for all } x \in \mathcal{X}. \tag{2.2}$$

*Proof.* If $\hat{x}$ is the mode of $X$ according to Definition 2.5, $f$ is a fictitious density and $g(\epsilon)$ is the fictitious uniform measure of an $\epsilon$-ball, then

$$\lim_{\epsilon \downarrow 0} \frac{P\big(d(X,x) < \epsilon\big)}{P\big(d(X,\hat{x}) < \epsilon\big)} = \frac{f(x)}{\max_{x' \in \mathcal{X}} f(x')} \leq 1 \qquad \text{for all } x \in \mathcal{X}.$$

conversely, if (2.2) is satisfied for some $\hat{x} \in \mathcal{X}$, then $X$ admits a fictitious density using $g(\epsilon) := P\big(d(X,\hat{x}) < \epsilon\big)$ as a fictitious uniform measure of $\epsilon$-balls. Futhermore, the maximum of the fictitious density, which is equal to one, is attained at $\hat{x}$, implying that it is a mode. □

Proposition 2.6 means that, if $\hat{x}$ is a mode and $x$ is not, then there exists an $\bar{\epsilon} \in \mathbb{R}_{>0}$ such that, for all $\epsilon \in (0, \bar{\epsilon}]$ the $\hat{x}$-centered $\epsilon$-ball has higher probability than the $x$-centered one. Additionally, if both $\hat{x}$ and $\tilde{x}$ are modes, then for all $\delta \in \mathbb{R}_{>0}$ there exists an $\bar{\epsilon} \in \mathbb{R}_{>0}$ such that, for all $\epsilon \in (0, \bar{\epsilon}]$,

$$(1-\delta)P\big(d(X,\tilde{x}) < \epsilon\big) < P\big(d(X,\hat{x}) < \epsilon\big) < (1+\delta)P\big(d(X,\tilde{x}) < \epsilon\big)\,,$$

that is, the probabilities of $\epsilon$-balls centered on both values is arbitrarily close.

### 2.1.3 Maximum *a posteriori* and Bayesian estimation

A maximum *a posteriori* (MAP) estimate of a random variable, given an observation, consists of the mode of its posterior distribution. From the definitions of the preceding section (Sec. 2.1.2), we have that it can be interpreted as the point of the variable's sample space around which the *posterior* probability is most concentrated. The MAP estimator is the Bayesian equivalent of the maximum likelihood estimator (Robert, 2001, Sec. 4.1.2; Migon and Gamerman, 1999, p. 83). It differs from the latter by the inclusion of prior information, and can also be interpreted as a penalized maximum likelihood estimator of classical statistics. Its definition is formalized below.

**Definition 2.7** (maximum *a posteriori* estimator)**.** *Let $X$ be an $\mathcal{X}$-valued random variable, where $(\mathcal{X}, d)$ is a metric space. A maximum* a posteriori *estimator of $X$, given an observation $y \in \mathcal{Y}$, if it exists, is any which returns its mode under the posterior distribution $P_X(\cdot \mid Y = y)$.*

When the variable of interest $X$ is a discrete random variable, i.e., when there is a countable subset of $\mathcal{X}$ with $P_X$-measure one, the MAP estimator is the Bayesian estimator associated with the 0–1 loss:

$$L_{01}(\hat{x}, x) := \begin{cases} 0, & \text{if } x = \hat{x} \\ 1, & \text{if } x \neq \hat{x}. \end{cases}$$

For these variables, the expected posterior 0–1 loss can be written in terms of the posterior probability mass function,

$$\mathrm{E}[L_{01}(\hat{x}, X) \mid Y = y] = 1 - P(X = \hat{x} \mid Y = y)\,, \tag{2.3}$$

and the posterior modes, according to Definition 2.2, are the Bayesian estimates.

For random variables over general metric spaces $(\mathcal{X}, d)$, however, the MAP estimator is usually not associated with any single loss function. It can, nonetheless, be associated with a series of loss functions which approximate the 0–1 loss:

$$L_\epsilon(\hat{x}, x) := \begin{cases} 0, & \text{if } d(\hat{x}, x) < \epsilon \\ 1, & \text{if } d(\hat{x}, x) \geq \epsilon. \end{cases}$$

Similarly to what was done to the 0–1 loss in (2.3), the expected posterior $L_\epsilon$ loss can be written in terms of the probability of $\epsilon$-balls:

$$\mathrm{E}[L_\epsilon(x, X) \mid Y = y] = 1 - P\big(d(x, X) < \epsilon \,\big|\, Y = y\big)\,.$$

Recalling the interpretation of the mode that followed from Proposition 2.6, we can then place MAP estimation in the framework of Bayesian decision theory and Bayesian estimation point estimation. If $\hat{x}$ is a MAP estimate and $x$ is not, then there exists an $\bar{\epsilon} \in \mathbb{R}_{>0}$ such that for all $\epsilon \in (0, \bar{\epsilon}]$ the expected posterior $L_\epsilon$ loss associated with $\hat{x}$ is smaller than that associated with $x$.

We note that many authors define the MAP estimate, in the context of Bayesian decision theory, as the limit of the Bayesian estimates associated with $L_\epsilon$, as $\epsilon \to 0$ (Robert, 2001, Sec. 4.1.2; Mitchell, 2012, Sec. 10.2; Webb, 1999, Sec. B.1.4). This definition is only applicable, however, when the posterior distribution of the variable of interest is unimodal. If the posterior distribution is multimodal, then the limit of any convergent subsequence of Bayesian estimates associated with $L_\epsilon$ is a MAP estimate, under some regularity conditions (Daniel, 1969, p. 30). However, it is possible that not all modes are, necessarily, limits of Bayesian estimates associated with $L_\epsilon$. Because of these limitations, we believe that Definition 2.7 is preferable, due to its greater generality. Hypo-convergence and variational analysis, which are used in Chapter 3 to relate sequences of maximization problems, are powerful tools to better understand the relation between the MAP estimator and the Bayesian estimators associated with $L_\epsilon$.

In the sequence, we derive the joint posterior fictitious density of state paths and parameters of systems described by stochastic differential equations and, consequently, their MAP estimators.

## 2.2 Joint MAP state path and parameter estimation in SDEs

In this section, we apply the concepts developed in Section 2.1 for joint MAP state-path and parameter estimation in systems described by SDEs. In Section 2.2.1 the prior joint fictitious density is derived and in Section 2.2.2 it is used to obtain the posterior joint fictitious density.

### 2.2.1 Prior fictitious density

We now show how the Onsager–Machlup functional can be used to construct the joint fictitious density of state paths and parameters of Itō diffusion processes, under some regularity conditions. This functional was initially proposed by Onsager and Machlup (1953) for Gaussian diffusions, in the context of thermodynamic systems. Tisza and Manning (1957) later interpreted the maxima of this functional as the "most probable region" of paths of diffusion processes, which has since then become a consensus in various fields such as physics (Adib, 2008; Dürr and Bach, 1978), mathematical statistics (Takahashi and Watanabe, 1981; Ikeda and Watanabe, 1981, Sec. VI.9; Zeitouni and Dembo, 1987) and engineering (Aihara and Bagchi, 1999a,b; Dutra et al., 2014).

Stratonovich (1971) extended these ideas to diffusions with nonlinear drift and provided a rigorous mathematical background. Other developments related the Onsager–Machlup functional include its generalization to diffusions on manifolds (Takahashi and Watanabe, 1981); the extension of its domain to the Cameron–Martin space of paths (Shepp and Zeitouni, 1992); and proving that it also corresponds to the mode, according to Definition 2.5, with respect to other metrics besides the one induced by the supremum norm (Capitaine, 1995, 2000; Shepp and Zeitouni, 1993). In this section, we extend the Onsager–Machlup functional for diffusions with unknown parameters in the drift function and rank-deficient diffusion matrices. Our approach for handling rank-deficient diffusion matrices is similar to that of Aihara and Bagchi (1999a,b), but we note that their proofs are incomplete.

To begin, let $X$ and $Z$ be $\mathbb{R}^n$- and $\mathbb{R}^q$-valued stochastic processes, respectively, representing the state of a dynamical system over the time interval $\mathcal{T} := [0, t_\mathrm{f}]$. We consider systems whose dynamics are given by a system of stochastic differential equations (SDEs) of the form

$$\mathrm{d}X_t = f(t, X_t, Z_t, \Theta)\,\mathrm{d}t + \boldsymbol{G}\,\mathrm{d}W_t, \tag{2.4a}$$
$$\mathrm{d}Z_t = h(t, X_t, Z_t, \Theta)\,\mathrm{d}t, \tag{2.4b}$$

where $t \in \mathcal{T}$ is the time instant; $f\colon \mathcal{T} \times \mathbb{R}^n \times \mathbb{R}^q \times \mathbb{R}^m \to \mathbb{R}^n$ and $h\colon \mathcal{T} \times \mathbb{R}^n \times \mathbb{R}^q \times \mathbb{R}^m \to \mathbb{R}^q$ are the drift functions; the $\mathbb{R}^m$-valued random variable $\Theta$ is the unknown parameter vector; $W$ is an $n$-dimensional Wiener process over $\mathbb{R}_{\geq 0}$, with respect to the filtration $\{\mathcal{E}_t\}_{t\geq 0}$, representing the process noise; and $\boldsymbol{G} \in \mathbb{R}^{n\times n}$ is the diffusion matrix. As the process $X$ is under direct influence of noise and $Z$ is not, they will be denoted the *noisy* and *clean* state processes, respectively. Note that since the diffusion matrix does not depend on the state, the SDE can be interpreted in both the Itō or Stratonovich senses.

The following assumptions will be made on the system dynamics and prior distribution.

**Assumption 2.8** (prior and system dynamics)**.**

*a. The initial states $X_0$ and $Z_0$ and the parameter vector $\Theta$ are $\mathcal{E}_0$-measurable.*

*b. The initial states $X_0$ and $Z_0$ and the parameter vector $\Theta$ admit a continuous joint prior density $\pi\colon \mathbb{R}^n \times \mathbb{R}^q \times \mathbb{R}^m \to \mathbb{R}_{\geq 0}$.*

*c. The drift functions $f$ and $h$ are uniformly continuous with respect to all of their arguments in $\mathcal{T} \times \mathbb{R}^n \times \mathbb{R}^q \times \operatorname{supp}(P_\Theta)$, where $\operatorname{supp}(P_\Theta)$ indicates the topological support of the probability measure induced by $\Theta$, i.e., there exists $\rho_f\colon \mathbb{R}_{\geq 0} \to \mathbb{R}_{\geq 0}$ and $\rho_h\colon \mathbb{R}_{\geq 0} \to \mathbb{R}_{\geq 0}$ such that $\lim_{\epsilon\downarrow 0} \rho_f(\epsilon) + \rho_h(\epsilon) = 0$ and*

$$\left|f(t, x, z, \theta) - f(t', x', z', \theta')\right| \leq \rho_f(\epsilon), \tag{2.5a}$$
$$\left|h(t, x, z, \theta) - h(t', x', z', \theta')\right| \leq \rho_h(\epsilon), \tag{2.5b}$$

*for all $t, t' \in \mathcal{T}$, $x, x' \in \mathbb{R}^n$, $z, z' \in \mathbb{R}^q$ and $\theta, \theta' \in \operatorname{supp}(P_\Theta)$ such that $|t - t'| \leq \epsilon$, $|x - x'| \leq \epsilon$, $|z - z'| \leq \epsilon$ and $|\theta - \theta'| \leq \epsilon$.*

*d. For each $\theta \in \operatorname{supp}(P_\Theta)$, the drift functions $f$ and $h$ are Lipschitz continuous with respect to their second and third arguments $x$ and $z$, uniformly over their first argument $t$, i.e., there exist $L_{\mathrm{f}}^\theta, L_{\mathrm{h}}^\theta \in \mathbb{R}_{>0}$ such that for all $t \in \mathcal{T}$, $x, x' \in \mathbb{R}^n$ and $z, z' \in \mathbb{R}^q$,*

$$\left|f(t, x, z, \theta) - f(t, x', z', \theta)\right| \leq L_{\mathrm{f}}^\theta \left(\left|x - x'\right| + \left|z - z'\right|\right), \tag{2.6a}$$
$$\left|h(t, x, z, \theta) - h(t, x', z', \theta)\right| \leq L_{\mathrm{h}}^\theta \left(\left|x - x'\right| + \left|z - z'\right|\right). \tag{2.6b}$$

*e. For all $t \in \mathcal{T}$ and $\theta \in \operatorname{supp}(P_\Theta)$, the noisy drift function $f$ is twice differentiable with respect to its second argument $x$ and differentiable with respect to its first and third arguments $t$ and $z$. Furthermore, its first and second derivatives mentioned above are continuous with respect to all their arguments.*

f. *The diffusion matrix $\boldsymbol{G}$ has full rank.*

g. *The drift function $f$, the diffusion matrix $\boldsymbol{G}$, the state processes $X$ and $Z$, and the parameter vector $\Theta$ are such that*

$$\mathrm{E}\left[\exp\left(\int_0^1 \left|\boldsymbol{G}^{-1} f(t, X_t, Z_t, \Theta)\right|^2 \mathrm{d}t\right)\right] < \infty.$$

Assumption 2.8a ensures that the state process $X$ is adapted to the filtration $\{\mathcal{E}_t\}_{t\geq 0}$. Assumption 2.8b, like the requirement of continuity of densities in Definition 2.3, makes $\pi$ the *best* joint density for $X_0$, $Z_0$ and $\Theta$, as its Lebesgue points are its whole domain. Assumption 2.8d ensures, together with Lemma A.9 of Appendix A, that there exists, almost surely, a unique strong solution to the system of SDEs (2.4). Assumption 2.8f serves to enforce the consistency of the division of the clean and noisy states of (2.4). Assumption 2.8g is known as Novikov's condition, and can be interpreted as a requirement that the expected "energy" exerted by and upon the system during the experiment is finite (cf. Stummer, 1993). The easiest way to guarantee that it holds is ensuring that $f$ is bounded. More general requirements on the drift $f$ to guarantee that Novikov's condition is satisfied can be obtained by using Corollary 3.5.16 of Karatzas and Shreve (1991) together with Lemma A.9. The remaining assumptions are used to ensure the existence of the fictitious density.

Having stated the assumptions we have made on the system dynamics, we are now ready to derive the joint fictitious density of the $X$ path and the $Z_0$ and $\Theta$ vectors. We first state the theorem then a series of lemmas which are used for its proof, which is presented at the end of the subsection, on page 38. As mentioned in the beginning of this section, the fictitious density of paths of diffusion processes can be taken with respect to many norms of function spaces, all of them leading to the same density. However, the assumptions used to guarantee the existence of the fictitious density depends on the norm. In this thesis, we choose the supremum norm (also known as the uniform or infinity norm), as it requires the weakest assumptions. This metric is a natural choice as it is the metric of the classical Wiener space.

The supremum norm is defined by

$$|||x||| := \sup\nolimits_{t\in\mathcal{T}} \|x(t)\|, \tag{2.7}$$

where any finite-dimensional norm can be used on the right-hand side. Under this metric the $\epsilon$-balls represent tubes (sometimes also referred to as sausages) of radius $\epsilon$, like the ones represented graphically in Figure 2.1. The choice of the underlying finite-dimensional norm in (2.7) defines the geometry of the tubes' transversal cross section. The Euclidean distance is implied unless otherwise noted, leading to circular cross sections.

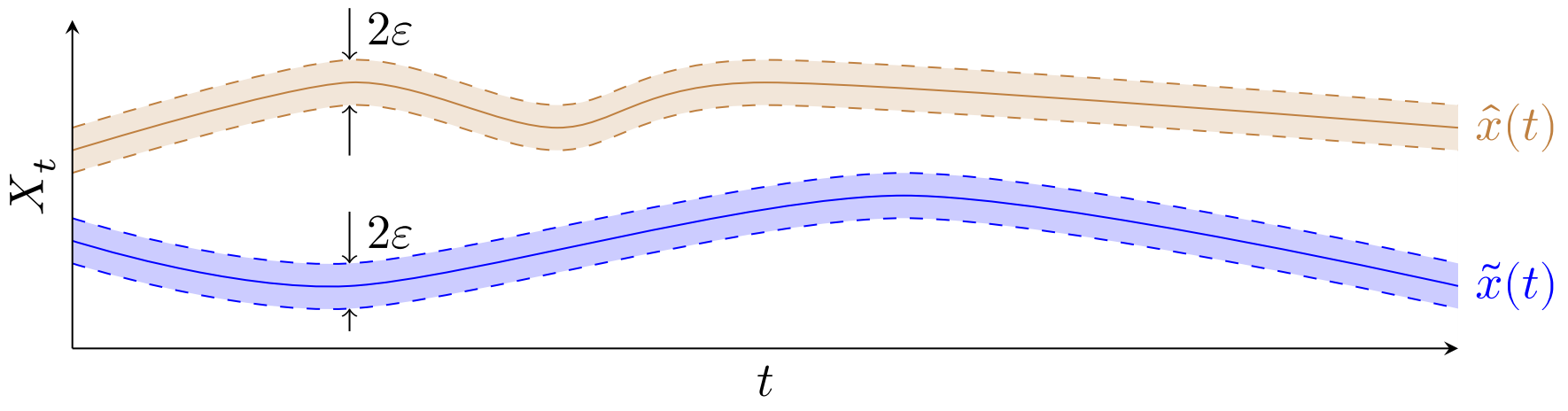

Figure 2.1: Graphical representation of $\epsilon$-balls under the supremum norm. The dashed lines correspond to the centers and all sample paths of the $\mathbb{R}$-valued process in the shaded region belong to the ball.

There are two main approaches to derive the Onsager–Machlup functional. One is acomplished by performing a Taylor expansion of the drift process of the SDE, and the other by using the stochastic Stokes' theorem, as pioneered by Takahashi and Watanabe (1981). Once more, in this thesis we use the stochastic Stokes route since it requires weaker assumptions, leading to a wider applicability of the theorem. However, to use the stochastic Stokes' theorem, the underlying finite-dimensional norm of the supremum norm must be the Mahalanobis distance associated with the diffusion matrix, i.e., the norm with respect to which the fictitious density of the $X$ paths is taken must be

$$\|x\|_{\boldsymbol{G}} := \sup_{t\in\mathcal{T}} \left|\boldsymbol{G}^{-1}x(t)\right|. \tag{2.8}$$

Most of the results of this chapter still hold for the supremum norm with other underlying finite-dimensional norms, albeit with different assumptions. In particular, the Jacobian matrix of the noisy drift $f$ with respect to $x$ must be Lipschitz continuous with respect to $x$ (cf. Capitaine, 1995).

Throughout this section, we will denote by $\{\mathbb{B}_\epsilon\}_{\epsilon\in\mathbb{R}_{>0}} \subset \mathcal{E}$ the family of events, indexed by $\epsilon \in \mathbb{R}_{>0}$, corresponding to the outcomes in which $X$, $Z_0$ and $\Theta$ are inside $\epsilon$-balls centered around the test points $x \in \mathcal{C}(\mathcal{T}, \mathbb{R}^n)$, $z_0 \in \mathbb{R}^q$, and $\theta \in \mathbb{R}^m$, respectively, i.e.,

$$\mathbb{B}_\epsilon := \Big\{\omega \in \Omega \Big| \ \|X(\omega) - x\|_{\boldsymbol{G}} < \epsilon, |Z_0(\omega) - z_0| < \epsilon, |\Theta(\omega) - \theta| < \epsilon\Big\}. \tag{2.9}$$

In addition, $\mathcal{W}_n^2$ will denote the Sobolev space of absolutely continuous functions from $\mathcal{T}$ to $\mathbb{R}^n$ whose weak derivative is in $L_n^2(\mathcal{T})$. The joint fictitious density of $X$, $Z_0$ and $\Theta$ is then given by the theorem below.

**Theorem 2.9** (joint fictitious density of state paths and parameters)**.** *If Assumption 2.8 is satisfied, then there exist* $a_1, a_2 \in \mathbb{R}_{>0}$ *such that for all* $x \in \mathcal{W}_n^2$*,* $z_0 \in \mathbb{R}^q$*, and* $\theta \in \mathbb{R}^m$*,*

$$\lim_{\epsilon\downarrow 0} \frac{P(\mathbb{B}_\epsilon)}{a_1 \exp(-\frac{a_2}{\epsilon^2})\,\epsilon^{m+n+q}} = \pi(x(0), z_0, \theta) \exp\Big(J(x, z_0, \theta)\Big), \tag{2.10}$$

*where $\mathbb{B}_\epsilon$ is the $\epsilon$-ball centered in $x$, $z_0$, and $\theta$, defined in (2.9), and the the Onsager–Machlup functional $J\colon \mathcal{W}_n^2 \times \mathbb{R}^q \times \mathbb{R}^m \to \mathbb{R}$ is defined by*

$$J(x, z_0, \theta) := -\frac{1}{2} \int_{\mathcal{T}} \left|\boldsymbol{G}^{-1} \left[f(t, x(t), z(t), \theta) - \dot{x}(t)\right]\right|^2 \mathrm{d}t - \frac{1}{2} \int_{\mathcal{T}} \operatorname{div}_{\mathrm{x}} f(t, x(t), z(t), \theta) \, \mathrm{d}t, \quad (2.11)$$

*and $z \in \mathcal{W}_q^2$ is the solution to the initial value problem*

$$\dot{z}(t) = h(t, x(t), z(t), \theta)\,, \qquad z(0) = z_0. \quad (2.12)$$

We now present a series of helping lemmas which culminate in the proof of Theorem 2.9 on page 38. We begin by remarking that if the theorem holds for the unit time interval, then it holds for any $t_\mathrm{f} \in \mathbb{R}_{>0}$. This is done so that the subsequent lemmas can be made simpler by considering $t_\mathrm{f} = 1$.

**Remark 2.10** (reduction to the unit interval)**.** *If Theorem 2.9 holds for $t_\mathrm{f} = 1$, then it holds for any $t_\mathrm{f} \in \mathbb{R}_{>0}$.*

*Proof.* If $\mathcal{T}_\mathrm{u} := [0, 1]$ denotes the unit time interval, then the linear map $u(\tau) := \tau t_\mathrm{f}$ is a diffeomorphism between $\mathcal{T}_\mathrm{u}$ and $\mathcal{T}$. Defining the processes $\hat{X}\colon \mathcal{T}_\mathrm{u} \times \Omega \to \mathbb{R}^n$ and $\hat{Z}\colon \mathcal{T}_\mathrm{u} \times \Omega \to \mathbb{R}^q$ by

$$\hat{X}(\tau, \omega) := X(u(\tau), \omega), \qquad \hat{Z}(\tau, \omega) := Z(u(\tau), \omega),$$

we have that they satisfy the following SDEs:

$$\begin{aligned} \mathrm{d}\hat{X}_\tau &:= t_\mathrm{f} f\Big(u(\tau), \hat{X}_\tau, \hat{Z}_\tau, \Theta\Big) \, \mathrm{d}\tau + \sqrt{t_\mathrm{f}}\, \boldsymbol{G} \, \mathrm{d}\hat{W}_\tau, \\ \mathrm{d}\hat{Z}_\tau &:= t_\mathrm{f} h\Big(u(\tau), \hat{X}_\tau, \hat{Z}_\tau, \Theta\Big) \, \mathrm{d}\tau, \end{aligned}$$

where, because of the Wiener process scaling property, the process

$$\hat{W}(\tau, \omega) := \tfrac{1}{\sqrt{t_\mathrm{f}}} W(t_\mathrm{f}\tau, \, \omega)$$

is an $n$-dimensional Wiener process with respect to a time-changed filtration. With this transformation we have that, for all $x\colon \mathcal{T} \to \mathbb{R}^n$ and $\epsilon \in \mathbb{R}_{>0}$,

$$P\Big(\sup\nolimits_{t \in \mathcal{T}} |\boldsymbol{G}^{-1}[x(t) - X_t]| < \epsilon\Big) = P\Big(\sup\nolimits_{\tau \in \mathcal{T}_\mathrm{u}} |\boldsymbol{G}^{-1}[x(\tau t_\mathrm{f}) - \hat{X}_\tau]| < \epsilon\Big)\,,$$

implying that the joint fictitious density of $X$ over $\mathcal{T}$, $Z_0$ and $\Theta$ is the same as that of $\hat{X}$ over $\mathcal{T}_\mathrm{u}$, $\hat{Z}_0$ and $\Theta$.

Next, note that if the original system satisfies Assumption 2.8, then the time-changed one satisfies it as well. If Theorem 2.9 holds for the unit interval, we then have that the joint fictitious density of both systems is given by (2.10), whith the Onsager–Machlup functional $J$ given by

$$J(x, z_0, \theta) := -\frac{1}{2}\int_{\mathcal{T}_{\mathrm{u}}} \operatorname{div}_{\mathrm{x}} f(u(\tau), z(u(\tau), x(u(\tau)), \theta)\, t_{\mathrm{f}}\, \mathrm{d}\tau \\ - \frac{1}{2}\int_{\mathcal{T}_{\mathrm{u}}} \left|\boldsymbol{G}^{-1}\left[f(u(\tau), z(u(\tau)), x(u(\tau)), \theta) - \dot{x}(u(\tau))\right]\right|^2 t_{\mathrm{f}}\, \mathrm{d}\tau.$$

Performing a change in the region of integration using $u$, we obtain the same expression for $J$ as in (2.11), since $\mathrm{d}t = t_{\mathrm{f}}\mathrm{d}\tau$. □

Using Remark 2.10, we will now consider $t_{\mathrm{f}} = 1$ for simplicity. In addition, to further simplify the notation and the equations that follow, we define the $\mathbb{R}^n$- and $\mathbb{R}^q$-valued stochastic processes $F$ and $H$ as

$$F_t := f(t, X_t, Z_t, \Theta)\,, \qquad H_t := h(t, X_t, Z_t, \Theta)\,, \tag{2.13}$$

and the functions $\hat{f}\colon \mathcal{T} \to \mathbb{R}^n$ and $\hat{h}\colon \mathcal{T} \to \mathbb{R}^q$ as

$$\hat{f}(t) := f(t, x(t), z(t), \theta)\,, \qquad \hat{h}(t) := h(t, x(t), z(t), \theta)\,, \tag{2.14}$$

where $z$ is the solution to the initial value problem (2.12).

We now present a lemma which states that if $X$, $Z_0$ and $\Theta$ are contained in $\epsilon$-balls, then the clean state path $Z$ is also contained in an $\varepsilon$-ball with respect to the supremum norm. The center of this ball is the solution to the initial value problem (2.12) and its radius $\varepsilon$ vanishes when $\epsilon$ vanishes.

**Lemma 2.11** (clean states growth bound)**.** *For a system satisfying Assumption 2.8, then for all $x \in \mathcal{C}(\mathcal{T}, \mathbb{R}^n)$, $z_0 \in \mathbb{R}^q$, and $\theta \in \operatorname{supp}(P_\Theta)$ there exists $a_3 \in \mathbb{R}_{>0}$ such that, almost surely,*

$$|\!|\!| Z - z |\!|\!| \le a_3\left[\rho_f(|\Theta - \theta|) + |Z_0 - z_0| + |\!|\!| X - x |\!|\!|\right],$$

*where $z \in \mathcal{W}_q^2$ is the solution to the initial value problem* (2.12) *and $\rho_f$ is the modulus of continuity of $f$, which is assumed to exist by Assumption 2.8c.*

*Proof.* Using the process $H$ and the function $\hat{h}$ defined in (2.13) and (2.14), the differential equations (2.4b) and (2.12) admit the following representation in their integral forms:

$$Z_t = Z_0 + \int_0^t H_\tau\, \mathrm{d}\tau, \qquad z(t) = z_0 + \int_0^t \hat{h}(\tau)\, \mathrm{d}\tau.$$

Using the triangle inequality we then have that, for all $t \in \mathcal{T}$, the difference between $Z$ and $z$ is bounded by

$$|Z_t - z(t)| \le |Z_0 - z_0| + \int_0^t \left|H_\tau - \hat{h}(\tau)\right| \mathrm{d}\tau. \tag{2.15}$$

The difference between $H$ and $\hat{h}$, on the other hand, can be bounded almost surely:

$$\begin{aligned}|H_t - \hat{h}(t)| &= |H_t - h(t, X_t, Z_t, \theta) + h(t, Z_t, X_t, \theta) - \hat{h}(t)| \\ &\leq |H_t - h(t, X_t, Z_t, \theta)| + |h(t, X_t, Z_t, \theta) - \hat{h}(t)| \qquad (2.16a) \\ &\leq \rho_h(|\Theta - \theta|) + L_{\mathrm{h}}^{\theta} |X_t - x(t)| + L_{\mathrm{h}}^{\theta} |Z_t - z(t)|, \qquad (2.16b)\end{aligned}$$

where to obtain (2.16a) the triangle inequality was used and to obtain (2.16b) the Assumptions 2.8c–d on the uniform and Lipschitz continuity assumptions of $h$ were used. Also note that (2.16b) only holds almost surely because $h$ is only assumed to be uniformly continuous for $\theta$ on the support of $\Theta$.

Substituting (2.16b) into (2.15), we then obtain

$$|Z_t - z(t)| \leq L_{\mathrm{h}}^{\theta} \, |\!|\!|X - x|\!|\!| + \rho_h(|\Theta - \theta|) + |Z_0 - z_0| + L_{\mathrm{h}}^{\theta} \int_0^t |Z_\tau - z(t)| \, \mathrm{d}\tau$$

Applying the Grönwall–Bellman inequality, Lemma A.8, we then have that

$$|\!|\!|Z - z|\!|\!| \leq [L_{\mathrm{h}}^{\theta} \, |\!|\!|X - x|\!|\!| + \rho_h(|\Theta - \theta|) + |Z_0 - z_0|] \exp(L_{\mathrm{h}}^{\theta}). \qquad \square$$

The next helping lemma is a change of measure which is used to represent the $\epsilon$-ball probabilities as conditional expectations of a random variable.

**Lemma 2.12** (change of measure)**.** *For a system satisfying Assumption 2.8 and a given $x \in \mathcal{W}_n^2$, define the $\mathbb{R}$-valued random variable $M$ and the $\mathbb{R}^n$-valued stochastic processes $U$, $\tilde{W}$ and $\tilde{W}^{\circ}$ as*

$$U_t := \boldsymbol{G}^{-1} [F_t - \dot{x}(t)] \qquad (2.17a)$$

$$\tilde{W}_t := \boldsymbol{G}^{-1} [X_t - x(t)] \qquad (2.17b)$$

$$\tilde{W}_t^{\circ} := \tilde{W}_t - \tilde{W}_0 \qquad (2.17c)$$

$$M := \exp\Big(-\int_0^1 U_t^{\mathsf{T}} \, \mathrm{d}W_t - \frac{1}{2} \int_0^1 |U_t|^2 \, \mathrm{d}t\Big) \qquad (2.17d)$$

*and the measure $\tilde{P} \colon \mathcal{E} \to \mathbb{R}_{\geq 0}$ as*

$$\tilde{P}(\mathbb{E}) := \int_{\mathbb{E}} M(\omega) \, \mathrm{d}P(\omega).$$

*Then $\tilde{P}$ is a probability measure equivalent*[4] *to $P$. Furthermore, $\tilde{W}^{\circ}$ is a $n$-dimensional Wiener process with respect to the filtration $\{\mathcal{E}_t\}_{t \geq 0}$ under $\tilde{P}$ and*

$$M^{-1} = \exp\Big(\textstyle\int_0^1 U_t^{\mathsf{T}} \, \mathrm{d}\tilde{W}_t - \frac{1}{2} \int_0^1 |U_t|^2 \, \mathrm{d}t\Big) \quad \text{almost surely,} \qquad (2.18a)$$

$$P(\mathbb{E}) = \textstyle\int_{\mathbb{E}} M^{-1}(\omega) \, \mathrm{d}\tilde{P}(\omega) \quad \text{for all } \mathbb{E} \in \mathcal{E}, \qquad (2.18b)$$

$$\tilde{P}(\mathbb{E}) = P(\mathbb{E}) \quad \text{for all } \mathbb{E} \in \mathcal{E}_0. \qquad (2.18c)$$

[4]That is, mutually absolutely continuous with respect to one another.

*Proof.* This is simply an application of the Girsanov transformation of measures (Øksendal, 2003, Sec. 8.6). By substituting (2.4a) into (2.17b), we have that $\tilde{W}$ satisfies the SDE

$$\mathrm{d}\tilde{W}_t = \boldsymbol{G}^{-1}\left[f(t, Z_t, X_t, \Theta) - \dot{x}(t)\right]\mathrm{d}t + \mathrm{d}W_t. \tag{2.19}$$

Novikov's condition for the the removal of the drift of (2.19) is

$$\begin{aligned}
&\mathrm{E}\Big[\exp\Big(\tfrac{1}{2}\textstyle\int_0^1 \big|\boldsymbol{G}^{-1}[F_t - \dot{x}(t)]\big|^2\,\mathrm{d}t\Big)\Big] \\
&\quad \leq \mathrm{E}\Big[\exp\Big(\tfrac{1}{2}\textstyle\int_0^1 [|\boldsymbol{G}^{-1}F_t| + |\boldsymbol{G}^{-1}\dot{x}(t)|]^2\mathrm{d}t\Big)\Big] && (2.20\text{a}) \\
&\quad \leq \mathrm{E}\Big[\exp\Big(\textstyle\int_0^1 \big|\boldsymbol{G}^{-1}f(t, X_t, Z_t, \Theta)\big|^2\,\mathrm{d}t\Big)\Big]\exp\Big(\textstyle\int_0^1 |\boldsymbol{G}^{-1}\dot{x}(t)|^2\mathrm{d}t\Big) && (2.20\text{b}) \\
&\quad < \infty. && (2.20\text{c})
\end{aligned}$$

Equation (2.20a) is obtained from the triangle inequality, (2.20b) by applying the simple identity of Lemma A.2, and (2.20c) from Assumption 2.8g and the fact that $x \in \mathcal{W}_n^2$. Applying the Girsanov transformation of measures (Øksendal, 2003, Sec. 8.6) we then obtain the results stated in the lemma. □

Next, we obtain the fictitious uniform measure of $\epsilon$-balls with respect to which the joint fictitious density of the state-paths and parameters is taken, i.e., the denominator of the left-hand side of (2.1) in the definition of fictitious densities (Definition 2.4).

**Lemma 2.13** ($\epsilon$-ball fictitious uniform measure)**.** *For a system satisfying Assumption 2.8, there exist $a_1, a_2 \in \mathbb{R}_{>0}$ such that, for all $x \in \mathcal{W}_n^2$ with $x(0) \in \operatorname{supp}(P_{X_0})$, $z_0 \in \operatorname{supp}(P_{Z_0})$, and $\theta \in \operatorname{supp}(P_\Theta)$,*

$$\lim_{\epsilon\downarrow 0}\frac{P(\mathbb{B}_\epsilon)}{a_1\exp(-\frac{a_2}{\epsilon^2})\,\epsilon^{m+n+q}} = \pi(x(0), z_0, \theta)\lim_{\epsilon\downarrow 0}\tilde{\mathrm{E}}\big[M^{-1}\,\big|\,\mathbb{B}_\epsilon\big], \tag{2.21}$$

*whenever the limit on the right-hand side of* (2.21) *exists. Additionally, $P(\mathbb{B}_\epsilon) > 0$ for all $\epsilon \in \mathbb{R}_{>0}$, so the conditional expectation is well-defined.*

*Proof.* Applying Lemma 2.12 we can obtain the probability of $\mathbb{B}_\epsilon$ in terms of $M^{-1}$ and the measure $\tilde{P}$. For any event $\mathbb{A} \in \mathcal{E}$ such that $\tilde{P}(\mathbb{A}) > 0$, (2.18b) implies that

$$P(\mathbb{A}) = \int_{\mathbb{A}} M^{-1}\mathrm{d}\tilde{P}(\omega) = \tilde{\mathrm{E}}\big[M^{-1}\,\big|\,\mathbb{A}\big]\,\tilde{P}(\mathbb{A}).$$

Hence, to prove that (2.21) holds we need to show that, under the conditions specified in the Lemma's statement, $\tilde{P}(\mathbb{B}_\epsilon) > 0$ for all $\epsilon \in \mathbb{R}_{>0}$ and the following limit holds:

$$\lim_{\epsilon\downarrow 0}\frac{\tilde{P}(\mathbb{B}_\epsilon)}{\exp(-\frac{a_2}{\epsilon^2})\,\epsilon^{m+n+q}} = \pi(x(0), z_0, \theta)\,a_1. \tag{2.22}$$

To see that (2.22) indeed holds, we begin by obtaining an expression for $\tilde{P}(\mathbb{B}_\epsilon)$. From Assumption 2.8a, we have that $X_0$, $Z_0$, and $\Theta$ are $\mathcal{E}_0$-measurable. Additionally, from (2.18c) of the change of measure lemma (Lemma 2.12), we know that $\tilde{P}$ and $P$ coincide on $\mathcal{E}_0$-events. Consequently, Assumption 2.8b implies that the joint prior density of $X_0$, $Z_0$, and $\Theta$ under $\tilde{P}$ is also given by $\pi$. Finally, from Lemma 2.12 we have that, given $X_0$, $X$ is independent of $\mathcal{E}_0$ under $\tilde{P}$ and, as a consequence, independent of $Z_0$ and $\Theta$ under $\tilde{P}$. Using the definition of the $\mathbb{B}_\epsilon$ event in (2.9), we then have its $\tilde{P}$ measure is given by

$$\tilde{P}(\mathbb{B}_\epsilon) = \iiint\limits_{\mathbb{Q}_\epsilon \mathbb{R}_\epsilon \mathbb{S}_\epsilon} \pi(\hat{x} + x(0), \hat{z} + z_0, \hat{\theta} + \theta) \times \tilde{P}\Big( |\!|\!| X - x |\!|\!|_{\boldsymbol{G}} < \epsilon \,\Big|\, X_0 = \hat{x} + x(0) \Big) \,\mathrm{d}\hat{x}\,\mathrm{d}\hat{z}\,\mathrm{d}\hat{\theta}, \quad (2.23)$$

where the family of sets $\mathbb{Q}_\epsilon$, $\mathbb{R}_\epsilon$, and $\mathbb{S}_\epsilon$ are zero-centered $\epsilon$-balls defined as

$$\begin{aligned} \mathbb{Q}_\epsilon &:= \left\{ x \in \mathbb{R}^n \,\middle|\, |\boldsymbol{G}^{-1} x| < \epsilon \right\}, \\ \mathbb{R}_\epsilon &:= \left\{ z \in \mathbb{R}^q \,\middle|\, |z| < \epsilon \right\}, \\ \mathbb{S}_\epsilon &:= \left\{ \theta \in \mathbb{R}^m \,\middle|\, |\theta| < \epsilon \right\}. \end{aligned}$$

From the definition of the $|\!|\!|\cdot|\!|\!|_{\boldsymbol{G}}$ norm in (2.8) and the definition of the $\tilde{W}$ process in (2.17b) we have that

$$|\!|\!| X - x |\!|\!|_{\boldsymbol{G}} = |\!|\!| \tilde{W} |\!|\!|,$$

implying that

$$\tilde{P}\Big( |\!|\!| X - x |\!|\!|_{\boldsymbol{G}} < \epsilon \,\Big|\, X_0 = \hat{x} + x(0) \Big) = \tilde{P}\Big( |\!|\!| \tilde{W} |\!|\!| < \epsilon \,\Big|\, \tilde{W}_0 = \boldsymbol{G}^{-1} \hat{x} \Big). \quad (2.24)$$

Using the Wiener process scaling property, we have that the process $\epsilon W_{t\epsilon^{-2}}$ has the same distribution as $W_t$, so the right-hand side of (2.24) can be further simplified

$$\begin{aligned} \tilde{P}\Big( |\!|\!| \tilde{W} |\!|\!| < \epsilon \,\Big|\, \tilde{W}_0 = \boldsymbol{G}^{-1} \hat{x} \Big) &= P\Big( \sup\nolimits_{0 \le t \le 1} |\epsilon W_{t\epsilon^{-2}}| < \epsilon \,\Big|\, W_0 = \boldsymbol{G}^{-1} \hat{x} \Big) \\ &= P\Big( \sup\nolimits_{0 \le \tau \le \epsilon^{-2}} |W_\tau| < 1 \,\Big|\, W_0 = \boldsymbol{G}^{-1} \hat{x} \Big). \quad (2.25) \end{aligned}$$

As it so happens, the right-hand side of (2.25) is the probability that the Wiener process starting at $\boldsymbol{G}^{-1}\hat{x}$ sojourns in the unit sphere

$$\mathbb{U} := \left\{ u \in \mathbb{R}^n \,\middle|\, |u| < 1 \right\}$$

for longer than $\epsilon^{-2}$, which according to Theorem A.15 satisfies a boundary value problem whose solution is given by Proposition A.17:

$$\tilde{P}\Big( |\!|\!| \tilde{W} |\!|\!| < \epsilon \,\Big|\, \tilde{W}_0 = \hat{w} \Big) = \sum_{k=1}^{\infty} \exp(-\lambda_k \epsilon^{-2}) \gamma_k(\hat{w}) \int_{\mathbb{U}} \gamma_k(w) \,\mathrm{d}w, \quad (2.26)$$

where $\{\gamma_k\}_{k=1}^{\infty}$ are the eigenfunctions and $\{\lambda_k\}_{k=1}^{\infty}$ the corresponding eigenvalues of the Dirichlet eigenvalue problem on the unit sphere $\mathbb{U}$, as defined in Lemma A.16.

From (2.26) we can see that, for all $\epsilon \in \mathbb{R}_{>0}$ and $\hat{w} \in \mathbb{U}$,

$$\tilde{P}\Big(|\!|\!|\tilde{W}|\!|\!| < \epsilon \,\Big|\, \tilde{W}_0 = \hat{w}\Big) > 0.$$

Consequently, from (2.23) we can then conclude that $\tilde{P}(\mathbb{B}_\epsilon) > 0$, as since $x(0)$, $z_0$, and $\theta$ are in their prior's support, $P(X_0 \in \mathbb{Q}_\epsilon, Z_0 \in \mathbb{R}_\epsilon, \Theta \in \mathbb{S}_\epsilon) > 0$ and the integral of a nonzero function over a positive measure set is always positive.

Finally, performing a change of variables, we have that (2.23) can be simplified to

$$\begin{aligned}\tilde{P}(\mathbb{B}_\epsilon) = \epsilon^{m+n+q}\,|\det \boldsymbol{G}| \iiint\limits_{\mathbb{U}\,\mathbb{R}_1\mathbb{S}_1} &\pi(\epsilon \boldsymbol{G}\hat{w} + x(0), \epsilon\hat{z} + z_0, \epsilon\hat{\theta} + \theta) \\ &\times \left(\sum_{k=1}^{\infty} \exp\Big(-\frac{\lambda_k}{\epsilon^2}\Big)\,\gamma_k(\hat{w}) \int_{\mathbb{U}} \gamma_k(w)\,\mathrm{d}w\right) \mathrm{d}\hat{w}\,\mathrm{d}\hat{z}\,\mathrm{d}\hat{\theta},\end{aligned}$$

Using the bounded convergence theorem we then have that (2.22) and (2.21) hold with

$$a_1 = \mu(\mathbb{R}_1)\mu(\mathbb{S}_1)\,|\det \boldsymbol{G}| \left(\int_{\mathbb{G}} \gamma_1(w)\,\mathrm{d}w\right)^2, \qquad a_2 = \lambda_1, \tag{2.27}$$

where $\mu$ denotes the Lebesgue measure. □

It should be noted that the constants obtained in (2.27) are consistent with those of the Onsager–Machlup with initial condition (Fujita and Kotani, 1982, p. 129). From Lemma 2.13 we have that for (2.10) in Theorem 2.9 to hold, the expectation in the right-hand side of (2.21) should converge to the exponential of the Onsager–Machlup functional. From the expression of $M^{-1}$ in (2.18a), we should expect that its quadratic term converges to the quadratic term of the Onsager–Machlup functional in (2.11). This is shown with the help of the Lemma below.

**Lemma 2.14** (energy term of the Onsager–Machlup functional). *If Assumption 2.8 is satisfied, then for all $c \in \mathbb{R}$, $x \in \mathcal{W}_n^2$ with $x(0) \in \operatorname{supp}(P_{X_0})$, $z_0 \in \operatorname{supp}(P_{Z_0})$, and $\theta \in \operatorname{supp}(P_\Theta)$,*

$$\limsup_{\epsilon\downarrow 0} \tilde{\mathrm{E}}\left[\exp\left(c\int_0^1 \Big(|U_t|^2 - \Big|\boldsymbol{G}^{-1}\,\big[\hat{f}(t) - \dot{x}(t)\big]\Big|^2\Big)\,\mathrm{d}t\right)\,\middle|\, \mathbb{B}_\epsilon\right] \le 1. \tag{2.28}$$

*Proof.* To begin, define $\varepsilon \in \mathbb{R}_{>0}$ by

$$\varepsilon := a_3[\rho_f(\epsilon) + 2\epsilon] + \epsilon.$$

From Lemma 2.11 we then have that it bounds both $\|\!|X-x|\!\|$ and $\|\!|Z-z|\!\|$ for almost all $\omega \in \mathbb{B}_\epsilon$, i.e.,

$$\|\!|X-x|\!\| + \|\!|Z-z|\!\| \leq \varepsilon \qquad \text{for } P\text{-almost all } \omega \in \mathbb{B}_\epsilon.$$

Expanding the integrand of (2.28), we obtain

$$\begin{aligned} |U_t|^2 - \left|\boldsymbol{G}^{-1}\left[\hat{f}(t) - \dot{x}(t)\right]\right|^2 = \\ \left(|U_t| - \left|\boldsymbol{G}^{-1}\left[\hat{f}(t) - \dot{x}(t)\right]\right|\right)\left(|U_t| + \left|\boldsymbol{G}^{-1}\left[\hat{f}(t) - \dot{x}(t)\right]\right|\right). \end{aligned} \tag{2.29}$$

Using the reverse triangle inequality, (2.5a) and (2.13), we then have that, for almost all $\omega \in \mathbb{B}_\epsilon$,

$$\Big||U_t| - |\boldsymbol{G}^{-1}[\hat{f}(t) - \dot{x}(t)]|\Big| \leq \left|\boldsymbol{G}^{-1}[F_t - \hat{f}(t)]\right| \leq \left\|\boldsymbol{G}^{-1}\right\| \rho_f(\varepsilon), \tag{2.30}$$

where the induced operator norm is used for $\boldsymbol{G}^{-1}$. Next, using the triangle inequality, we have that

$$\Big||U_t| + |\boldsymbol{G}^{-1}[\hat{f}(t) - \dot{x}(t)]|\Big| \leq \|\boldsymbol{G}^{-1}\|\left(|F_t| + 2\,|\dot{x}| + \left|\hat{f}(t)\right|\right). \tag{2.31}$$

A bound for the $F$ process can be found using the triangle inequality and the uniform continuity of $f$ as well:

$$|F_t| \leq |F_t - \hat{f}(t)| + |\hat{f}(t)| \leq \rho_f(\varepsilon) + |\hat{f}(t)|. \tag{2.32}$$

Next noting that $\hat{f}$ is bounded since it is continuous over a compact space and that $\dot{x} \in L^1_n$ according to Lemma A.6, define

$$a_4 := 2\|\!|\hat{f}|\!\| + 2\|\dot{x}\|_{L^n_1}, \qquad\qquad a_5 := \|\boldsymbol{G}^{-1}\|^2.$$

Then, collecting (2.29) to (2.32), we have that for almost all $\omega \in \mathbb{B}_\epsilon$

$$\int_0^1 \left(|U_t|^2 - \left|\boldsymbol{G}^{-1}\left[\hat{f}(t) - \dot{x}(t)\right]\right|^2\right)\mathrm{d}t \leq a_5 \rho_f(\varepsilon)[a_4 + \rho_f(\varepsilon)].$$

Letting $\epsilon$ vanish we can then see that (2.28) holds. □

Next, from the expression of $M^{-1}$ in (2.18a), we should expect that its Itō integral term converges to the drift divergence term of the Onsager–Machlup functional in (2.11). This is proved with the help of the lemmas that follow.

**Lemma 2.15** (conversion to Stratonovich integral)**.** *If Assumption 2.8 is satisfied, then*

$$\int_0^1 [\boldsymbol{G}^{-1}F_t]^\mathsf{T}\,\mathrm{d}\tilde{W}_t = \int_0^1 [\boldsymbol{G}^{-1}F_t]^\mathsf{T} \circ \mathrm{d}\tilde{W}_t - \frac{1}{2}\int_0^1 \mathrm{div}_\mathrm{x}\, f(t, X_t, Z_t, \Theta)\,\mathrm{d}t, \tag{2.33}$$

*where $\int A \circ \mathrm{d}B$ denotes the Stratonovich integral of $A$ with respect to $B$.*

*Proof.* From the definition of the $F$ and $H$ processes in (2.13) we have that its stochastic differential is given by

$$\begin{aligned}\mathrm{d}F_t = \tfrac{\partial f}{\partial t}(t, X_t, Z_t, \Theta)\,\mathrm{d}t + \nabla_{\mathrm{x}} f(t, X_t, Z_t, \Theta)\,\mathrm{d}X_t \\ + \nabla_{\mathrm{z}} f(t, X_t, Z_t, \Theta) H_t\,\mathrm{d}t\end{aligned}$$

In addition, from the definition of the $\tilde{W}$ process in (2.17b) we have that

$$X_t = \boldsymbol{G}\tilde{W}_t + x(t), \qquad \mathrm{d}X_t = \boldsymbol{G}\,\mathrm{d}\tilde{W}_t + \dot{x}(t)\,\mathrm{d}t, \tag{2.34}$$

Consequently,

$$\begin{aligned}\mathrm{d}\tilde{W}_t^{\mathsf{T}}\boldsymbol{G}^{-1}\,\mathrm{d}F_t &= \mathrm{d}\tilde{W}_t^{\mathsf{T}}\boldsymbol{G}^{-1}\nabla_{\mathrm{x}} f(t, X_t, Z_t, \Theta)\boldsymbol{G}\,\mathrm{d}\tilde{W}_t \\ &= \mathrm{tr}(\boldsymbol{G}^{-1}\nabla_{\mathrm{x}} f(t, X_t, Z_t, \Theta)\boldsymbol{G})\ \mathrm{d}t, \\ &= \mathrm{tr}(\nabla_{\mathrm{x}} f(t, X_t, Z_t, \Theta))\ \mathrm{d}t,\end{aligned} \tag{2.35}$$

where to obtain (2.35) we used the fact that the trace is invariant to similarity transformations. The quadratic covariation process of $\tilde{W}^{\mathsf{T}}$ and $\boldsymbol{G}^{-1}F$ is then given by

$$\left[\!\left[\boldsymbol{G}^{-1}F^{\mathsf{T}}, \tilde{W}\right]\!\right]_t = \int_0^t \mathrm{div}_{\mathrm{x}}\, f(t, X_t, Z_t, \Theta)\,\mathrm{d}t.$$

Converting the Itō integral to a Stratonovich integral we then obtain (2.33). □

**Lemma 2.16** (convergence of divergences)**.** *If Assumption 2.8 is satisfied, then for all* $c \in \mathbb{R}$*,* $x \in \mathcal{W}_n^2$ *with* $x(0) \in \mathrm{supp}(P_{X_0})$*,* $z_0 \in \mathrm{supp}(P_{Z_0})$*, and* $\theta \in \mathrm{supp}(P_\Theta)$*,*

$$\limsup_{\epsilon \downarrow 0} \tilde{\mathrm{E}}\Big[\exp\Big(c\int_0^1 \mathrm{tr}\,[\nabla_{\mathrm{x}} f(t, X_t, Z_t, \Theta) - \nabla_{\mathrm{x}} f(t, x(t), z(t), \theta)]\,\mathrm{d}t\Big)\ \Big|\ \mathbb{B}_\epsilon\Big] \le 1.$$

*Proof.* Without loss of generality, consider that $\epsilon \le 1$. Lemma 2.11 and the definition of the $\epsilon$-ball (2.9) then imply that $X$, $Z$, and $\Theta$ are bounded for almost all $\omega \in \mathbb{B}_\epsilon$. From Assumption 2.8e, we then have that the drift divergence is continuous, and continuous functions over compact spaces are uniformly continuous. Consequently,

$$\lim_{\epsilon \downarrow 0} \sup_{\omega \in \mathbb{B}_\epsilon} \sup_{t \in \mathcal{T}} |\mathrm{div}_{\mathrm{x}}\, f(t, X_t, Z_t, \Theta) - \mathrm{div}_{\mathrm{x}}\, f(t, x(t), z(t), \theta)| = 0. \qquad \square$$

**Lemma 2.17.** *If Assumption 2.8 is satisfied, then for all $x \in \mathcal{W}_n^2$ with $x(0) \in \operatorname{supp}(P_{X_0})$, $z_0 \in \operatorname{supp}(P_{Z_0})$, and $\theta \in \operatorname{supp}(P_\Theta)$,*

$$\begin{aligned}\int_0^1 [\boldsymbol{G}^{-1} F_t]^\mathsf{T} \circ \mathrm{d}\tilde{W}_t &= \bar{f}(1, \tilde{W}_1, \omega)^\mathsf{T} \tilde{W}_1 - \bar{f}(0, \tilde{W}_0, \omega)^\mathsf{T} \tilde{W}_0 \\ &\quad - \int_0^1 \frac{\partial \bar{f}}{\partial t}(t, \tilde{W}_t, \omega)^\mathsf{T} \tilde{W}_t \, \mathrm{d}t - \sum_{i,j=1}^n \int_0^1 \frac{\partial \bar{f}^{(i)}}{\partial x^{(j)}}(t, \tilde{W}_t, \omega) \, \mathrm{d}S_t^{(ij)} \\ &\quad - \frac{1}{2} \sum_{i,j=1}^n \int_0^1 \left( \frac{\partial^2 \bar{f}^{(i)}}{\partial x^{(j)} \partial x^{(j)}}(t, \tilde{W}_t, \omega) \tilde{W}^{(i)} - \frac{\partial^2 \bar{f}^{(i)}}{\partial x^{(j)} \partial x^{(i)}}(t, \tilde{W}_t, \omega) \tilde{W}^{(j)} \right) \mathrm{d}t,\end{aligned} \tag{2.36}$$

*where the function $\bar{f} \colon \mathcal{T} \times \mathbb{R}^n \times \Omega$ is defined as*

$$\bar{f}(t, \hat{x}, \omega) := \int_0^1 \boldsymbol{G}^{-1} f(t, \tau \boldsymbol{G} \hat{x} + x(t), Z_t, \Theta) \, \mathrm{d}\tau$$

*and the process $S \colon \mathcal{T} \times \Omega \to \mathbb{R}^{n \times n}$ is the stochastic area enclosed by $\tilde{W}$:*

$$S_t^{(ij)} := \int_0^1 \tilde{W}_\tau^{(i)} \circ \mathrm{d}\tilde{W}_\tau^{(j)} - \int_0^1 \tilde{W}_\tau^{(j)} \circ \mathrm{d}\tilde{W}_\tau^{(i)}.$$

*Proof.* From (2.34) and definition of the $F$ process in (2.13) we have that

$$F_t = f\Big(t, \boldsymbol{G}\tilde{W}_t + x(t), Z_t, \Theta\Big).$$

Consequently, applying the stochastic Stokes' theorem (Lemma A.18) on the space-time 1-form

$$\beta(t, \hat{x}) := \sum_{i=1}^n \left[\boldsymbol{G}^{-1} f(t, \boldsymbol{G}\hat{x} + x(t), Z_t, \Theta)\right]^{(i)} \mathrm{d}\hat{x}^{(i)},$$

we obtain that

$$\begin{aligned}\int_0^1 [\boldsymbol{G}^{-1} F_t]^\mathsf{T} \circ \mathrm{d}\tilde{W}_t &= \bar{f}(1, \tilde{W}_1, \omega)^\mathsf{T} \tilde{W}_1 - \bar{f}(0, \tilde{W}_0, \omega)^\mathsf{T} \tilde{W}_0 \\ &\quad - \int_0^1 \frac{\partial \bar{f}}{\partial t}(t, \tilde{W}_t, \omega)^\mathsf{T} \tilde{W}_t \, \mathrm{d}t - \sum_{i,j=1}^n \int_0^1 \frac{\partial \bar{f}^{(i)}}{\partial x^{(j)}}(t, \tilde{W}_t, \omega) \circ \mathrm{d}S_t^{(ij)}.\end{aligned} \tag{2.37}$$

Next, note that since the quadratic covariation $\left[\!\left[\tilde{W}^{(i)}, \tilde{W}^{(j)}\right]\!\right] = 0$ for $i \neq j$,

$$S_t^{(ij)} = \int_0^1 \tilde{W}_\tau^{(i)} \, \mathrm{d}\tilde{W}_\tau^{(j)} - \int_0^1 \tilde{W}_\tau^{(j)} \, \mathrm{d}\tilde{W}_\tau^{(i)}. \tag{2.38}$$

Consequently, if we define the process $\bar F$ as

$$\bar F_t^{(ij)} := \frac{\partial \bar f^{(i)}}{\partial x^{(j)}}(t, \tilde W_t, \omega),$$

then its quadratic covariation process with the area process $S$ is

$$\llbracket \bar F^{(ij)}, S^{(ij)} \rrbracket_t = \int_0^t \frac{\partial^2 \bar f^{(i)}}{\partial x^{(j)} \partial x^{(j)}}(\tau, \tilde W_\tau, \omega) \tilde W_\tau^{(i)} \,\mathrm{d}\tau \\ - \int_0^t \frac{\partial^2 \bar f^{(i)}}{\partial x^{(j)} \partial x^{(i)}}(\tau, \tilde W_\tau, \omega) \tilde W_\tau^{(j)} \,\mathrm{d}\tau.$$

Converting the Stratonovich integral on the right-hand side of (2.37) to its Itō form we then obtain (2.36). □

**Lemma 2.18.** *If Assumption 2.8 is satisfied, then for all* $c \in \mathbb{R}$, $i, j, k \in \{1, \dots, n\}$, $x \in \mathcal{W}_n^2$ *with* $x(0) \in \operatorname{supp}(P_{X_0})$, $z_0 \in \operatorname{supp}(P_{Z_0})$, *and* $\theta \in \operatorname{supp}(P_\Theta)$,

$$\limsup_{\epsilon \downarrow 0} \tilde{\mathrm{E}}\Big[\exp\Big(c \bar f(1, \tilde W_1, \omega)^\intercal \tilde W_1 - c \bar f(0, \tilde W_0, \omega)^\intercal \tilde W_0\Big) \,\Big|\, \mathbb{B}_\epsilon\Big] \le 1, \tag{2.39a}$$

$$\limsup_{\epsilon \downarrow 0} \tilde{\mathrm{E}}\left[\exp\left(c \int_0^1 \frac{\partial \bar f}{\partial t}(t, \tilde W_t, \omega)^\intercal \tilde W_t \,\mathrm{d}t\right) \,\middle|\, \mathbb{B}_\epsilon\right] \le 1, \tag{2.39b}$$

$$\limsup_{\epsilon \downarrow 0} \tilde{\mathrm{E}}\left[\exp\left(c \int_0^1 \frac{\partial^2 \bar f^{(i)}}{\partial x^{(j)} \partial x^{(k)}}(t, \tilde W_t, \omega) \tilde W_t^{(j)} \,\mathrm{d}t\right) \,\middle|\, \mathbb{B}_\epsilon\right] \le 1. \tag{2.39c}$$

*Proof.* Lemma 2.11 and the definition of the $\epsilon$-ball (2.9) imply that, for all $\epsilon \le 1$, the arguments of $f$ and its derivatives in (2.39) are bounded for almost all $\omega \in \mathbb{B}_\epsilon$. From continuity it then follows that $\bar f$ and its derivatives are bounded for almost all $\omega \in \mathbb{B}_\epsilon$. Consequently, (2.39) then holds as $|||\tilde W||| \to 0$. □

**Lemma 2.19.** *If Assumption 2.8 is satisfied, we then have that for all* $c \in \mathbb{R}$, $i, j \in \{1, \dots, n\}$, $x \in \mathcal{W}_n^2$ *with* $x(0) \in \operatorname{supp}(P_{X_0})$, $z_0 \in \operatorname{supp}(P_{Z_0})$, *and* $\theta \in \operatorname{supp}(P_\Theta)$,

$$\limsup_{\epsilon \downarrow 0} \tilde{\mathrm{E}}\left[\exp\left(c \int_0^1 \frac{\partial \bar f^{(i)}}{\partial x^{(j)}}(t, \tilde W_t, \omega) \,\mathrm{d}S_t^{(ij)}\right) \,\middle|\, \mathbb{B}_\epsilon\right] \le 1. \tag{2.40}$$

*Proof.* Define the $\mathbb{R}^3$-valued process $B$ by

$$B_t^{(1)} := |\tilde W_0| + \int_0^t \frac{\tilde W_\tau^\intercal}{|\tilde W_\tau|} \,\mathrm{d}\tilde W_\tau, \qquad B_t^{(2)} := |Z_0|, \qquad B_t^{(3)} := |\Theta|. \tag{2.41}$$

Using the Itō isometry, we then have that

$$\llbracket B^{(1)} \rrbracket_t = \sum_{i=0}^n \int_0^t \frac{\tilde W_t^{(i)2}}{|\tilde W_t|^2} \,\mathrm{d}t = t,$$

implying that $B^{(1)} - B_0^{(1)}$ is a Wiener process, by Lévy's characterization. Consequently, $B$ is a square-integrable martingale with the predictable representation property adapted with respect to the filtration $\{\mathcal{E}_t\}_{t\geq 0}$. Furthermore, by Itō's lemma we have that

$$\mathrm{d}|\tilde{W}_t|^2 = n\,\mathrm{d}t + 2|\tilde{W}_t|\,\mathrm{d}B_t^{(1)},$$

implying that the random variable $|||X - x|||_{\boldsymbol{G}} = |||\tilde{W}|||$ and the event $\mathbb{B}_\epsilon$ are measurable with respect to the $\sigma$-algebra generated by $B$.

Next, from (2.38) and (2.41) note that

$$\mathrm{d}B_t^{(1)}\,\mathrm{d}S_t^{(ij)} = \frac{\tilde{W}_t^{(i)}\tilde{W}_t^{(j)}\,\mathrm{d}\tilde{W}_t^{(j)}\,\mathrm{d}\tilde{W}_t^{(j)}}{|\tilde{W}_t|} - \frac{\tilde{W}_t^{(j)}\tilde{W}_t^{(i)}\,\mathrm{d}\tilde{W}_t^{(i)}\,\mathrm{d}\tilde{W}_t^{(i)}}{|\tilde{W}_t|} = 0,$$

implying that $[\![ B^{(1)}, \mathrm{d}S^{(ij)} ]\!] = 0$, i.e., the martingales are orthogonal. Additionally, the quadratic variation of $S$ satisfies

$$[\![ S^{(ij)} ]\!]_t = \int_0^t \left[\tilde{W}_\tau^{(i)2} + \tilde{W}_\tau^{(j)2}\right] \mathrm{d}\tau \leq \epsilon^2 \qquad \text{for all } \omega \in \mathbb{B}_\epsilon$$

and from continuity of the derivatives of $f$ we have that the integrand of (2.40) is bounded for almost all $\omega \in \mathbb{B}_\epsilon$. Applying Lemma A.19 we then have that

$$\begin{aligned}\tilde{\mathrm{E}}&\left[\exp\left(c\int_0^1 \frac{\partial \bar{f}^{(i)}}{\partial x^{(j)}}(t, \tilde{W}_t, \omega)\,\mathrm{d}S_t^{(ij)}\right)\,\middle|\,\mathbb{B}_\epsilon\right] \\ &\leq \left(\tilde{\mathrm{E}}\left[\exp\left(2c\int_0^1 \frac{\partial \bar{f}^{(i)}}{\partial x^{(j)}}(t, \tilde{W}_t, \omega)^2\,\mathrm{d}\,[\![ S^{(ij)} ]\!]_t\right)\,\middle|\,\mathbb{B}_\epsilon\right]\right)^{1/2}. \end{aligned} \tag{2.42}$$

Consequently, the integral on the right-hand side of (2.42) vanishes as $\epsilon \downarrow 0$ and (2.40) holds. □

We are now ready to prove the main result of this section.

*Proof of Theorem 2.9.* First, note that if either $x(0) \notin \operatorname{supp}(P_{\boldsymbol{X}_0})$, $z_0 \notin \operatorname{supp}(P_{\boldsymbol{Z}_0})$, or $\theta \notin \operatorname{supp}(P_{\boldsymbol{\Theta}})$, that is, if these variables lie outside the prior support of their corresponding random variables, then (2.10) is satisfied trivially, as $\pi(x(0), z_0, \theta) = 0$ and there exists some $\bar{\epsilon} \in \mathbb{R}_{>0}$ such that, for all $\epsilon \in (0, \bar{\epsilon}]$,

$$P(\mathbb{B}_\epsilon) \leq P(|X_0 - x(0)| < \epsilon, |Z_0 - z_0| < \epsilon, |\Theta - \theta| < \epsilon) = 0. \tag{2.43}$$

We now prove that (2.10) holds for all $x(0) \in \operatorname{supp}(P_{\boldsymbol{X}_0})$, $z_0 \in \operatorname{supp}(P_{\boldsymbol{Z}_0})$, and $\theta \in \operatorname{supp}(P_{\boldsymbol{\Theta}})$. From Lemma 2.13 we have that it suffices to prove that

$$\lim_{\epsilon\downarrow 0} \tilde{\mathrm{E}}[M^{-1} \,|\, \mathbb{B}_\epsilon] = \exp\Big(J(z_0, x, \theta)\Big)$$

or, equivalently, that

$$\lim_{\epsilon\downarrow 0}\tilde{\mathrm{E}}\left[\exp\left(\int_0^1 U_t^\intercal\,\mathrm{d}\tilde{W}_t - \frac{1}{2}\int_0^1 |U_t|^2\,\mathrm{d}t - J(z_0,x,\theta)\right)\middle|\,\mathbb{B}_\epsilon\right] = 1. \tag{2.44}$$

Using Lemmas 2.15 and 2.17 we then have that the exponent of (2.44) can be expanded to

$$\begin{aligned}\int_0^1 U_t^\intercal\,\mathrm{d}\tilde{W}_t - \frac{1}{2}\int_0^1 |U_t|^2\,\mathrm{d}t - J(z_0,x,\theta) = \\ -\frac{1}{2}\int_0^1\left(|U_t|^2 - \left|\boldsymbol{G}^{-1}\left[\hat{f}(t)-\dot{x}(t)\right]\right|^2\right)\mathrm{d}t \\ -\frac{1}{2}\int_0^1 \operatorname{tr}\left[\nabla_{\mathrm{x}} f(t,X_t,Z_t,\Theta) - \nabla_{\mathrm{x}} f(t,x(t),z(t),\theta)\right]\mathrm{d}t \\ -\int_0^1 [\boldsymbol{G}^{-1}\dot{x}(t)]^\intercal\,\mathrm{d}\tilde{W}_t - \bar{f}(1,\tilde{W}_1,\omega)^\intercal\tilde{W}_1 - \bar{f}(0,\tilde{W}_0,\omega)^\intercal\tilde{W}_0 \\ -\int_0^1 \frac{\partial\bar{f}}{\partial t}(t,\tilde{W}_t,\omega)^\intercal\tilde{W}_t\,\mathrm{d}t - \sum_{i,j=1}^n\int_0^1 \frac{\partial\bar{f}^{(i)}}{\partial x^{(j)}}(t,\tilde{W}_t,\omega)\,\mathrm{d}S_t^{(ij)} \\ -\frac{1}{2}\sum_{i,j=1}^n\int_0^1\left(\frac{\partial^2\bar{f}^{(i)}}{\partial x^{(j)}\partial x^{(j)}}(t,\tilde{W}_t,\omega)\tilde{W}^{(i)} - \frac{\partial^2\bar{f}^{(i)}}{\partial x^{(j)}\partial x^{(i)}}(t,\tilde{W}_t,\omega)\tilde{W}^{(j)}\right)\mathrm{d}t.\end{aligned} \tag{2.45}$$

Using Lemma A.13 we then have that (2.44) will hold if (A.13) holds for every element of the right-hand side of (2.45). For the Itō integral of $\dot{x}$, we have that from the theorem of Shepp and Zeitouni (1992),

$$\lim_{\epsilon\downarrow 0}\tilde{\mathrm{E}}\left[\exp\left(c\int_0^1 [\boldsymbol{G}^{-1}\dot{x}(t)]^\intercal\,\mathrm{d}\tilde{W}_t\right)\middle|\,\mathbb{B}_\epsilon\right] = 1 \qquad \text{for all } c\in\mathbb{R}.$$

For the remaining terms of (2.45), we have that the conditional exponential moments vanish from Lemmas 2.14, 2.16, 2.18, and 2.19. □

For paths around $\hat{x}\notin\mathcal{W}_n^2$, outside the Cameron–Martin space, we assume that the fictitious density is null. This assumption was explicitly made in a derivative work of this thesis (Dutra et al., 2014) and is implicit whenever the Onsager–Machlup functional is used for MAP estimation (Aihara and Bagchi, 1999a,b; Zeitouni and Dembo, 1987), as the search space is, at most, the Cameron–Martin space. This is formalized in the conjecture below. We also note that a draft proof of this conjecture was done by George Lowther in the research mathematics forum MathOverflow.[5]

[5] See `http://mathoverflow.net/q/160599`.

**Conjecture 2.20.** *For all $\hat{x} \notin \mathcal{W}_n^2$, $\hat{z}_0 \in \mathbb{R}^q$, and $\hat{\theta} \in \mathbb{R}^m$,*

$$\lim_{\epsilon \downarrow 0} \frac{P(\mathbb{C}_\epsilon)}{a_1 \exp(-\frac{a_2}{\epsilon^2})\, \epsilon^{m+n+q}} = 0$$

*where $a_1, a_2 \in \mathbb{R}_{>0}$ are the constants of Theorem 2.9 and $\mathbb{C}_\epsilon$ is the $\epsilon$-ball centered in $\hat{x}$, $\hat{z}_0$, $\hat{\theta}$, as defined below:*

$$\mathbb{C}_\epsilon := \left\{ \omega \in \Omega \,\middle|\, \|\!|X(\omega) - \hat{x}|\!\|_{\boldsymbol{G}} < \epsilon, |Z(\omega) - \hat{z}_0| < \epsilon, |\Theta(\omega) - \hat{\theta}| < \epsilon \right\}.$$

Finally, we note that a wide range of systems of interest do not satisfy Assumptions 2.8d–c, but it still makes sense applying the joint MAP state path and parameter estimator to them. This is because the model functions are only valid on a small envelope, and the system dynamics looses its meaning outside of it. A rigid body airplane model, for example, is only valid for small accelerations. For large accelerations flexibility effects become apparent and the rigid-body model is not representative. Larger accelerations still would cause structural failure. Consequently, the model function must only have the specified form inside the validity envelope and the regularization of Remark A.11 of Appendix A can be used to ensure the model satisfies the required regularity conditions.

### 2.2.2 Joint posterior fictitious density of state paths and parameters

We now show how the joint *prior* fictitious density derived in the previous subsection can be used to construct the joint *posterior* fictitious density and, with it, the MAP estimator. We begin by presenting the measurement model and its assumptions in a general abstract setting which covers many widespread use cases.

**Assumption 2.21** (measurements)**.**

*a. The $\mathcal{Y}$-valued random variable $Y$ is observed, where $\mathcal{Y}$ is a metric space.*

*b. For all $x \in \mathcal{C}(\mathcal{T}, \mathbb{R}^n)$, $z \in \mathcal{C}(\mathcal{T}, \mathbb{R}^q)$, and $\theta \in \mathbb{R}^m$, the conditional probability measure induced by $Y$, conditioned on $X = x$, $Z = z$ and $\Theta = \theta$, is absolutely continuous with respect to a measure $\nu$ over $(\mathcal{Y}, \mathcal{B}_\mathcal{Y})$ and admits a conditional density $\psi$ with respect to it, i.e., for all $\mathbb{B} \in \mathcal{B}_\mathcal{Y}$*

$$P(Y \in \mathbb{B} \,|\, X = x, Z = z, \Theta = \theta) = \int_{\mathbb{B}} \psi(y \,|\, x, z, \theta)\, \mathrm{d}\nu(y). \tag{2.46}$$

*c. For the observed value $y \in \mathcal{Y}$ of $Y$, the likelihood $\psi$ is continuous with respect to its second argument $x$ (with respect to the supremum norm) and with respect to its third argument $\theta$.*

*d. For the observed value $y \in \mathcal{Y}$ of $Y$,* $\mathrm{E}[\psi(y \mid X, Z, \varTheta)] > 0$.

Assumption 2.21 and its underlying representation cover many measurement models and distributions. For $\mathcal{Y} = \mathbb{R}^N$ and $\nu$ as the Lebesgue measure, for example, $Y$ can be a continuous random variable with $\psi$ its conditional probability density function, as examplified in Sections 4.1.1 and 4.1.2. Similarly, $\nu$ can be the counting measure, $Y$ a discrete random variable and $\psi$ its probability *mass* function, as examplified in Section 4.1.3. To represent continuous-time measurements when $Y$ is a diffusion depending on $X$, $Z$, and $\varTheta$, $\nu$ can be the Gaussian measure over $\mathcal{C}(\mathcal{T}, \mathbb{R}^d)$ of the driving process and $\psi$ is given by variations of the Kallianpur–Striebel formula (Kallianpur and Striebel, 1968, 1969; Kallianpur, 1980, Sec. 11.3; van Handel, 2007, Lem. 1.1.5).

Using these assumptions, the joint *posterior* fictitious density of $X$, $Z_0$ and $\varTheta$ is then given by the theorem below.

**Theorem 2.22** (joint posterior fictitious density of state paths and parameters)**.** *If a system of the form of* (2.4) *satisfies Assumption 2.8 and has measurements satisfying Assumption 2.21, then there exist $a_1, a_2 \in \mathbb{R}_{>0}$ such that for all $x \in \mathcal{W}^2_n$, $z_0 \in \mathbb{R}^q$, and $\theta \in \mathbb{R}^m$,*

$$\lim_{\epsilon \downarrow 0} \frac{P(\mathbb{B}_\epsilon \mid Y = y)}{a_1 \exp(-\frac{a_2}{\epsilon^2})\, \epsilon^{m+n+q}} = \frac{\psi(y \mid x, z, \theta)\, \pi(x(0), z_0, \theta) \exp\Big(J(z_0, x, \theta)\Big)}{\mathrm{E}[\psi(y \mid X, Z, \varTheta)]}, \tag{2.47}$$

*where $\mathbb{B}_\epsilon$ is the $\epsilon$-ball centered in $x$, $z_0$, and $\theta$, as defined in* (2.9)*; $z \in \mathcal{W}^2_q$ is the solution to the initial value problem* (2.12)*; and $J$ is the Onsager–Machlup functional defined in* (2.11)*.*

*Proof.* First, note that (2.47) is satisfied trivially for $x(0) \notin \operatorname{supp}(P_{X_0})$, $z_0 \notin \operatorname{supp}(P_{Z_0})$, or $\theta \notin \operatorname{supp}(P_\varTheta)$, as $\pi(x(0), z_0, \theta) = 0$ and there exists some $\bar\epsilon \in \mathbb{R}_{>0}$ such that (2.43) holds for all $\epsilon \in (0, \bar\epsilon]$.

Next, assume that $x(0) \in \operatorname{supp}(P_{X_0})$, $z_0 \in \operatorname{supp}(P_{Z_0})$, and $\theta \in \operatorname{supp}(P_\varTheta)$. Using the measure-theoretic formulation of Bayes' theorem (Schervish, 1995, Thm. 1.31), we have that

$$P(\mathbb{B}_\epsilon \mid Y = y) = \frac{\int_{\mathbb{B}_\epsilon} \psi(y \mid X(\omega), Z(\omega), \varTheta(\omega))\, \mathrm{d}P(\omega)}{\mathrm{E}[\psi(y \mid X, Z, \varTheta)]}. \tag{2.48}$$

As Lemma 2.13 guarantees that $P(\mathbb{B}_\epsilon) > 0$, (2.48) can be simplified to

$$P(\mathbb{B}_\epsilon \mid Y = y) = \frac{\mathrm{E}[\psi(y \mid X, Z, \varTheta) \mid \mathbb{B}_\epsilon]\, P(\mathbb{B}_\epsilon)}{\mathrm{E}[\psi(y \mid X, Z, \varTheta)]}.$$

Due to the continuity of $\psi$ (Assumption 2.21c),

$$\lim_{\epsilon \downarrow 0} \mathrm{E}[\psi(y \mid X, Z, \varTheta) \mid \mathbb{B}_\epsilon] = \psi(y \mid x, z, \theta)\,,$$

which together with Theorem 2.9 implies that (2.47) holds. □

From Definition 2.7 and Theorem 2.22, we have that the joint MAP estimator of $X$, $Z_0$ and $\Theta$ is obtained by maximizing posterior joint fictitious density, the left-hand side of (2.47), subject to having $z$ satisfy the initial value problem (2.12). For analysis and implementation of the estimator, however, it is more tractable to work with the logarithm of fictitious posterior density, which is known as the *log-posterior* for short. In addition, the denominator can be dropped since it is constant and does not influence the location of maxima.

The joint MAP state-path and parameter estimator is then the solution to the following optimization problem:

$$\begin{aligned} \underset{x\in\mathcal{W}_n^2,\; z\in\mathcal{W}_q^2,\; \theta\in\mathbb{R}^m}{\text{maximize}} \quad & \ell(x,z,\theta) \\ \text{subject to} \quad & \dot z(t) = h\big(t, x(t), z(t), \theta\big)\,, \end{aligned}$$

where the log-posterior $\ell\colon \mathcal{W}_n^2 \times \mathcal{W}_q^2 \times \mathbb{R}^m \to \overline{\mathbb{R}}$ is given by

$$\begin{aligned} \ell(x,z,\theta) := \ln\psi(y\,|\,x,z,\theta) + \ln\pi(x(0),z(0),\theta) \\ -\frac{1}{2}\int_{\mathcal{T}} \left|\boldsymbol{G}^{-1}\left[f(t,x(t),z(t),\theta) - \dot x(t)\right]\right|^2 \mathrm{d}t \\ -\frac{1}{2}\int_{\mathcal{T}} \operatorname{div}_{\mathrm{x}} f\big(t,x(t),z(t),\theta\big)\,\mathrm{d}t. \end{aligned} \tag{2.49}$$

Alternatively, the search space can be expanded and more constraints added to obtain the following equivalent optimization problem:

$$\begin{aligned} \underset{x,w\in\mathcal{W}_n^2,\; z\in\mathcal{W}_q^2,\; \theta\in\mathbb{R}^m}{\text{maximize}} \quad & \ell_{\mathrm{a}}(x,z,\theta,w) \\ \text{subject to} \quad & \dot x(t) = f\big(t,x(t),z(t),\theta\big) + \boldsymbol{G}\dot w(t), \\ & \dot z(t) = h\big(t,x(t),z(t),\theta\big) \end{aligned} \tag{2.50}$$

where the alternative log-posterior $\ell_{\mathrm{a}}: \mathcal{W}_n^2 \times \mathcal{W}_q^2 \times \mathbb{R}^m \times \mathcal{W}_n^2 \to \overline{\mathbb{R}}$ is given by

$$\begin{aligned} \ell_{\mathrm{a}}(x,z,\theta,w) := \ln\psi(y\,|\,x,z,\theta) + \ln\pi(x(0),z(0),\theta) - \tfrac{1}{2}\|w\|^2_{\mathcal{W}_n^2} \\ -\frac{1}{2}\int_{\mathcal{T}} \operatorname{div}_{\mathrm{x}} f\big(t,x(t),z(t),\theta\big)\,\mathrm{d}t. \end{aligned} \tag{2.51}$$

The optimization in the form (2.50) is analogous to an optimal control problem and is more tractable for numerical implementation and some theoretical analyses.

## 2.3 Minimum-energy state path and parameter estimation

Throughout the literature, a commonly used alternative to the MAP state-path estimator is obtained by using a functional equal to the log-posterior without the drift divergence term (Bryson and Frazier, 1963; Cox, 1963, Chap. III; Mortensen, 1968; Jazwinski, 1970, p. 155). This was denoted the *minimum energy estimator* by Hijab (1980), since its estimates minimize the input energy of the associated control system.

For many years, the minimum energy estimator was believed to be a good MAP estimator (Zeitouni and Dembo, 1987, p. 234). Hijab (1984) showed that, for systems with continuous-time measurements under Gaussian noise, the minimum energy merit is the limit optimal smoother when the intensity of both the process and measurement noise vanishes. In a derivative work of this thesis (Dutra et al., 2014), we proved that the minimum energy estimates are state paths corresponding to MAP noise paths. In this section, we extend our previous work by also covering the case where there are unknown parameters and not all states are under direct influence of noise.

We begin by stating the fictitious density of Wiener processes with the following theorem, a variant of Theorem 2.9 with no drift and a fixed initial state.

**Theorem 2.23** (Fujita and Kotani, 1982, Capitaine, 1995)**.** *Let $W$ be an $n$-dimensional Wiener process. Then there exists $a_{11}, a_{12} \in \mathbb{R}_{>0}$ such that, for all $w \in \mathcal{W}_n^2$ with $w(0) = 0$,*

$$\lim_{\epsilon \downarrow 0} \frac{P(|||W - w||| < \epsilon)}{a_{11} \exp(-\frac{a_{12}}{\epsilon^2})} = \exp\Big(-\tfrac{1}{2}\|\dot{w}\|^2_{L^2_n}\Big).$$

We are now ready to derive the prior joint fictitious density of $W$, $X_0$, $Z_0$, and $\Theta$. We will denote by $\mathbb{B}^{\mathrm{n}}_\epsilon \in \mathcal{E}$ the event that $W$, $X_0$, $Z_0$ and $\Theta$ are inside an $\epsilon$-ball centered in some $w \in \mathcal{C}(\mathcal{T}, \mathbb{R}^n)$, $x_0 \in \mathbb{R}^n$, $z_0 \in \mathbb{R}^q$ and $\theta \in \mathbb{R}^m$, for $\epsilon \in \mathbb{R}_{>0}$:

$$\mathbb{B}^{\mathrm{n}}_\epsilon := \Big\{\omega \in \Omega \,\Big|\, |||w - W(\omega)||| < \epsilon, |X_0(\omega) - x_0| < \epsilon, \\ |Z_0(\omega) - z_0| < \epsilon, |\Theta(\omega) - \theta| < \epsilon\Big\}. \tag{2.52}$$

**Theorem 2.24.** *If Assumption 2.8 is satisfied, then there exist $a_{11}, a_{12} \in \mathbb{R}_{>0}$ such that, for all $w \in \mathcal{W}_n^2$ with $w(0) = 0$, $x_0 \in \mathbb{R}^n$, $z_0 \in \mathbb{R}^q$, and $\theta \in \mathbb{R}^m$,*

$$\lim_{\epsilon \downarrow 0} \frac{P(\mathbb{B}^{\mathrm{n}}_\epsilon)}{a_{11} \exp(-\frac{a_{12}}{\epsilon^2})\, \epsilon^{m+n+q}} = \pi(x_0, z_0, \theta) \exp\Big(-\tfrac{1}{2}\|\dot{w}\|^2_{L^2_n}\Big), \tag{2.53}$$

*where $\mathbb{B}^{\mathrm{n}}_\epsilon$ is the $\epsilon$-ball defined in* (2.52)*.*

*Proof.* From Assumption 2.8a we have that $X_0$, $Z_0$, and $\Theta$ are $\mathcal{E}_0$-measurable and, consequently, independent of $W$. Applying the Lebesgue differentiation theorem and Theorem 2.23 we then have that (2.53) holds. □

Next, to calculate the fictitious joint posterior density of noise paths, we show that for each $w$, $x_0$, $z_0$ and $\theta$ there exists an associated state path such that its largest distance from $X$, for all $\omega \in \mathbb{B}^{\mathrm{n}}_\epsilon$, vanishes with $\epsilon$. This lemma is similar to Lemma 2.11.

**Lemma 2.25** (associated state path)**.** *If Assumption 2.8 is satisfied then, for all for all $w \in \mathcal{W}^2_n$ with $w(0) = 0$, $x_0 \in \mathbb{R}^n$, $z_0 \in \mathbb{R}^q$, and $\theta \in \operatorname{supp}(P_\Theta)$,*

$$\lim_{\epsilon \downarrow 0} \sup_{\omega \in \mathbb{B}^{\mathrm{n}}_\epsilon} |||X(\omega) - x||| = 0, \qquad \lim_{\epsilon \downarrow 0} \sup_{\omega \in \mathbb{B}^{\mathrm{n}}_\epsilon} |||Z(\omega) - z||| = 0,$$

*where $\mathbb{B}^{\mathrm{n}}_\epsilon$ is the $\epsilon$-ball defined in* (2.52) *and $x \in \mathcal{W}^2_n$ and $z \in \mathcal{W}^2_q$ are the unique solutions to the following integral equations, for all $t \in \mathcal{T}$:*

$$x(t) = x_0 + \int_0^t f(\tau, x(\tau), z(\tau), \theta)\, \mathrm{d}\tau + \boldsymbol{G} w(t) \tag{2.54a}$$

$$z(t) = z_0 + \int_0^t h(\tau, x(\tau), z(\tau), \theta)\, \mathrm{d}\tau. \tag{2.54b}$$

*Proof.* Taking the difference between the integral representation of $x$ and $z$ in (2.54) and $X$ and $Z$ in (2.4) we have that

$$X_t - x(t) = X_0 - x_0 + \int_0^t \left[F_\tau - \hat{f}(\tau)\right] \mathrm{d}\tau + \boldsymbol{G}[W_t - w(t)],$$
$$Z_t - z(t) = Z_0 - z_0 + \int_0^t \left[H_\tau - \hat{h}(\tau)\right] \mathrm{d}\tau,$$

where the helping functions and processes $\hat{f}$, $\hat{h}$, $F$, and $H$ of (2.13) and (2.14) were used. Taking the norm on both sides an applying the triangle inequality we have that, for all $\omega \in \mathbb{B}^{\mathrm{n}}_\epsilon$,

$$|X_t - x(t)| < \epsilon + \int_0^t \left|F_\tau - \hat{f}(\tau)\right| \mathrm{d}\tau + \epsilon \left\|\boldsymbol{G}\right\|,$$
$$|Z_t - z(t)| < \epsilon + \int_0^t \left|H_\tau - \hat{h}(\tau)\right| \mathrm{d}\tau,$$

where the matrix norm of $\|\boldsymbol{G}\|$ is the one induced by the Euclidean norm.

Next, note that

$$\begin{aligned} \left|F_t - \hat{f}(t)\right| &= \left|F_t - f(t, X_t, Z_t, \theta) + f(t, X_t, Z_t\theta) - \hat{f}(t)\right|, \\ &\leq |F_t - f(t, X_t, Z_t, \theta)| + \left|f(t, X_t, Z_t, \theta) - \hat{f}(t)\right|, && (2.55a) \\ &\leq \rho_f(|\Theta - \theta|) + L^\theta_{\mathrm{f}} |X_t - x(t)| + L^\theta_{\mathrm{f}} |Z_t - z(t)|, && (2.55b) \end{aligned}$$

where (2.55a) is obtained by using the triangle inequality and (2.55b) by using the uniform and Lipschitz continuity of $f$ from Assumptions 2.8c–d. Similarly, for the clean states we have that

$$\left|H_t - \hat{h}(t)\right| \leq \rho_h(|\Theta - \theta|) + L_{\mathrm{h}}^{\theta}\,|X_t - x(t)| + L_{\mathrm{h}}^{\theta}\,|Z_t - z(t)|\,.$$

Consequently, we have that for all $\omega \in \mathbb{B}_\epsilon^{\mathrm{n}}$,

$$\begin{aligned}|X_t - x(t)| + |Z_t - z(t)| < 2\epsilon + \epsilon\,\|\boldsymbol{G}\| + \rho_f(\epsilon) + \rho_h(\epsilon)\\ + (L_{\mathrm{f}}^{\theta} + L_{\mathrm{h}}^{\theta})\int_0^t (|X_\tau - x(\tau)| + |Z_\tau - z(\tau)|)\,\mathrm{d}\tau.\end{aligned}$$

Applying the Grönwall–Bellman inequality, Lemma A.8, we then have that

$$\begin{aligned}\sup_{\omega\in\mathbb{B}_\epsilon^{\mathrm{n}}} \|\!|X(\omega) - x(t)|\!\| + \|\!|Z(\omega) - z(t)|\!\|\\ \leq [2\epsilon + \epsilon\,\|\boldsymbol{G}\| + \rho_f(\epsilon) + \rho_h(\epsilon)]\exp(L_{\mathrm{f}}^{\theta} + L_{\mathrm{h}}^{\theta}).\end{aligned}$$ □

Since the state path, given $\mathbb{B}_\epsilon^{\mathrm{n}}$, converges to the associated state path, the posterior fictitious density of $X_0$, $Z_0$ $\Theta$ and $W$ is simply the product of the prior fictitious density and the likelihood evaluated at the associated state path. The statement and proof of the following theorem are analogous to those of Theorem 2.22.

**Theorem 2.26** (joint posterior fictitious density of noise paths and parameters)**.** *If Assumptions 2.8 and 2.21 are satisfied, then there exist* $a_{11}, a_{12} \in \mathbb{R}_{>0}$ *such that, for all* $w \in \mathcal{W}_n^2$ *with* $w(0) = 0$*,* $x_0 \in \mathbb{R}^n$*,* $z_0 \in \mathbb{R}^q$*, and* $\theta \in \mathbb{R}^m$*,*

$$\lim_{\epsilon\downarrow 0} \frac{P(\mathbb{B}_\epsilon^{\mathrm{n}} \,|\, Y = y)}{a_{11}\exp\!\left(-\frac{a_{12}}{\epsilon^2}\right)\epsilon^{m+n+q}} = \psi(y \,|\, x, z, \theta)\,\pi(x_0, z_0, \theta)\exp\!\Big(-\tfrac{1}{2}\|\dot{w}\|_{L_n^2}^2\Big)\,, \tag{2.56}$$

*where* $\mathbb{B}_\epsilon^{\mathrm{n}}$ *is the* $\epsilon$*-ball defined in* (2.52) *and* $x$ *and* $z$ *are the associated noise paths satisfying the integral equations* (2.54)*.*

*Proof.* As in the proof of Theorem 2.22, we assume that $x_0 \in \operatorname{supp}(P_{X_0})$, $z_0 \in \operatorname{supp}(P_{Z_0})$, and $\theta \in \operatorname{supp}(P_\Theta)$ as (2.56) is satisfied trivially otherwise. Then, using the measure-theoretic formulation of Bayes' theorem (Schervish, 1995, Thm. 1.31), we have that

$$P(\mathbb{B}_\epsilon^{\mathrm{n}} \,|\, Y = y) = \frac{\int_{\mathbb{B}_\epsilon^{\mathrm{n}}} \psi(y \,|\, X(\omega), Z(\omega), \Theta(\omega))\,\mathrm{d}P(\omega)}{\mathrm{E}[\psi(y \,|\, X, Z, \Theta)]}. \tag{2.57}$$

As $P(\mathbb{B}_\epsilon) > 0$, (2.57) can be simplified to

$$P(\mathbb{B}_\epsilon^{\mathrm{n}} \,|\, Y = y) = \frac{\mathrm{E}[\psi(y \,|\, X, Z, \Theta) \,|\, \mathbb{B}_\epsilon]\,P(\mathbb{B}_\epsilon^{\mathrm{n}})}{\mathrm{E}[\psi(y \,|\, X, Z, \Theta)]}.$$

From Lemma 2.25 and the continuity of $\psi$ (Assumption 2.21c),

$$\lim_{\epsilon\downarrow 0} \mathrm{E}[\psi(y \,|\, X, Z, \Theta) \,|\, \mathbb{B}^{\mathrm{n}}_{\epsilon}] = \psi(y \,|\, x, z, \theta)\,,$$

which together with Theorem 2.24 implies that (2.56) holds. □

Having an expression for the the joint posterior fictitious density, we can take its logarithm and construct the minimum energy estimation problem:

$$\begin{aligned} \underset{x,w\in\mathcal{W}^2_n,\ z\in\mathcal{W}^2_q,\ \theta\in\mathbb{R}^m}{\text{maximize}} \quad & \ln\psi(y \,|\, x, z, \theta) + \ln\pi(x(0), z(0), \theta) - \tfrac{1}{2}\|w\|_{\mathcal{W}^2_n} \\ \text{subject to} \quad & \dot{x}(t) = f(t, x(t), z(t), \theta) + \boldsymbol{G}\dot{w}(t), \\ & \dot{z}(t) = h(t, x(t), z(t), \theta)\,. \end{aligned} \tag{2.58}$$

Note that for any $w$ satisfying the constraints,

$$\dot{w}(t) = \boldsymbol{G}^{-1}\left[\dot{x}(t) - f(t, x(t), z(t), \theta)\right],$$

implying that (2.58) is equivalent to

$$\begin{aligned} \underset{x\in\mathcal{W}^2_n,\ z\in\mathcal{W}^2_q,\ \theta\in\mathbb{R}^m}{\text{maximize}} \quad & \ell_{\mathrm{e}}(x, z, \theta) \\ \text{subject to} \quad & \dot{z}(t) = h(t, x(t), z(t), \theta)\,, \end{aligned} \tag{2.59}$$

where the energy log-posterior $\ell_{\mathrm{e}} \colon \mathcal{W}^2_n \times \mathcal{W}^2_q \times \mathbb{R}^m \to \overline{\mathbb{R}}$ is defined as

$$\begin{aligned} \ell_{\mathrm{e}}(x, z, \theta) := {} & \ln\psi(y \,|\, x, z, \theta) + \ln\pi(x(0), z(0), \theta) \\ & - \frac{1}{2}\int_{\mathcal{T}} \left|\boldsymbol{G}^{-1}\left[\dot{x}(t) - f(t, x(t), z(t), \theta)\right]\right|^2 \,\mathrm{d}t. \end{aligned} \tag{2.60}$$

# Chapter 3

# MAP estimation in discretized SDEs

In the context of discrete-time dynamical systems, maximum *a posteriori* (MAP) state-path estimation has become viable in the recent years due to the availability of efficient software packages for nonlinear optimization of large-scale problems. It can be applied to nonlinear discrete-time systems with general initial, transition and measurement densities, giving it a larger applicability than nonlinear Kalman smoothers and filters.[1] Furthermere, the possibility to use heavy-tailed measurement distributions such as Student's $t$ and the $\ell_1$-Laplace with these estimators lends them robustness against outlier measurements (Aravkin et al., 2011, 2012a,b,c, 2013; Dutra et al., 2014; Farahmand et al., 2011).

Many systems and phenomena of engineering interest are continuous-time in nature and can be modelled by stochastic differential equations. To apply discrete-time MAP state-path estimation to these systems, their dynamics first needs to be discretized. However, while the discretization of linear systems can be exact, for general nonlinear SDEs it is necessary to employ approximations which improve as the discretization step vanishes; see Kloeden and Platen (1992) for a thourough review of many such discretization schemes. Applications of MAP state-path estimation to discretized systems are presented by Bell et al. (2009) and Aravkin et al. (2011, 2012c).

Nevertheless, the estimates of the discretized state-paths are only meaningful if they have some statistical interpretation under the original continuous-time model. In particular one would expect that, as the discretization step decreases and the discretization improves, the discretized MAP state-paths converge to the MAP state-path of the continuous-time system, as defined in Chapter 2. In a derivative work of this thesis (Dutra et al., 2014) we proved that, under some regularity conditions, MAP estimators discretized with both the Euler and trapezoidal schemes converge hypographically[2] as the discretization step vanishes. However, the hypographical limit of these discretized estimators is not, in general, the same; the Euler-discretized estimator hypo-converges to the minimum energy estimator and the trapezoidally-discretized estimator hypo-converges to the MAP state-path estimator. Hypographical convergence of log-densities is used as the mode of convergence because it has implications on the convergence of the MAP estimates.

[1] In this thesis, we group the robust *Kalman* smoothers of Bell et al. (2009) and Aravkin et al. (2011, 2012a,b,c, 2013) with MAP state-path estimators instead of with nonlinear Kalman smoothers, as they maximize the posterior log-density of the state-paths instead of approximating the posterior mean and covariance of the states with heuristics and then applying Kalman smoothing algorithms (as do the extended and unscented Kalman smoothers, for example).

In this chapter, we extend the results of Dutra et al. (2014) to joint MAP state-path and parameter estimation with possibly singular diffusion matrices and more general measurements. The results are analogous, i.e., the Euler-discretized estimator hypo-converges to the joint minimum energy state-path and parameter estimator and the trapezoidally-discretized estimator hypo-converges to the joint MAP state-path and parameter estimator.

The remainder of this chapter is organized as follows: in Section 3.1 we define hypographical convergence and state its implications. Then, in Sections 3.2 and 3.3 we obtain the hypographical limits of the Euler- and trapezoidally-discretized joint MAP state-path and parameter estimators, respectively.

## 3.1 Hypo-convergence

Hypographical convergence is a mode of convergence of functions which is widely used for the analysis of approximations to optimization functions. Its importance lies in the fact that if a sequence of functions converge hypographically, then any limit point of their sequence of maximizers is a maximizer of the hypographical limit. Hypo-convergence is formally defined as the Kuratowski convergence of hypographs, but many equivalent and easier-to-work-with definitions exist (see Attouch, 1984, Sec. 1.2). The definition in the lemma below, adapted from Attouch (1984, Prop. 1.14), is the one used in this thesis.

**Lemma 3.1** (hypo-convergence)**.** *Let* $(\mathcal{X}, \tau)$ *be a first-countable topological space and* $\{f_i\}_{i=1}^{\infty}$ *be a sequence of functions* $f_i\colon \mathcal{X} \to \overline{\mathrm{I\!R}}$*, where* $\overline{\mathrm{I\!R}}$ *is the extended real line. Then* $f_i$ *is said to converge hypographically to a function* $f\colon \mathcal{X} \to \overline{\mathrm{I\!R}}$ *if and only if*

*a. for every convergent sequece* $\{x_i\}_{i=1}^{\infty}$ *of* $x_i \in \mathcal{X}$

$$\limsup_{i\to\infty} f_i(x_i) \leq f(x), \qquad \text{where } x := \lim_{i\to\infty} x_i;$$

[2] In a slight abuse of notation, we will say that MAP *estimators* converge hypographically if their *log-posteriors* do so.

*b. for every $x \in \mathcal{X}$ there exists a convergent sequence $\{x_i\}_{i=1}^{\infty}$ of $x_i \in \mathcal{X}$ such that*

$$\lim_{i\to\infty} x_i = x, \qquad \liminf_{i\to\infty} f_i(x_i) \geq f(x).$$

The importance of hypo-convergence lies in the following lemma, adapted from Polak (1997, Thm. 3.3.3), which relates the maximizers of a sequence of functions to the maximizer of its hypographical limit.

**Lemma 3.2.** *Let $(\mathcal{X}, \tau)$ be a first-countable space, and $\{f_i\}_{i=1}^{\infty}$ be a sequence of functions $f_i\colon \mathcal{X} \to \overline{\mathbb{R}}$ which converges hypographically to $f\colon \mathcal{X} \to \overline{\mathbb{R}}$. Then, for any convergent sequence $\{x_i^*\}_{i\in\mathcal{N}}$ indexed over $\mathcal{N} \subset \mathbb{N}$ of maximizers $x_i^* \in \mathcal{X}$ of $f_i$, i.e.,*

$$\sup_{x\in\mathcal{X}} f_i(x) = f_i(x_i^*) \qquad \text{for all } i \in \mathcal{N},$$

*we have that their limit point $x^* := \lim_{i\to\infty} x_i$ is a maximizer of $f$, i.e.,*

$$\sup_{x\in\mathcal{X}} f(x) = f(x^*).$$

Next we apply these concepts to the discretized joint MAP state-path and parameter estimation.

## 3.2 Euler-discretized estimator

The stochastic Euler scheme (Kloeden and Platen, 1992, Sec. 9.1), also known as the Euler–Maruyama scheme, is the simplest and one of the most widely used SDE discretization schemes. Throughout this section, we will assume that $f$, $h$, $X$, $Z$, $\Theta$, and $W$ are the same as defined in Section 2.2 and that the system and measurements satisfy Assumptions 2.8 and 2.21. Next, let $\mathcal{P} := \{t_0, \dots, t_N\}$ be a partition of $\mathcal{T}$, with $t_0 = 0$, $t_N = t_\mathrm{f}$, and $t_k < t_{k+1}$. The Euler scheme approximates the SDEs (2.4) at the partition points with the following system of difference equations:

$$\tilde{X}_{t_{k+1}} = \tilde{X}_{t_k} + f(t_k, \tilde{X}_{t_k}, \tilde{Z}_{t_k}, \Theta)\delta_k + \boldsymbol{G}\Delta W_k, \tag{3.1a}$$

$$\tilde{Z}_{t_{k+1}} = \tilde{Z}_{t_k} + h(t_k, \tilde{X}_{t_k}, \tilde{Z}_{t_k}, \Theta)\delta_k, \tag{3.1b}$$

where $\delta_k := t_{k+1} - t_k$ is the time increment, $\Delta W_k := W_{t_{k+1}} - W_{t_k}$ is the Wiener process increment, and the processes $\tilde{X}$ and $\tilde{Z}$ are the approximations to $X$ and $Z$ with $\tilde{X}_0 := X_0$ and $\tilde{Z}_0 := Z_0$.

For the remaining time-points, we consider $\tilde{X}$ to be the piecewise linear interpolation of the $\tilde{X}_{t_k}$, i.e., for all $t \in [t_k, t_{k+1}]$,

$$\tilde{X}_t = \frac{t_{k+1} - t}{\delta_k}\tilde{X}_{t_k} + \frac{t - t_k}{\delta_k}\tilde{X}_{t_{k+1}},$$

$$\tilde{Z}_t = \frac{t_{k+1} - t}{\delta_k}\tilde{Z}_{t_k} + \frac{t - t_k}{\delta_k}\tilde{Z}_{t_{k+1}}.$$

This interpolation is often used in association with the Euler scheme (Kloeden and Platen, 1992, p. 307) and was chosen due to its simplicity and the fact that the resulting functions are absolutely continuous. The results of this subsection hold, nonetheless, for any other absolutely continuous interpolation scheme which is a convex combination of the adjacent interpolation points. In what follows, we will denote by $\mathrm{PL}(\mathcal{P}, \mathbb{R}^n)$ the space of piecewise linear functions from $\mathcal{T}$ to $\mathbb{R}^n$ with breaks over the partition $\mathcal{P}$.

The joint density of $\tilde{X}_{t_0}, \dots, \tilde{X}_{t_N}$, $\tilde{Z}_0$ and $\Theta$, which will be denoted the Euler-discretized joint state-path and parameter density, can be found from the joint density of $\tilde{X}_0$, $\tilde{Z}_0$, $\Theta$ and $\Delta W_0, \dots, \Delta W_{N-1}$ through a change of variables (Lemma A.20), since from (3.1) we have that both groups of variables are related by a diffeomorphism. Recalling the definition of the Wiener process, we have that the $\Delta W_k$ are independent normally distributed random variables with zero mean and variance $\delta_k \boldsymbol{I}_n$. Furthermore, $\Delta W$ is independent of $\tilde{X}_0$, $\tilde{Z}_0$, and $\Theta$. Consequently, the Euler-discretized joint posterior state-path and parameter density is given by

$$
\begin{aligned}
p(\tilde{x}, \tilde{z}_0, \vartheta \,|\, y) = {} & \frac{\psi(y \,|\, \tilde{x}, \tilde{z}, \vartheta)}{\mathrm{E}[\psi(y \,|\, X, Z, \Theta)]} \pi(\tilde{x}(0), \tilde{z}_0, \vartheta) \\
& \times \prod_{k=0}^{N-1} \frac{\exp\Big(-\delta_k \frac{1}{2} \left|\boldsymbol{G}^{-1} \left[\frac{\Delta \tilde{x}_k}{\delta_k} - f(t_k, \tilde{x}(t_k), \tilde{z}(t_k), \vartheta)\right]\right|^2\Big)}{|\det G| \sqrt{\delta_k (2\pi)^n}},
\end{aligned} \tag{3.2}
$$

where $\Delta \tilde{x}_k := \tilde{x}(t_{k+1}) - \tilde{x}(t_k)$ and $\tilde{z} \in \mathrm{PL}(\mathcal{P}, \mathbb{R}^q)$ is the Euler-discretized clean-state path associated with $\tilde{x}$, $\tilde{z}_0$ and $\vartheta$, i.e., $\tilde{z}(0) = \tilde{z}_0$ and

$$
\tilde{z}(t) = \tilde{z}(t_k) + \int_{t_k}^{t} h(t_k, \tilde{x}(t_k), \tilde{z}(t_k), \vartheta) \, \mathrm{d}s \qquad \text{for all } t \in (t_k, t_{k+1}].
$$

It should be noted that any fictitious density of $\tilde{X}$ with respect to the supremum norm (with whatever underlying finite-dimensional norm) evaluated at $\tilde{x} \in \mathrm{PL}(\mathcal{P}, \mathbb{R}^n)$ is proportional to the joint density of $\tilde{X}_{t_0}, \dots, \tilde{X}_{t_N}$, with respect to the Lebesgue measure, evaluated at $\tilde{x}(t_0), \dots, \tilde{x}(t_N)$. Consequently, we may refer to the joint fictitious density of $\tilde{X}$, $\tilde{Z}$, and $\Theta$ and to the joint density of $\tilde{X}_{t_0}, \dots, \tilde{X}_{t_N}$, $\tilde{Z}_0$, and $\Theta$, interchangeably, as the Euler-discretized state-path and parameter density, and analogously to the posterior densities.

Due to a better numerical tractability, it is more convenient to minimize the logarithm of the posterior density instead of the density itself. The logarithm over $\mathbb{R}_{\geq 0}$ is a strictly increasing function and as such does not change the location of maxima. Taking the logarithm of (3.2) and removing the constant terms, which do not change the location of maxima as well, we obtain $\tilde{\ell}\colon \mathrm{PL}(\mathcal{P}, \mathbb{R}^n) \times \mathbb{R}^q \times \mathbb{R}^m \to \overline{\mathbb{R}}$, the Euler-discretized joint state-path and parameter log-posterior:

$$\begin{aligned}\tilde{\ell}(\tilde{x},\tilde{z}_0,\vartheta) := \ln\psi(y\,|\,\tilde{x},\tilde{z},\vartheta) + \ln\pi(\tilde{x}(0),\tilde{z}_0,\vartheta) \\ -\frac{1}{2}\sum_{k=0}^{N-1}\delta_k\left|\boldsymbol{G}^{-1}\left[\tfrac{\Delta\tilde{x}_k}{\delta_k} - f(t_k,\tilde{x}(t_k),\tilde{z}(t_k),\vartheta)\right]\right|^2,\end{aligned} \tag{3.3}$$

where, as in (3.2), $\tilde{z}$ is the Euler-discretized clean-state path associated with $\tilde{x}$, $\tilde{z}_0$ and $\vartheta$.

Comparing the above expression of the Euler-discretized log-posterior (3.3) to the energy merit (2.59) in page 46, a great similarity is apparent. We will now prove that for a sequence of partitions with a vanishing mesh, the Euler-discretized log-posterior converges hypographically to the energy log-posterior.

### 3.2.1 Hypo-convergence of the Euler-discretized log-posterior

Let $\{\mathcal{P}_i\}_{i=1}^{\infty}$ be a sequence of nested partitions $\mathcal{P}_i := \{t_{k,i}\}_{k=0}^{N_i}$ of $\mathcal{T}$ with a vanishing mesh, i.e., $\mathcal{P}_i \subset \mathcal{P}_{i+1}$,

$$\delta_{ki} := t_{k+1,i} - t_{ki}, \qquad \bar{\delta}_i := \max_{0\le k<N_i}\delta_{ki}, \qquad \lim_{i\to\infty}\bar{\delta}_i = 0.$$

The variable $\bar{\delta}_i$ represents the mesh of $\mathcal{P}_i$ and the comma in the variables with double subscript may be dropped if unambiguous ($t_{k,i} = t_{ki}$).

For each $\mathcal{P}_i$, we consider its corresponding Euler-discretized joint state-path and parameter log-posterior $\tilde{\ell}_i\colon \mathcal{W}_n^2(\mathcal{T})\times\mathbb{R}^q\times\mathbb{R}^m \to \overline{\mathbb{R}}$. We extend its domain to non piecewise linear functions by setting its value to negative infinity, i.e., $\tilde{\ell}_i(x,\zeta_i,\vartheta_i) := -\infty$ for $x\notin \mathrm{PL}(\mathcal{P}_i,\mathbb{R}^n)$, while for $\tilde{x}_i \in \mathrm{PL}(\mathcal{P}_i,\mathbb{R}^n)$

$$\begin{aligned}\tilde{\ell}_i(\tilde{x}_i,\zeta_i,\vartheta_i) := \ln\psi(y\,|\,\tilde{x}_i,\tilde{z}_i,\vartheta_i) + \ln\pi(\tilde{x}_i(0),\zeta_i,\vartheta_i) \\ -\frac{1}{2}\sum_{k=0}^{N_i-1}\delta_{ki}\left|\boldsymbol{G}^{-1}\left[\tfrac{\Delta\tilde{x}_{ki}}{\delta_{ki}} - f(t_{ki},\tilde{x}_i(t_{ki}),\tilde{z}_i(t_{ki}),\vartheta_i)\right]\right|^2,\end{aligned} \tag{3.4}$$

where $\Delta\tilde{x}_{ki} := \tilde{x}_i(t_{k+1,i}) - \tilde{x}_i(t_{ki})$ and $\tilde{z}_i \in \mathrm{PL}(\mathcal{P}_i,\mathbb{R}^q)$ is the Euler-discretized clean-state path associated with $\tilde{x}_i$, $\zeta_i$ and $\vartheta_i$, i.e., $\tilde{z}_i(0) = \zeta_i$ and, for all $t\in(t_{ki},t_{k+1,i}]$,

$$\tilde{z}_i(t) = \tilde{z}_i(t_{ki}) + \int_{t_{ki}}^{t} h(t_{ki},\tilde{x}_i(t_{ki}),\tilde{z}_i(t_{ki}),\vartheta_i)\,\mathrm{d}s. \tag{3.5}$$

The convergence of the discretized paths will be taken with respect to the topology of the reproducing kernel Hilbert space $\mathcal{W}_n^2$. This topology can be generated by the following inner product and its induced norm:

$$\langle a, b \rangle_{\mathcal{W}_n^2} := a(0)^\intercal b(0) + \int_0^{t_\mathrm{f}} \dot{a}(t)^\intercal \dot{b}(t) \, \mathrm{d}t,$$

$$\|a\|_{\mathcal{W}_n^2} := \sqrt{\langle a, a \rangle_{\mathcal{W}_n^2}} = \left( |a(0)|^2 + \int_0^{t_\mathrm{f}} |\dot{a}(t)|^2 \, \mathrm{d}t \right)^{1/2}.$$

A useful property of this norm is that it dominates the supremum norm, which is proved in Lemma A.7 of Appendix A. Next, we prove that all $x \in \mathcal{W}_n^2$ can be written as $\mathcal{W}_n^2$-limits of piecewise linear functions over the sequence of partitions.

**Lemma 3.3.** *For all $x \in \mathcal{W}_n^2$, there exists a sequence of $\tilde{x}_i \in \mathrm{PL}(\mathcal{P}_i, \mathbb{R}^n)$ such that $\lim_{i\to\infty} \|\tilde{x}_i - x\|_{\mathcal{W}_n^2} = 0$.*

*Proof.* Let $\{\phi_j\}_{j=1}^\infty$ be a sequence of continuous functions $\phi_j \in \mathcal{C}(\mathcal{T}, \mathbb{R}^n)$ such that

$$\lim_{j\to\infty} \left\| \phi_j - \dot{x} \right\|_{L_n^2} = 0.$$

Such a sequence exists because the space of continuous functions over a compact interval is dense in $L_n^2$. Next, let $\varphi_{ij} \colon \mathcal{T} \to \mathbb{R}^n$ be the right-continuous piecewise constant interpolation of $\phi_j$ over the partition $\mathcal{P}_i$, i.e.,

$$\varphi_{ij}(t) := \sum_{k=0}^{N_i - 1} \phi_j(t_{ki}) I_{[t_{ki}, t_{k+1,i})}(t) + \phi_j(t_\mathrm{f}) I_{\{t_\mathrm{f}\}}(t).$$

Since continuous functions over a compact interval such as $\mathcal{T}$ are absolutelly continuous, $\lim_{i\to\infty} \left\| \varphi_{ij} - \phi_j \right\|_{L_n^2} = 0$.

Now let $\xi_i$ be the $\varphi_{kj}$ closest to $\dot{x}$, for all $k \le i$ and $j \le i$, i.e.,

$$\xi_i := \operatorname*{arg\,min}_{\{\varphi_{kj} \,|\, k,j \in \{1,\dots,i\}\}} \left\| \varphi_{kj} - \dot{x} \right\|_{L_n^2}.$$

Consequently, for all $\varepsilon \in \mathbb{R}_{>0}$ there exist $k, j \in \mathbb{N}$ such that $\left\| \phi_j - \dot{x} \right\|_{L_n^2} \le \frac{\epsilon}{2}$ and $\left\| \varphi_{kj} - \phi_j \right\|_{L_n^2} \le \frac{\epsilon}{2}$, implying that $\|\xi_i - \dot{x}\|_{L_n^2} \le \epsilon$ for all $i \ge \max(k, j)$ and $\lim_{i\to\infty} \|\xi_i - \dot{x}\|_{L_n^2} = 0$. If we then define

$$\tilde{x}_i(t) := x(0) + \int_0^t \xi_i(s) \, \mathrm{d}s,$$

then $\tilde{x}_i$ is piecewise linear and $\lim_{i\to\infty} \|\tilde{x}_i - x\|_{\mathcal{W}_n^2} = 0$. □

Next, we show that the Euler-discretized clean-state path converges uniformly to the clean-state path when the noisy-state path, inital clean state and parameter vector converge. We begin by proving that the sequence of Euler-discretized clean-state paths is bounded.

**Lemma 3.4.** *Let $\{\tilde{x}_i, \zeta_i, \vartheta_i\}_{i=1}^{\infty}$ be a sequence of $\tilde{x}_i \in \mathrm{PL}(\mathcal{P}_i, \mathbb{R}^n)$, $\zeta_i \in \mathbb{R}^q$, and $\vartheta_i \in \mathbb{R}^m$ such that $\lim_{i\to\infty} \|\tilde{x}_i - x\|_{\mathcal{W}_n^2} + |\zeta_i - z_0| + |\vartheta_i - \theta| = 0$ for some $x \in \mathcal{W}_n^2$, $z_0 \in \mathbb{R}^q$, and $\theta \in \operatorname{int}\operatorname{supp}(P_\Theta)$. Then the sequence $\{\|\!|\tilde{z}_i|\!\|\}_{i=1}^{\infty}$ is bounded, where $\tilde{z}_i \in \mathrm{PL}(\mathcal{P}_i)$ is the Euler-discretized clean-state path with $\tilde{z}_i(0) = \zeta_i$ and associated with $\tilde{x}_i$ and $\vartheta_i$, satisfying* (3.5).

*Proof.* First, note that since $\tilde{z}_i$ is piecewise linear,

$$\|\!|\tilde{z}_i|\!\| = \max_{0\le k\le N_i} |\tilde{z}_i(t_{ki})|\,. \tag{3.6}$$

Next, from the definition of the Euler scheme and the triangle inequality,

$$\left|\tilde{z}_i(t_{ji})\right| \le |\zeta_i| + \sum_{k=0}^{j-1} |h(t_{ki}, \tilde{x}_i(t_{ki}), \tilde{z}_i(t_{ki}), \vartheta_i)|\,\delta_{ki}. \tag{3.7}$$

Using the triangle inequality again, we have that the argument of the summation can be expanded as

$$\begin{aligned}
&|h(t_{ki}, \tilde{x}_i(t_{ki}), \tilde{z}_i(t_{ki}), \vartheta_i)| \\
&\qquad \le |h(t_{ki}, \tilde{x}_i(t_{ki}), \tilde{z}_i(t_{ki}), \vartheta_i) - h(t_{ki}, \tilde{x}_i(t_{ki}), \tilde{z}_i(t_{ki}), \theta)| \\
&\qquad\quad + |h(t_{ki}, \tilde{x}_i(t_{ki}), \tilde{z}_i(t_{ki}), \theta) - h(t_{ki}, \tilde{x}_i(t_{ki}), 0, \theta)| \\
&\qquad\quad + |h(t_{ki}, \tilde{x}_i(t_{ki}), 0, \theta)|\,.
\end{aligned} \tag{3.8}$$

The first term of right-hand side of (3.8) is bounded since as $\theta$ lies in the interior of the support of $\Theta$, there exists some $j \in \mathbb{N}$ such that $\vartheta_i \in \operatorname{supp}(P_\Theta)$ for all $i > j$. The function $h$ is assumed to be uniformly continuous in the support of $\Theta$, so

$$\lim_{i\to\infty} |h(t_{ki}, \tilde{x}_i(t_{ki}), \tilde{z}_i(t_{ki}), \vartheta_i) - h(t_{ki}, \tilde{x}_i(t_{ki}), \tilde{z}_i(t_{ki}), \theta)| = 0,$$

and every convergent sequence in $\mathbb{R}$ is bounded. For the second term of right-hand side of (3.8), using the Lipschitz continuity assumption on $h$ we obtain the bound

$$|h(t_{ki}, \tilde{x}_i(t_{ki}), \tilde{z}_i(t_{ki}), \theta) - h(t_{ki}, \tilde{x}_i(t_{ki}), 0, \theta)| \le L_{\mathrm{h}}^{\theta} \, |\tilde{z}_i(t_{ki})|\,.$$

Finally, for the third term of the right-hand side of (3.8), we have that since $\|\!|\tilde{x}_i|\!\|$ converges, it is bounded. Consequently, as continuous funtions over a compact space are bounded, $|h(t_{ki}, \tilde{x}_i(t_{ki}), 0, \theta)|$ is bounded.

Returning to (3.7) and noting that $|\zeta_i|$ is also bounded since it converges, we have that there exists some $a_{17} \in \mathbb{R}_{>0}$ such that

$$\left|\tilde{z}_i(t_{ji})\right| \le a_{17} + \sum_{k=0}^{j-1} L_{\mathrm{h}}^{\theta} \delta_{ki} \, |\tilde{z}_i(t_{ki})|\,.$$

Applying the discrete Grönwall inequality (Clark, 1987), we have that

$$\left|\tilde{z}_i(t_{ji})\right| \leq a_{17} \prod_{k=0}^{j-1} (1 + L_{\mathrm{h}}^{\theta} \delta_{ki}).$$

From Lemma A.3 and (3.6) we then obtain

$$|\!|\!|\tilde{z}_i|\!|\!| \leq a_{17} \exp(L_{\mathrm{h}}^{\theta} t_{\mathrm{f}}) \qquad \text{for all } i \in \mathbb{N}. \qquad \square$$

**Lemma 3.5.** *Let* $\{\tilde{x}_i, \zeta_i, \vartheta_i\}_{i=1}^{\infty}$ *be a sequence of* $\tilde{x}_i \in \mathrm{PL}(\mathcal{P}_i, \mathbb{R}^n)$, $\zeta_i \in \mathbb{R}^q$, *and* $\vartheta_i \in \mathbb{R}^m$ *such that* $\lim_{i\to\infty} \|\tilde{x}_i - x\|_{\mathcal{W}_n^2} + |\zeta_i - z_0| + |\vartheta_i - \theta| = 0$ *for some* $x \in \mathcal{W}_n^2$, $z_0 \in \mathbb{R}^q$, *and* $\theta \in \operatorname{int}\operatorname{supp}(P_\Theta)$. *Then,*

$$\lim_{i\to\infty} |\!|\!|\tilde{z}_i - z|\!|\!| = 0, \tag{3.9}$$

*where* $\tilde{z}_i \in \mathrm{PL}(\mathcal{P}_i)$ *is the Euler-discretized clean-state path with* $\tilde{z}_i(0) = \zeta_i$ *and associated with* $\tilde{x}_i$ *and* $\vartheta_i$, *satisfying* (3.5), *and* $z \in \mathcal{W}_q^2$ *is the unique solution to the initial value problem*

$$\dot{z}(t) = h(t, x(t), z(t), \theta)\,, \qquad z(0) = z_0.$$

*Proof.* From Picard's lemma (Lem. A.10), we have that the distance between $z$ and $\tilde{z}_i$, with respect to the supremum norm, is bounded by

$$|\!|\!|z - \tilde{z}_i|\!|\!| \leq \exp(L_{\mathrm{h}}^{\theta} t_{\mathrm{f}}) \left( |z_0 - \zeta_i| + \int_0^{t_{\mathrm{f}}} \left|\dot{\tilde{z}}_i(t) - h(t, x(t), \tilde{z}_i(t), \theta)\right| \mathrm{d}t \right). \tag{3.10}$$

As $|z_0 - \zeta_i| \to 0$ trivially, to prove that (3.9) holds we have only to show that the integral in the right-hand side of (3.10) vanishes as $i \to \infty$. From (3.5) we have the expression for $\dot{\tilde{z}}_i$, leading to

$$\begin{aligned} &\int_0^{t_{\mathrm{f}}} \left|\dot{\tilde{z}}_i(t) - h(t, x(t), \tilde{z}_i(t), \theta)\right| \mathrm{d}t \\ &\qquad \leq \sum_{k=0}^{N_i - 1} \int_{t_{ki}}^{t_{k+1,i}} |h(t_{ki}, \tilde{x}_i(t_{ki}), \tilde{z}_i(t_{ki}), \vartheta_i) - h(t, x(t), \tilde{z}_i(t), \theta)| \,\mathrm{d}s. \end{aligned} \tag{3.11}$$

Using the triangle inequality, the integrand of (3.11) can be expanded into

$$\begin{aligned} &|h(t_{ki}, \tilde{x}_i(t_{ki}), \tilde{z}_i(t_{ki}), \vartheta_i) - h(t, x(t), \tilde{z}_i(t), \theta)| \\ &\qquad\qquad \leq |h(t_{ki}, \tilde{x}_i(t_{ki}), \tilde{z}_i(t_{ki}), \vartheta_i) - h(t, x(t), \tilde{z}_i(t_{ki}), \theta)| \\ &\qquad\qquad\quad + |h(t, x(t), \tilde{z}_i(t_{ki}), \theta) - h(t, x(t), \tilde{z}_i(t), \theta)|\,. \end{aligned} \tag{3.12}$$

Letting $\rho_x$ denote the modulus of continuity of $x$ and using the triangle inequality again we have that, for all $t \in [t_{ki}, t_{k+1,i}]$,

$$|\tilde{x}_i(t_{ki}) - x(t)| \leq |\tilde{x}_i(t_{ki}) - x(t_{ki})| + |x(t_{ki}) - x(t)| \leq \|\tilde{x}_i - x\| + \rho_x(\bar{\delta}_i).$$

If we then denote by $\rho_h$ the modulus of continuity of $h$, we have that the first term of the right-hand side of (3.12), for all $t \in [t_{ki}, t_{k+1,i}]$, is bounded by

$$\begin{aligned}&|h(t_{ki}, \tilde{x}_i(t_{ki}), \tilde{z}_i(t_{ki}), \vartheta_i) - h(t, x(t), \tilde{z}_i(t_{ki}), \theta)| \\ &\qquad \leq \rho_h\big(\max(\bar{\delta}_i, \|\tilde{x}_i - x\| + \rho_x(\bar{\delta}_i), |\vartheta_i - \theta|)\big),\end{aligned}$$

implying that for (3.9) to hold it suffices to prove that

$$\lim_{i\to\infty} \sum_{k=0}^{N_i-1} \int_{t_{ki}}^{t_{k+1,i}} |h(t, x(t), \tilde{z}_i(t_{ki}), \theta) - h(t, x(t), \tilde{z}_i(t), \theta)| \,\mathrm{d}s = 0. \tag{3.13}$$

Using the Lipschitz property of $h$ we have that the integrand of (3.13) satisfies

$$|h(t, x(t), \tilde{z}_i(t_{ki}), \theta) - h(t, x(t), \tilde{z}_i(t), \theta)| \leq L_{\mathrm{h}}^{\theta} \left|\tilde{z}_i(t_{ki}) - \tilde{z}_i(t)\right|. \tag{3.14}$$

The right-hand side of (3.14) can be further simplified using (3.5):

$$|\tilde{z}_i(t_{ki}) - \tilde{z}_i(t)| \leq \int_{t_{ki}}^{t} |h(t_{ki}, \tilde{x}_i(t_{ki}), \tilde{z}_i(t_{ki}), \vartheta_i)| \,\mathrm{d}s, \tag{3.15}$$

for all $t \in [t_{ki}, t_{k+1,i}]$. Lemmas A.7 and 3.4, together with the fact that the $\tilde{x}_i$, $\zeta_i$ and $\vartheta_i$ converge, imply that the arguments of $h$ in (3.15) are bounded. As $h$ is continuous, this implies that there exists some $a_{18} \in \mathbb{R}_{>0}$ such that, for all $i \in \mathbb{N}$ and $k \in \{0, \dots, N_i\}$,

$$|h(t_{ki}, \tilde{x}_i(t_{ki}), \tilde{z}_i(t_{ki}), \vartheta_i)| \leq a_{18},$$

which implies that, for all $t \in [t_{ki}, t_{k+1,i}]$,

$$|\tilde{z}_i(t_{ki}) - \tilde{z}_i(t)| \leq \int_{t_{ki}}^{t} a_{18} \mathrm{d}s \leq a_{18}\bar{\delta}_i. \tag{3.16}$$

Substituting (3.14) and (3.16) into (3.13), we obtain

$$\sum_{k=0}^{N_i-1} \int_{t_{ki}}^{t_{k+1,i}} |h(t, x(t), \tilde{z}_i(t_{ki}), \theta) - h(t, x(t), \tilde{z}_i(t), \theta)| \,\mathrm{d}s \leq \bar{\delta}_i L_{\mathrm{h}}^{\theta} a_{18} t_{\mathrm{f}} \to 0. \quad \square$$

We are now ready to prove prove that, for convergent sequences of piecewise linear functions, the Euler-discretized log-posterior converges to the energy log-posterior.

**Lemma 3.6.** *Let* $\{\tilde{x}_i, \zeta_i, \vartheta_i\}_{i=1}^{\infty}$ *be a sequence of* $\tilde{x}_i \in \mathrm{PL}(\mathcal{P}_i, \mathbb{R}^n)$, $\zeta_i \in \mathbb{R}^q$, *and* $\vartheta_i \in \mathbb{R}^m$ *such that* $\lim_{i\to\infty} \|\tilde{x}_i - x\|_{\mathcal{W}_n^2} + |\zeta_i - z_0| + |\vartheta_i - \theta| = 0$ *for some* $x \in \mathcal{W}_n^2$, $z_0 \in \mathbb{R}^q$, *and* $\theta \in \mathbb{R}^m$. *Then*

$$\lim_{i\to\infty} \tilde{\ell}(\tilde{x}_i, \zeta_i, \vartheta_i) = \ell_{\mathrm{e}}(x, z_0, \theta). \tag{3.17}$$

*Proof.* First, consider the case where $\theta \notin \operatorname{int}\operatorname{supp}(P_\Theta)$, for which $\ell_{\mathrm{e}}(x, z_0, \theta) = -\infty$. As the summation on the right-hand side of (3.4) is nonpositive, we have that

$$\tilde{\ell}_i(\tilde{x}_i, \zeta_i, \vartheta_i) \leq \ln \psi(y \,|\, \tilde{x}_i, \tilde{z}_i, \vartheta_i) + \ln \pi(\tilde{x}_i(0), \zeta_i, \vartheta_i)\,.$$

Then, due to the continuity of $\pi$ and $\psi$,

$$\limsup_{i\to\infty} \tilde{\ell}(\tilde{x}_i, \zeta_i, \vartheta_i) \leq \limsup_{i\to\infty} \ln \psi(y \,|\, \tilde{x}_i, \tilde{z}_i, \vartheta_i) + \ln \pi(\tilde{x}_i(0), \zeta_i, \vartheta_i) = -\infty.$$

As the limit superior dominates the limit inferior, both coincide and (3.17) holds.

Next, consider the case where $\theta \in \operatorname{int}\operatorname{supp}(P_\Theta)$. From Lemma 3.5 we then have that $\tilde{z}_i \to z$ uniformly. Additionally, from the continuity of $\pi$ and $\psi$, we have that

$$\lim_{i\to\infty} \ln \psi(y \,|\, \tilde{x}_i, \tilde{z}_i, \vartheta_i) + \ln \pi(\tilde{x}_i(0), \zeta_i, \vartheta_i) = \ln \psi(y \,|\, x, z, \theta) + \ln \pi(x(0), z_0, \theta)\,,$$

so for (3.17) to hold it suffices to prove that the summation of the right-hand side of (3.4) converges to the integral of the right-hand side of (2.60). Let $\tilde{h}_i : \mathcal{T} \to \mathbb{R}^n$ be the right-continuous piecewise constant interpolation, with pieces defined by $\mathcal{P}_i$, of $h(t, \tilde{x}_i(t), \tilde{z}_i(t), \vartheta_i)$ using the left endpoint. Then

$$-\frac{1}{2}\sum_{k=0}^{N_i-1} \delta_{ki} \left|\boldsymbol{G}^{-1}\left[\frac{\Delta \tilde{x}_{ki}}{\delta_{ki}} - f(t_{ki}, \tilde{x}_i(t_{ki}), \tilde{z}_i(t_{ki}), \vartheta_i)\right]\right|^2 = -\frac{1}{2}\left\|\dot{\tilde{x}}_i - \tilde{h}_i\right\|_{L_n^2}^2\,.$$

Consequently, as exponentiatin and the norm and continuous operations, for (3.17) to hold it suffices to prove that $\dot{\tilde{x}}_i \to \dot{x}$ in $L_n^2$ and $\tilde{h}_i(t) \to h(t, x(t), z(t), \theta)$ uniformly. The former holds trivially from the definition of the $\mathcal{W}_n^2$ norm. For the latter, note that for all $t \in [t_{ki}, t_{k+1,i}]$,

$$\begin{aligned} |\tilde{x}_i(t_{ki}) - x(t)| &\leq \|\!|\tilde{x}_i - x|\!\| + \rho_x(\bar{\delta}_i), \\ |\tilde{z}_i(t_{ki}) - x(t)| &\leq \|\!|\tilde{z}_i - z|\!\| + \rho_z(\bar{\delta}_i), \end{aligned}$$

where $\rho_x$ and $\rho_z$ are the moduli of continuity of $x$ and $z$, respectively. Together with the uniform continuity property of $h$, this implies that

$$\lim_{i\to\infty} \sup_{t\in\mathcal{T}} \left|\tilde{h}_i(t) - h(t, x(t), z(t), \theta)\right| = 0. \qquad \square$$

We are now ready to prove hypo-convergence of the Euler-discretized log-density.

**Theorem 3.7.** *The Euler-discretized joint state-path and parameter log-posterior* $\tilde{\ell}_i$*, defined in* (3.4)*, hypo-converges to the energy log-posterior* $\ell_\mathrm{e}$ *defined in* (2.60).

*Proof.* From Lemma 3.6 we have that for any convergent sequence $\{\tilde{x}_i, \zeta_i, \vartheta_i\}_{i=1}^\infty$ of $\tilde{x}_i \in \mathrm{PL}(\mathcal{P}_i, \mathbb{R}^n)$, $\zeta_i \in \mathbb{R}^q$, and $\vartheta_i \in \mathbb{R}^m$, then the Euler-discretized log-posterior converges to the energy log-posterior as in (3.17). It suffices that this holds for sequences of $\tilde{x}_i \in \mathrm{PL}(\mathcal{P}_i, \mathbb{R}^n)$ as $\tilde{\ell}_i$ equals negative infinity whenever its first argument lies outside of $\mathrm{PL}(\mathcal{P}_i, \mathbb{R}^n)$. Furthermore, from Lemma 3.3 we have that for all $x \in \mathcal{W}_n^2$, $z_0 \in \mathbb{R}^q$, and $\theta \in \mathbb{R}^m$ there exists such a sequence which converges to $x$, $z_0$, and $\theta$. □

A direct corollary of Theorem 3.7 and Lemma 3.2 is then that any cluster point of any sequence of Euler-discretized MAP estimates is a minimum energy estimate.

## 3.3 Trapezoidally-discretized estimator

The trapezoidal scheme for SDEs (Kloeden and Platen, 1992, p. 500) is converges weakly to the true solution of the SDE with order 2, making it more appropriate to approximate the state-path density then the Euler scheme, as argued by Horsthemke and Bach (1975, p. 191). As in the beginning of Section 3.2, let $\mathcal{P} := \{t_0, \dots, t_N\}$ be a partition of $\mathcal{T}$, with $t_0 = 0$, $t_N = t_\mathrm{f}$, and $t_k < t_{k+1}$. The trapezoidal scheme approximates the SDEs (2.4) at the partition points with the following implicit difference equations:

$$\Delta\hat{X}_k = \frac{1}{2}\left[f(t_k, \hat{X}_{t_k}, \hat{Z}_{t_k}, \Theta) + f(t_{k+1}, \hat{X}_{t_{k+1}}, \hat{Z}_{t_{k+1}}, \Theta)\right]\delta_k + \boldsymbol{G}\Delta W_k, \tag{3.18a}$$

$$\Delta\hat{Z}_k = \frac{1}{2}\left[h(t_k, \hat{X}_{t_k}, \hat{Z}_{t_k}, \Theta) + h(t_{k+1}, \hat{X}_{t_{k+1}}, \hat{Z}_{t_{k+1}}, \Theta)\right]\delta_k, \tag{3.18b}$$

where, as in the previous subsection, the processes $\hat{X}$ and $\hat{Z}$ are the approximations to $X$ and $Z$ with $\hat{X}_0 := X_0$ and $\hat{Z}_0 := Z_0$ and

$$\Delta\hat{X}_k := \hat{X}_{t_{k+1}} - \hat{X}_{t_k}, \qquad \delta_k := t_{k+1} - t_k,$$
$$\Delta\hat{Z}_k := \hat{Z}_{t_{k+1}} - \hat{Z}_{t_k}, \qquad \Delta W_k := W_{t_{k+1}} - W_{t_k}.$$

For the remaining time-points, we consider the trapezoidal approximations to be the piecewise linear interpolation of the values at the partition points, i.e., for all $t \in [t_k, t_{k+1}]$,

$$\hat{X}_t = \hat{X}_{t_k} + \frac{t - t_k}{\delta_k}\Delta\hat{X}_k, \qquad \hat{Z}_t = \hat{Z}_{t_k} + \frac{t - t_k}{\delta_k}\Delta\hat{Z}_k.$$

To ensure that there exists a unique solution to the difference equation (3.18), almost surely, we make the following additional assumptions.

**Assumption 3.8** (trapezoidal scheme)**.**

*a. Restricted to the support of $\Theta$, the functions $f$ and $h$ are Lipschitz continuous with respect to their second and third arguments $x$ and $z$, uniformly over their first and fourth arguments $t$ and $\theta$, i.e., there exist $L_\mathrm{f}, L_\mathrm{h} \in \mathbb{R}_{>0}$ such that, for all $t \in \mathcal{T}$, $x', x'' \in \mathbb{R}^n$, $z', z'' \in \mathbb{R}^q$ and $\theta \in \operatorname{supp}(P_\Theta)$,*

$$\big|f(t, x', z', \theta) - f(t, x'', z'', \theta)\big| \le \big(|x' - x''| + |z' - z''|\big)L_\mathrm{f}, \tag{3.20a}$$
$$\big|h(t, x', z', \theta) - h(t, x'', z'', \theta)\big| \le \big(|x' - x''| + |z' - z''|\big)L_\mathrm{h}. \tag{3.20b}$$

*b. The partition is sufficiently fine such that $(L_\mathrm{f} + L_\mathrm{h})\bar{\delta} < 2$, where $\bar{\delta} := \max_{0 \le k < N} \delta_k$ is the mesh of the partition $\mathcal{P}$.*

The next lemma follows as a corollary of Assumption 3.8.

**Lemma 3.9** (existence and uniqueness of solutions to the trapezoidal scheme)**.** *If Assumption 3.8 is satisfied, then a unique solution to* (3.18) *exists almost surely.*

*Proof.* First, note that $\hat{X}_0 = X_0$ and $\hat{Z}_0 = Z_0$ by definition. Next, assume that there exist unique $\hat{X}_{t_0}, \dots, \hat{X}_{t_k}$ and $\hat{Z}_{t_0}, \dots, \hat{Z}_{t_k}$ that satisfy (3.18). Then it suffices to prove that there exist unique $\hat{X}_{t_{k+1}}$ and $\hat{Z}_{t_{k+1}}$ that satisfy (3.18); the proposition then follows by induction.

If we define $q\colon \mathbb{R}^n \times \mathbb{R}^q \to \mathbb{R}^n \times \mathbb{R}^q$ by

$$q(x, z) := \begin{bmatrix} X_{t_k} + \frac{1}{2} f(t_k, \hat{X}_{t_k}, \hat{Z}_{t_k}, \Theta)\delta_k + \frac{1}{2} f(t_{k+1}, x, z, \Theta)\delta_k + \boldsymbol{G}\Delta W_k \\ Z_{t_k} + \frac{1}{2} h(t_k, \hat{X}_{t_k}, \hat{Z}_{t_k}, \Theta)\delta_k + \frac{1}{2} h(t_{k+1}, x, z, \Theta)\delta_k \end{bmatrix},$$

then (3.18) is satisfied if and only if

$$[\hat{X}_{t_{k+1}}, \hat{Z}_{t_{k+1}}] = q(\hat{X}_{t_{k+1}}, \hat{Z}_{t_{k+1}}). \tag{3.21}$$

Under the norm $\|[x, z]\| := |x| + |z|$ in $\mathbb{R}^n \times \mathbb{R}^q$, the function $q$ is Lipschitz continuous with Lipschitz constant $\delta_k(L_\mathrm{f} + L_\mathrm{h})/2$. Assumption 3.8b then implies that $q$ is a contraction mapping which, by the Banach fixed-point theorem, admits a unique solution that satisfies (3.21). □

Having proved that the trapezoidal scheme is well-defined, we can then proceed to obtain the joint posterior state path and parameter density for the discretized process. The joint density of $\hat{X}_{t_0}, \dots, \hat{X}_{t_N}$, $\hat{Z}_0$ and $\Theta$, which will be denoted the trapezoidally-discretized joint state-path and parameter density, can be found from the joint density of $\hat{X}_0$, $\hat{Z}_0$, $\Theta$, and $\Delta W_0, \dots, \Delta W_{N-1}$ by a change of variables (Lemma A.20), since from the definition of the trapezoidal scheme in (3.18) we have that both groups of variables are related by a diffeomorphism. Consequently, the trapezoidally-discretized joint posterior state-path and parameter density is given by

$$\begin{aligned} p(\hat{x}, \hat{z}_0, \vartheta \,|\, y) = {} & \frac{\psi(y \,|\, \hat{x}, \hat{z}, \vartheta)}{\mathrm{E}[\psi(y \,|\, X, Z, \Theta)]} \pi(\hat{x}(0), \hat{z}_0, \vartheta) \\ & \times \prod_{k=0}^{N-1} \det\bigl(\boldsymbol{I} - \tfrac{1}{2} \nabla_{\mathrm{x}} f(t_{k+1}, \hat{x}(t_{k+1}), \hat{z}(t_{k+1}), \vartheta)\, \delta_k\bigr) \\ & \times \prod_{k=0}^{N-1} \frac{\exp\Bigl(-\delta_k \frac{1}{2} \left|\boldsymbol{G}^{-1} \left[\frac{\Delta \hat{x}_k}{\delta_k} - \frac{1}{2}\hat{f}_k - \frac{1}{2}\hat{f}_{k+1}\right]\right|^2\Bigr)}{|\det G| \sqrt{\delta_k (2\pi)^n}}, \end{aligned} \quad (3.22)$$

where $\Delta \hat{x}_k := \hat{x}(t_{k+1}) - \hat{x}(t_k)$,

$$\hat{f}_k := f(t_k, \hat{x}(t_k), \hat{z}(t_k), \vartheta)\,, \qquad \hat{h}_k := h(t_k, \hat{x}(t_k), \hat{z}(t_k), \vartheta)\,, \quad (3.23)$$

and $\tilde{z}_i \in \mathrm{PL}(\mathcal{P}, \mathbb{R}^q)$ is the Euler-discretized clean-state path associated with $\hat{x}$, $\hat{z}_0$ and $\vartheta$, i.e., $\hat{z}(0) = \hat{z}_0$ and, for all $t \in (t_k, t_{k+1}]$, $\hat{z}$ is the unique solution to the implicit equation

$$\hat{z}(t) = \hat{z}(t_k) + \frac{1}{2} \int_{t_k}^{t} \left(\hat{h}_k + \hat{h}_{k+1}\right) \mathrm{d}s.$$

Because of its implicit nonlinear formulation, the trapezoidally-discretized state-path may depend nonlinearly on the associated noise path. Consequently, its density may not be proportional to the associated noise path density; it also depends on the change of the volume element, quantified by the jacobian determinant of the diffeomorphism which relates both groups of variables. This appears in the joint posterior state path and parameter density as the multiplication on the second line of the right-hand side of (3.22).

Due to better numerical tractability, it is more convenient to minimize the logarithm of the posterior density instead of the density itself. Taking the logarithm of (3.22) and removing the constant terms, which do not change the location of maxima, we obtain $\hat{\ell}\colon \mathrm{PL}(\mathcal{P}, \mathbb{R}^n) \times \mathbb{R}^q \times \mathbb{R}^m \to \overline{\mathbb{R}}$, the trapezoidally-discretized joint state-path and parameters log-posterior:

$$\begin{aligned} \hat{\ell}(\hat{x}, \hat{z}_0, \vartheta) := {} & \ln \psi(y \,|\, \hat{x}, \hat{z}, \vartheta) + \ln \pi(\hat{x}(0), \hat{z}_0, \vartheta) \\ & + \sum_{k=0}^{N-1} \ln \det\Bigl(\boldsymbol{I} - \frac{1}{2} \nabla_{\mathrm{x}} f(t_{k+1}, \hat{x}(t_{k+1}), \hat{z}(t_{k+1}), \vartheta)\, \delta_k\Bigr) \\ & - \frac{1}{2} \sum_{k=0}^{N-1} \delta_k \left|\boldsymbol{G}^{-1} \left[\frac{\Delta \hat{x}_k}{\delta_k} - \frac{1}{2}\hat{f}_k - \frac{1}{2}\hat{f}_{k+1}\right]\right|^2, \end{aligned}$$

where, as in (3.22), $\hat{f}_k$ is given by (3.23) and $\hat{z}$ is the trapezoidally-discretized clean-state path associated with $\hat{x}$, $\hat{z}_0$, and $\vartheta$. We now show that, as the partition mesh vanishes, the trapezoidally-discretized log-posterior hypo-converges to the log-posterior.

### 3.3.1 Hypo-convergence of the trapezoidally-discretized log-posterior

Let $\{\mathcal{P}_i\}_{i=1}^{\infty}$ be a sequence of nested partitions $\mathcal{P}_i := \{t_{k,i}\}_{k=0}^{N_i}$ of $\mathcal{T}$ with a vanishing mesh, i.e., $\mathcal{P}_i \subset \mathcal{P}_{i+1}$,

$$\delta_{ki} := t_{k+1,i} - t_{ki}, \qquad \bar{\delta}_i := \max_{0 \le k < N_i} \delta_{ki}, \qquad \lim_{i \to \infty} \bar{\delta}_i = 0.$$

The variable $\bar{\delta}_i$ represents the mesh of $\mathcal{P}_i$, and we assume that it is sufficiently fine such that Assumption 3.8b is satisfied for all $i \in \mathbb{N}$. The trapezoidal scheme is then well-defined for all partitions.

For each $\mathcal{P}_i$, we consider its corresponding trapezoidally-discretized joint state-path and parameter log-posterior $\hat{\ell}_i \colon \mathcal{W}_n^2(\mathcal{T}) \times \mathbb{R}^q \times \mathbb{R}^m \to \overline{\mathbb{R}}$. As explained in Section 3.1, we extend its domain to non piecewise linear functions by setting its value to $-\infty$, i.e., $\hat{\ell}_i(x, \zeta_i, \vartheta_i) := -\infty$ for $x \notin \mathrm{PL}(\mathcal{P}_i, \mathbb{R}^n)$, while for $\hat{x}_i \in \mathrm{PL}(\mathcal{P}_i, \mathbb{R}^n)$

$$\begin{aligned}\hat{\ell}_i(\hat{x}_i, \zeta_i, \vartheta_i) &:= \ln \psi(y \,|\, \hat{x}_i, \hat{z}_i, \vartheta_i) + \ln \pi(\hat{x}_i(0), \zeta_i, \vartheta_i) \\ &\quad + \sum_{k=0}^{N_i-1} \ln \det\big(\boldsymbol{I} - \tfrac{1}{2}\nabla_{\mathrm{x}} f\big(t_{k+1,i}, \hat{x}_i(t_{k+1,i}), \hat{z}_i(t_{k+1,i}), \vartheta_i\big)\, \delta_{ki}\big) \\ &\qquad - \frac{1}{2}\sum_{k=0}^{N_i-1} \delta_{ki} \left| \boldsymbol{G}^{-1} \left[ \tfrac{\Delta \hat{x}_{ki}}{\delta_{ki}} - \tfrac{1}{2}\hat{f}_{ki} - \tfrac{1}{2}\hat{f}_{k+1,i} \right] \right|^2, \end{aligned} \tag{3.24}$$

where $\Delta \hat{x}_{ki} := \hat{x}_i(t_{k+1}) - \tilde{x}(t_k)$,

$$\hat{f}_{ki} := f(t_{ki}, \hat{x}_i(t_{ki}), \hat{z}_i(t_{ki}), \vartheta_i)\,, \qquad \hat{h}_{ki} := h(t_{ki}, \hat{x}_i(t_{ki}), \hat{z}_i(t_{ki}), \vartheta_i)\,,$$

and $\tilde{z}_i \in \mathrm{PL}(\mathcal{P}_i, \mathbb{R}^q)$ is the trapezoidally-discretized clean-state path associated with $\tilde{x}_i$, $\zeta_i$ and $\vartheta_i$, i.e., $\tilde{z}_i(0) = \zeta_i$ and, for all $t \in (t_{ki}, t_{k+1,i}]$, $\tilde{z}_i$ is the unique solution to the implicit equation

$$\hat{z}_i(t) = \hat{z}_i(t_{ki}) + \frac{1}{2}\int_{t_{ki}}^{t} \left(\hat{h}_{ki} + \hat{h}_{k+1,i}\right) \mathrm{d}s. \tag{3.25}$$

As in the previous subsection, the convergence of the discretized paths will be taken with respect to the topology of the reproducing kernel Hilbert space $\mathcal{W}_n^2$. The proof of hypo-convergence of the trapezoidally-discretized log-posterior will follow roughly the same steps as that of the Euler-discretized log-posterior. We begin by proving that, for any convergent sequence of its arguments, the associated trapezoidally-discretized clean-state path sequence is uniformly bounded.

**Lemma 3.10.** *Let $\{\hat{x}_i, \zeta_i, \vartheta_i\}_{i=1}^{\infty}$ be a sequence of $\hat{x}_i \in \mathrm{PL}(\mathcal{P}_i, \mathbb{R}^n)$, $\zeta_i \in \mathbb{R}^q$, and $\vartheta_i \in \mathbb{R}^m$ such that*

$$\lim_{i\to\infty} \|\hat{x}_i - x\|_{\mathcal{W}_n^2} + |\zeta_i - z_0| + |\vartheta_i - \theta| = 0$$

*for some $x \in \mathcal{W}_n^2$, $z_0 \in \mathbb{R}^q$, and $\theta \in \operatorname{int}\operatorname{supp}(P_\Theta)$. Then the sequence $\{\|\hat{z}_i\|\}_{i=1}^{\infty}$ is bounded, where $\hat{z}_i \in \mathrm{PL}(\mathcal{P}_i)$ is the trapezoidally-discretized clean-state path with $\hat{z}_i(0) = \zeta_i$ and associated with $\hat{x}_i$ and $\vartheta_i$, satisfying* (3.25)*.*

*Proof.* For any $i \in \mathbb{N}$ and $k \in \{0, \dots, N_i - 1\}$, define $q\colon \mathbb{R}^q \to \mathbb{R}^q$ by

$$q(z) := \hat{z}_i(t_{ki}) + \tfrac{1}{2}\hat{h}_{ki}\delta_{ki} + \tfrac{1}{2}h\big(t_{k+1,i}, \hat{x}_i(t_{k+1,i}), z, \vartheta_i\big)\,\delta_{ki}. \tag{3.26}$$

We then have that $\hat{z}_i(t_{k+1,i})$ is the unique solution to

$$\hat{z}_i(t_{k+1,i}) = q(\hat{z}_i(t_{k+1,i}))$$

and $q$ is Lipschitz continuous with Lipschitz constant $L_\mathrm{q} = \frac{1}{2}L_\mathrm{h}\delta_{ki}$. Applying the Banach fixed-point theorem with the iteration starting at $\hat{z}_i(t_{ki})$, we have that

$$\big|\hat{z}_i(t_{k+1,i}) - \hat{z}_i(t_{ki})\big| \le \frac{1}{1 - L_\mathrm{q}}\,|\hat{z}_i(t_{ki}) - q(\hat{z}_i(t_{ki}))|\,,$$

and, from (3.26), and the triangle inequality,

$$|\hat{z}_i(t_{ki}) - q(\hat{z}_i(t_{ki}))| \le \tfrac{1}{2}\delta_{ki}\,\Big[\big|\hat{h}_{ki}\big| + \big|h\big(t_{k+1,i}, \hat{x}_i(t_{k+1,i}), \hat{z}_i(t_{ki}), \vartheta_i\big)\big|\Big]\,.$$

Assumption 3.8b implies that $L_\mathrm{q} < 1$, so defining $a_{19} := \frac{1}{1-L_\mathrm{q}}$ we have that

$$\big|\hat{z}_i(t_{ji})\big| \le |\zeta_i| + \sum_{k=0}^{j-1}\Big[\big|\hat{h}_{ki}\big| + \big|h\big(t_{k+1,i}, \hat{x}_i(t_{k+1,i}), \hat{z}_i(t_{ki}), \vartheta_i\big)\big|\Big]\,\tfrac{1}{2}a_{19}\delta_{ki}. \tag{3.27}$$

Furthermore, using the triangle inequality we have that

$$\big|\hat{h}_{ki}\big| \le \big|\hat{h}_{ki} - h\big(t_{k+1,i}, \hat{x}_i(t_{k+1,i}), 0, \vartheta_i\big)\big| + \big|h\big(t_{k+1,i}, \hat{x}_i(t_{k+1,i}), 0, \vartheta_i\big)\big|\,. \tag{3.28}$$

and

$$\begin{aligned}
&\big|h\big(t_{k+1,i}, \hat{x}_i(t_{k+1,i}), \hat{z}_i(t_{ki}), \vartheta_i\big)\big| \\
&\qquad\le \big|h\big(t_{k+1,i}, \hat{x}_i(t_{k+1,i}), \hat{z}_i(t_{ki}), \vartheta_i\big) - h\big(t_{k+1,i}, \hat{x}_i(t_{k+1,i}), 0, \vartheta_i\big)\big| \\
&\qquad\quad + \big|h\big(t_{k+1,i}, \hat{x}_i(t_{k+1,i}), 0, \vartheta_i\big)\big|\,.
\end{aligned} \tag{3.29}$$

From the Lipschitz continuity of $h$ (3.20b), we have that the first terms of the right-hand side of both (3.28) and (3.29) are bounded by $L_\mathrm{h}\,|\hat{z}_i(t_{ki})|$. The second terms of both equations, in turn, are bounded since $h$ is continuous and its arguments are bounded. Returning to (3.27) and noting that $|\zeta_i|$ is also bounded since it converges, we have that there exists some $a_{20} \in \mathbb{R}_{>0}$, indepencent of $i$, such that

$$\left|\hat{z}_i(t_{ji})\right| \leq a_{20} + \sum_{k=0}^{j-1} |\hat{z}_i(t_{ki})| \, L_{\mathrm{h}} a_{19} \delta_{ki}.$$

Applying the discrete Grönwall inequality (Clark, 1987), we have that

$$\left|\hat{z}_i(t_{ji})\right| \leq a_{20} \prod_{k=0}^{j-1} (1 + a_{19} L_{\mathrm{h}} \delta_{ki}).$$

From Lemma A.3 and (3.6) we then obtain

$$\|\hat{z}_i\| \leq a_{20} \exp(a_{19} L_{\mathrm{h}} t_{\mathrm{f}}). \qquad \text{for all } i \in \mathbb{N}. \qquad \square$$

**Corollary 3.11.** *Let $\{\hat{x}_i, \zeta_i, \vartheta_i\}_{i=1}^{\infty}$ be a sequence of $\hat{x}_i \in \mathrm{PL}(\mathcal{P}_i, \mathbb{R}^n)$, $\zeta_i \in \mathbb{R}^q$, and $\vartheta_i \in \mathbb{R}^m$ such that*

$$\lim_{i \to \infty} \|\hat{x}_i - x\|_{\mathcal{W}_n^2} + |\zeta_i - z_0| + |\vartheta_i - \theta| = 0$$

*for some $x \in \mathcal{W}_n^2$, $z_0 \in \mathbb{R}^q$, and $\theta \in \operatorname{int}\operatorname{supp}(P_\Theta)$. Furthermore, let $\{\hat{z}_i\}_{i=1}^{\infty}$ be the sequence of trapezoidally-discretized clean-state paths $\hat{z}_i \in \mathrm{PL}(\mathcal{P}_i)$ with $\hat{z}_i(0) = \zeta_i$ and associated with $\hat{x}_i$ and $\vartheta_i$, satisfying* (3.25)*. Then the $\hat{z}_i$ are Lipschitz equicontinuous, i.e., there exists some $L_{\hat{\mathrm{z}}} \in \mathbb{R}_{>0}$ such that, for all $t, \tau \in \mathcal{T}$ and $i \in \mathbb{N}$,*

$$|\hat{z}_i(t) - \hat{z}_i(\tau)| \leq L_{\hat{\mathrm{z}}} \, |t - \tau| \,. \tag{3.30}$$

*Proof.* From (3.25) we have that, for all $t, \tau \in [t_{ki}, t_{k+1,i}]$ such that $t \leq \tau$

$$|\hat{z}_i(t) - \hat{z}_i(\tau)| \leq \frac{1}{2} \int_t^{\tau} \left|\hat{h}_{ki} + \hat{h}_{k+1,i}\right| \mathrm{d}s.$$

Furthermore, note that $\mathcal{T}$ is bounded by definition; $\hat{x}_i$ is uniformly bounded in $i$ and $t$ since it is convergent in $\mathcal{W}_n^2$ and the supremum norm is dominated by the $\mathcal{W}_n^2$ norm (Lemma A.7); $\hat{z}_i$ is uniformly bounded in $i$ and $t$ from Lemma 3.10; and $\vartheta_i$ is uniformly bounded since it is convergent. Consequently, $\hat{h}_{ki}$ is uniformly bounded in $k$ and $i$ and (3.30) is satisfied with

$$L_{\hat{\mathrm{z}}} := \sup_{i \in \mathbb{N}} \sup_{0 \leq k \leq N_i} \left|\hat{h}_{ki}\right|. \qquad \square$$

We now prove that the trapezoidally-discretized clean-state path sequence converges to the clean-state path associated with the limits of the noisy-state path, initial clean state and parameter vector.

**Lemma 3.12.** *Let* $\{\hat{x}_i, \zeta_i, \vartheta_i\}_{i=1}^{\infty}$ *be a sequence of* $\hat{x}_i \in \mathrm{PL}(\mathcal{P}_i, \mathbb{R}^n)$, $\zeta_i \in \mathbb{R}^q$, *and* $\vartheta_i \in \mathbb{R}^m$ *such that*

$$\lim_{i\to\infty} \|\hat{x}_i - x\|_{\mathcal{W}_n^2} + |\zeta_i - z_0| + |\vartheta_i - \theta| = 0$$

*for some* $x \in \mathcal{W}_n^2$, $z_0 \in \mathbb{R}^q$, *and* $\theta \in \operatorname{int}\operatorname{supp}(P_\Theta)$. *Then,*

$$\lim_{i\to\infty} \|\hat{z}_i - z\| = 0, \tag{3.31}$$

*where* $\hat{z}_i \in \mathrm{PL}(\mathcal{P}_i)$ *is the trapezoidally-discretized clean-state path associated with* $\hat{x}_i$ *and* $\vartheta_i$, *satisfying* (3.25) *and* $\hat{z}_i(0) = \zeta_i$*; and* $z \in \mathcal{W}_q^2$ *is the unique solution to the initial value problem*

$$\dot{z}(t) = h(t, x(t), z(t), \theta)\,, \qquad z(0) = z_0.$$

*Proof.* From Picard's lemma (Lem. A.10), we have that the distance between $z$ and $\hat{z}_i$, with respect to the supremum norm, is bounded by

$$\|z - \hat{z}_i\| \le \exp(L_\mathrm{h}^\theta t_\mathrm{f}) \left( |z_0 - \zeta_i| + \int_0^{t_\mathrm{f}} \left|\dot{\hat{z}}_i(t) - h(t, x(t), \hat{z}_i(t), \theta)\right| \mathrm{d}t \right). \tag{3.32}$$

As $|z_0 - \zeta_i| \to 0$ trivially, to prove that (3.31) holds we have only to show that the integral in the right-hand side of (3.32) vanishes as $i \to \infty$. From (3.25) we have the expression for $\dot{\hat{z}}_i$, leading to

$$\int_0^{t_\mathrm{f}} \left|\dot{\hat{z}}_i(t) - h(t, x(t), \hat{z}_i(t), \theta)\right| \mathrm{d}t \le \sum_{k=0}^{N_i - 1} \int_{t_{ki}}^{t_{k+1,i}} \left|\tfrac{1}{2}\hat{h}_{ki} + \tfrac{1}{2}\hat{h}_{k+1,i} - h(t, x(t), \hat{z}_i(t), \theta)\right| \mathrm{d}t.$$

Letting $\rho_x$ denote the modulus of continuity of $x$ and using the triangle inequality we have that, for all $t \in [t_{ki}, t_{k+1,i}]$,

$$|\hat{x}_i(t_{ki}) - x(t)| \le |\hat{x}_i(t_{ki}) - x(t_{ki})| + |x(t_{ki}) - x(t)| \le \|\hat{x}_i - x\| + \rho_x(\bar{\delta}_i). \tag{3.33}$$

Together with Corollary 3.11 and the uniform continuity of $h$, this implies that for all $t \in [t_{ki}, t_{k+1,i}]$,

$$\left|\hat{h}_{ki} - h(t, x(t), \hat{z}_i(t), \theta)\right| \le \rho_h\big(\max(\bar{\delta}_i, \|\hat{x}_i - x\| + \rho_x(\bar{\delta}_i), L_{\hat{z}}\bar{\delta}_i, |\vartheta_i - \theta|)\big)\,,$$

with the same bound holding for $\left|\hat{h}_{k+1,i} - h(t, x(t), \hat{z}_i(t), \theta)\right|$, implying that

$$\int_0^{t_\mathrm{f}} \left|\dot{\hat{z}}_i(t) - h(t, x(t), \hat{z}_i(t), \theta)\right| \mathrm{d}t \le t_\mathrm{f}\rho_h\big(\max(\bar{\delta}_i, \|\hat{x}_i - x\| + \rho_x(\bar{\delta}_i), L_{\hat{z}}\bar{\delta}_i, |\vartheta_i - \theta|)\big) \to 0. \quad \square$$

We are now ready to prove prove that, for convergent sequences of piecewise linear functions, the trapezoidally-discretized log-posterior converges to the continuous log-posterior.

**Lemma 3.13.** *Let* $\{\hat{x}_i, \zeta_i, \vartheta_i\}_{i=1}^{\infty}$ *be a sequence of* $\hat{x}_i \in \mathrm{PL}(\mathcal{P}_i, \mathbb{R}^n)$, $\zeta_i \in \mathbb{R}^q$, *and* $\vartheta_i \in \mathbb{R}^m$ *such that*

$$\lim_{i\to\infty} \|\hat{x}_i - x\|_{\mathcal{W}_n^2} + |\zeta_i - z_0| + |\vartheta_i - \theta| = 0 \tag{3.34}$$

*for some* $x \in \mathcal{W}_n^2$, $z_0 \in \mathbb{R}^q$, *and* $\theta \in \mathbb{R}^m$. *Then*

$$\lim_{i\to\infty} \hat{\ell}(\hat{x}_i, \zeta_i, \vartheta_i) = \ell(x, z_0, \theta). \tag{3.35}$$

*Proof.* First, consider the case where $\theta \notin \operatorname{int}\operatorname{supp}(P_\Theta)$, for which $\ell_{\mathrm{e}}(x, z_0, \theta) = -\infty$. As the summation on the right-hand side of (3.24) is nonpositive, we have that

$$\hat{\ell}_i(\hat{x}_i, \zeta_i, \vartheta_i) \le \ln \psi(y \,|\, \hat{x}_i, \hat{z}_i, \vartheta_i) + \ln \pi(\hat{x}_i(0), \zeta_i, \vartheta_i)\,.$$

Then, due to the continuity of $\pi$ and $\psi$,

$$\limsup_{i\to\infty} \hat{\ell}(\hat{x}_i, \zeta_i, \vartheta_i) \le \limsup_{i\to\infty} \ln \psi(y \,|\, \hat{x}_i, \hat{z}_i, \vartheta_i) + \ln \pi(\hat{x}_i(0), \zeta_i, \vartheta_i) = -\infty.$$

As the limit superior dominates the limit inferior, both coincide and (3.35) holds.

Next, consider the case where $\theta \in \operatorname{int}\operatorname{supp}(P_\Theta)$. From Lemma 3.12 we then have that $\hat{z}_i \to z$ uniformly. Additionally, from the continuity of $\pi$ and $\psi$, we have that

$$\lim_{i\to\infty} \ln \psi(y \,|\, \hat{x}_i, \hat{z}_i, \vartheta_i) + \ln \pi(\hat{x}_i(0), \zeta_i, \vartheta_i) = \ln \psi(y \,|\, x, z, \theta) + \ln \pi(x(0), z_0, \theta)\,,$$

so for (3.35) to hold it suffices to prove that the summations on the second line and third lines of the right-hand side of (3.24) converge to the drift-divergence integral and the energy term of (2.51), respectively.

Let the noisy-state drift Jacobian at the partition points be denoted by

$$\hat{\boldsymbol{J}}_{ki} := \nabla_{\mathrm{x}} f(t_{ki}, \hat{x}_i(t_{ki}), \hat{z}_i(t_{ki}), \vartheta_i)\,. \tag{3.36}$$

As $f$ is assumed to be differentiable with respect to its second argument $x$, Assumption 3.8a implies that $\left|\hat{\boldsymbol{J}}_{ki}\right| \le L_{\mathrm{f}}$. As the spectral radius is dominated by consistent matrix norms, Assumption 3.8b implies that the spectral radius of $\frac{1}{2}\hat{\boldsymbol{J}}_{k+1,i}\delta_{ki}$ is smaller than unity for all $k$ and $i$. Lemmas A.4 and A.5 then imply that

$$\ln\det\!\left(\boldsymbol{I} - \tfrac{1}{2}\hat{\boldsymbol{J}}_{ki}\delta_{ki}\right) = \sum_{j=1}^{\infty} (-1)^j \frac{\delta_{ki}^j}{2j} \operatorname{tr}\!\left(\hat{\boldsymbol{J}}_{ki}^j\right).$$

Truncating the series, we have that by Assumption 3.8b there exists some $a_{22} \in \mathbb{R}_{>0}$ such that

$$\left|-\tfrac{1}{2}\operatorname{tr}(\hat{\boldsymbol{J}}_{ki})\delta_{ki} - \ln\det\!\left(\boldsymbol{I} - \tfrac{1}{2}\hat{\boldsymbol{J}}_{ki}\delta_{ki}\right)\right| \leq a_{22}\bar{\delta}_i^2.$$

In addition, from (3.34) and Lemmas 3.10 and A.7, we have that the arguments of the drift Jacobian in the right-hand side of (3.36) are bounded. As the Jacobian is assumed to be continuous, this implies that it is uniformly continuous on the subset of its domain being analysed. Denoting by $\rho_{\mathrm{J}}$ its modulus of continuity we have that, for all $t \in [t_{ki}, t_{k+1,i}]$,

$$\begin{aligned}&\left|\operatorname{tr}(\hat{\boldsymbol{J}}_{ki}) - \operatorname{div}_{\mathrm{x}} f(t, x(t), z(t), \theta)\right| \\ &\qquad\qquad \leq \rho_{\mathrm{J}}\big(\max(\bar{\delta}_i, \|\hat{x}_i - x\| + \rho_x(\bar{\delta}_i), L_{\hat{\mathrm{z}}}\bar{\delta}_i, |\vartheta_i - \theta|)\big), \quad (3.37)\end{aligned}$$

where (3.33) and Corollary 3.11 were used. Consequently,

$$\lim_{i\to\infty} \sum_{k=0}^{N_i-1} \ln\det\!\left(\boldsymbol{I} - \tfrac{1}{2}\hat{\boldsymbol{J}}_{k+1,i}\delta_{ki}\right) = -\frac{1}{2}\int_0^{t_{\mathrm{f}}} \operatorname{div}_{\mathrm{x}} f(t, x(t), z(t), \theta)\,.$$

Similarly to (3.37), for the energy term we have that for all $t \in [t_{ki}, t_{k+1,i}]$,

$$\left|\hat{f}_{ki} - f(t, x(t), \hat{z}_i(t), \theta)\right| \leq \rho_f\big(\max(\bar{\delta}_i, \|\hat{x}_i - x\| + \rho_x(\bar{\delta}_i), L_{\hat{\mathrm{z}}}\bar{\delta}_i, |\vartheta_i - \theta|)\big)\,,$$

with the same bound holding for $\left|\hat{f}_{k+1,i} - f(t, x(t), \hat{z}_i(t), \theta)\right|$, implying that

$$\begin{aligned}&\lim_{i\to\infty} \sum_{k=0}^{N_i-1} \left|\boldsymbol{G}^{-1}\left[\tfrac{\Delta\hat{x}_{ki}}{\delta_{ki}} - \tfrac{1}{2}\hat{f}_{ki} - \tfrac{1}{2}\hat{f}_{k+1,i}\right]\right|^2 \delta_{ki} \\ &\qquad\qquad = \int_0^{t_{\mathrm{f}}} \left|\boldsymbol{G}^{-1}[\dot{x}(t) - f(t, x(t), z(t), \theta)]\right|^2 \mathrm{d}t. \quad \square\end{aligned}$$

We are now ready to prove hypo-convergence of the trapezoidally-discretized log-density.

**Theorem 3.14.** *The trapezoidally-discretized joint state-path and parameter log-posterior $\hat{\ell}_i$, defined in* (3.24)*, hypo-converges to the continuous-time log-posterior $\ell$ defined in* (2.49)

*Proof.* For any convergent sequence $\{\hat{x}_i, \zeta_i, \vartheta_i\}_{i=1}^{\infty}$ of $\hat{x}_i \in \mathrm{PL}(\mathcal{P}_i, \mathbb{R}^n)$, $\zeta_i \in \mathbb{R}^q$, and $\vartheta_i \in \mathbb{R}^m$, by Lemma 3.13 we have that the trapezoidally-discretized log-posterior converges to the continuous log-posterior as in (3.35). Note that it suffices that this holds for sequences of $\hat{x}_i \in \mathrm{PL}(\mathcal{P}_i, \mathbb{R}^n)$ as $\hat{\ell}_i$ equals negative infinity whenever its first argument lies outside of $\mathrm{PL}(\mathcal{P}_i, \mathbb{R}^n)$. Furthermore, from Lemma 3.3 we have that for all $x \in \mathcal{W}_n^2$, $z_0 \in \mathbb{R}^q$, and $\theta \in \mathbb{R}^m$ there exists such a sequence which converges to $x$, $z_0$, and $\theta$. $\square$

A direct corollary of Theorem 3.14 and Lemma 3.2 is then that any cluster point of any sequence of trapezoidally-discretized MAP estimates is a MAP estimate of the continuous-time system.

# Chapter 4

# Example applications

> Figures often beguile me, particularly when I have the arranging of them myself; in which case the remark attributed to Disraeli would often apply with justice and force: "There are three kinds of lies: lies, damned lies, and statistics."
>
> Mark Twain, *Chapters from My Autobiography*

In this chapter, we demonstrate example applications of the proposed methods with both simulated and experimental data. Direct transcription methods were used to implement both the joint maximum *a posteriori* state path and parameter estimator (JMAPSPPE) and the minimum energy estimator (MEE). The variational optimization problems were translated to nonlinear programming problems using a direct transcription technique equivalent to the Hermite–Simpson method (Betts, 2010, Sec. 4.5), and the resulting nonlinear program was solved using the IPOPT solver of Wächter and Biegler (2006), part of the COIN-OR project. The large-scale sparse linear systems underlying the optimization were solved with the MA57 solver of the HSL Mathematical Software Library.[1]

As argued at the end of Section 2.2.1, the regularization of Remark A.11 can be used to ensure that the systems presented in this section satisfy Assumption 2.8. The states and parameters are meaningless too far away from the origin and the regularized systems can be used without loss of generality.

## 4.1 Simulated examples

In this section we present simulated example applications on benchmark nonlinear models. All stochastic differential equations (SDEs) were simulated using the strong explicit order 1.5 scheme (Kloeden and Platen, 1992, Sec. 11.2), with a time step of 0.005 and the initial states sampled from the initial-state prior. All realizations of each simulation were performed with the same nominal parameter values. Both the JMAPSPPE and MEE were applied to the simulated data and their estimates compared to those of the prediction error method (PEM, cf. Kristensen et al., 2004a) using the unscented Kalman filter (UKF).

[1]HSL(2013). A collection of Fortran codes for large scale scientific computation. `http://www.hsl.rl.ac.uk`

To implement the PEM, the unscented SDE prediction step of Arasaratnam et al. (2010) was used. The unscented Kalman smoother with the backward correction step of Särkkä (2008) was used to obtain the optimal state-path associated with the PEM-estimated MAP parameters. To be more favorable with the PEM and avoid local minima in its optimization, the nominal parameter values were used as the optimization's start point.

For the JMAPSPPE and the MEE, on the other hand, the initial guess for each optimization was obtained using only the measured data. As in all examples the measurements correspond to the noise-contaminated $z$ state, its guessed path was obtained with a least-squares spline approximation of the measurements. The guess for the $x$ state path, which is the derivative of the $z$ path in all examples, was then obtained as the derivative of the fitted spline. Finally, the parameter guess was obtained by a least-squares regression using the spline's second derivative and the guesses for the $z$ and $x$ paths.

The normalized integrated square error (ISE) metric was used to quantitatively evaluate the state-path estimation error:

$$\text{ISE} := \frac{1}{t_\text{f}} \int_0^{t_\text{f}} \Big( |X_t - x(t)|^2 + |Z_t - z(t)|^2 \Big) \, \mathrm{d}t \tag{4.1}$$

where $X$ and $Z$ are the simulated processes and $x$ and $z$ are the estimated state paths. We note, however, that the same qualitative behaviour of the state-path error observed in the examples is also observed with other metrics like the integrated absolute error (IAE).

The first example, presented in Section 4.1.1, is on the Duffing oscillator with Gaussian measurement noise. It is chosen so that the UKF and its PEM are applicable and can be used as a benchmark state and parameter estimator. Then, in Section 4.1.2, we show an example application on the Duffing oscillator with non-Gaussian measurements. We demonstrate how the MAP and minimum energy estimators can be used with heavy-tailed measurement distributions for robust state estimation and system identification in the presence of outlier measurements. Finally, in Section 4.1.3 we show an example on the Holmes–Rand oscillator with quantitized measurements, showing how the MAP and minimum energy estimators can be used to take into account analog to digital conversion in the modeling and estimation.

### 4.1.1 Duffing oscillator with Gaussian measurements

The first simulated application is made on the Duffing oscillator, a benchmark model for modeling nonlinear dynamics and chaos (Aguirre and Letellier, 2009, Sec. A.3) and state estimation in SDEs (Ghosh et al., 2008; Khalil et al., 2009; Namdeo and Manohar, 2007). The system has two states and its dynamics is given by the following SDEs:

$$\mathrm{d}X_t = [-AZ_t^3 - BZ_t - DX_t + \gamma\cos(t)]\,\mathrm{d}t + \sigma_{\mathrm{D}}\,\mathrm{d}W_t, \tag{4.2a}$$

$$\mathrm{d}Z_t = X_t\,\mathrm{d}t, \tag{4.2b}$$

where $A$, $B$, and $D$ are parameters considered unknown, to be estimated; $\gamma$ and $\sigma_{\mathrm{D}}$ are parameters considered known; and $W$ is a Wiener process representing the process noise. The initial states $X_0$ and $Z_0$ were drawn from independent normal distributions with zero mean and standard deviations $\sigma_{\mathrm{x}}$ and $\sigma_{\mathrm{z}}$ respectively, i.e.,

$$X_0 \sim \mathcal{N}(0, \sigma_{\mathrm{x}}^2)\,, \qquad Z_0 \sim \mathcal{N}(0, \sigma_{\mathrm{z}}^2)\,. \tag{4.3}$$

The total duration of the experiment $t_{\mathrm{f}}$ was varied across experiments, taking the values 50, 100, and 200. The nominal values of the parameters used in the simulation are presented in Table 4.1. The system exhibits chaos with the parameters at these nominal values and is characterized by a double-well potential. An example of the system's simulated state path is shown in Figure 4.1.

Discrete-time measurements of the $Z$ state, corrupted by independent Gaussian noise, were sampled with period $t_{\mathrm{s}}$. Given $Z$ and $\Theta$, each element $Y_k$ of the measurement vector $Y := [Y_0, \dots, Y_N]$ was drawn independentelly from a Gaussian distribution with mean $Z_{kt_{\mathrm{s}}}$ and standard deviation given by the unknown parameter $\Sigma_{\mathrm{y}}$, i.e.,

$$Y_k | Z, \Theta \sim \mathcal{N}\left(Z_{kt_{\mathrm{s}}}, \Sigma_{\mathrm{y}}^2\right).$$

For the estimation, we used the same dynamic and measurement model that was used for generating the data. For the drift parameters, non-informative Gaussian priors were used for the estimation, with zero mean and a large standard deviation:

$$A \sim \mathcal{N}(0, \sigma_{\theta}^2)\,, \qquad B \sim \mathcal{N}(0, \sigma_{\theta}^2)\,, \qquad D \sim \mathcal{N}(0, \sigma_{\theta}^2)\,. \tag{4.4}$$

Table 4.1: Nominal parameter values for the Duffing oscillator experiment

| symbol | $A$ | $B$ | $D$ | $\Sigma_{\mathrm{y}}$ | $\gamma$ | $\sigma_{\mathrm{D}}$ | $\sigma_{\mathrm{x}}$ | $\sigma_{\mathrm{z}}$ | $\sigma_{\theta}$ | $r$ | $s$ | $t_{\mathrm{s}}$ |
|---|---|---|---|---|---|---|---|---|---|---|---|---|
| value | 1.0 | -1.0 | 0.2 | 0.1 | 0.3 | 0.1 | 0.4 | 0.4 | 10 | 1.1 | 10 | 0.1 |

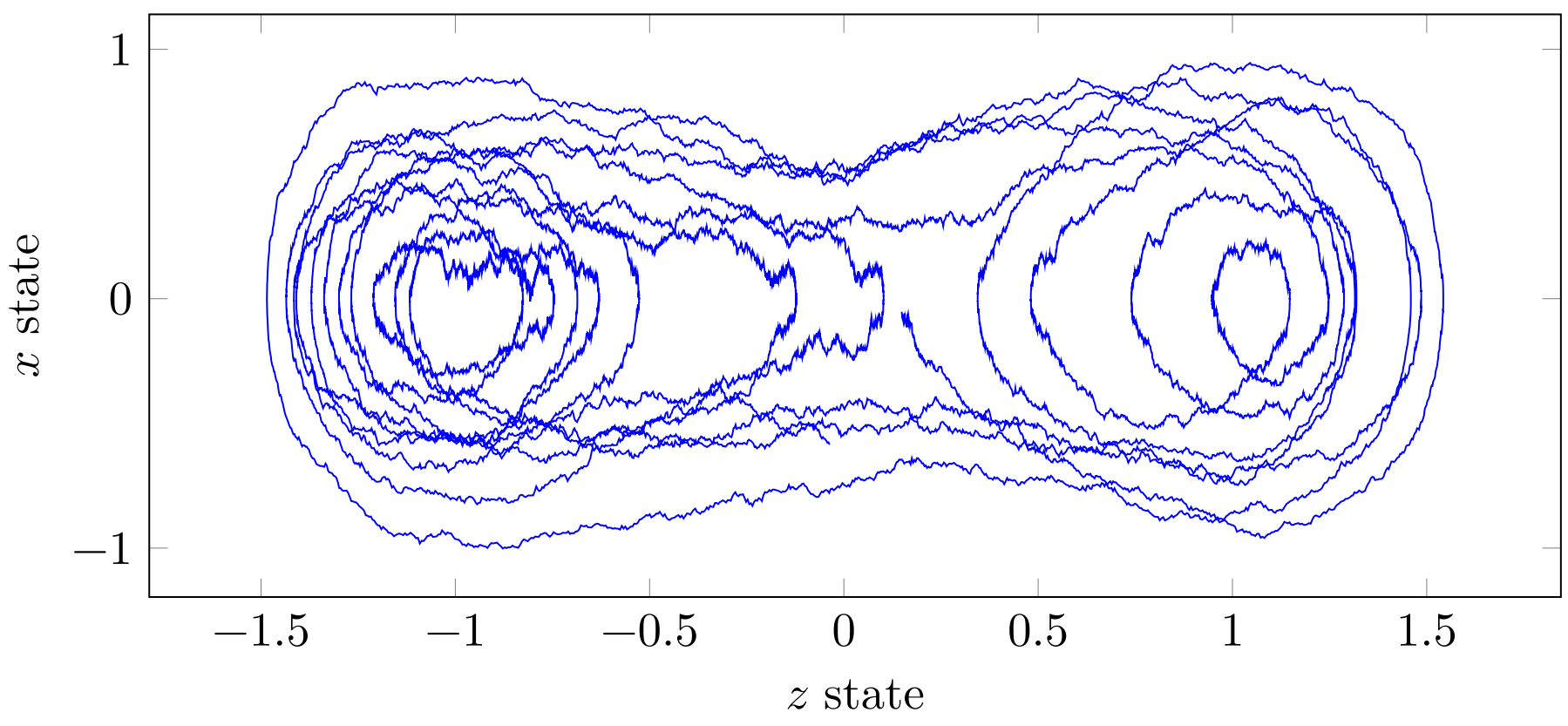

Figure 4.1: Example sample path of the Duffing oscillator, in the state space, with the double well visible.

For the measurement parameter $\varSigma_{\mathrm{y}}$, the gamma distribution with shape $r$ and a large scale parameter $s$ was used for the prior:

$$\varSigma_{\mathrm{y}} \sim \varGamma(r, s). \tag{4.5}$$

Putting the estimation model in the format of Chapter 2, we have that it is characterized by the unknown parameter vector $\theta := [a, b, d, \sigma_{\mathrm{y}}]$, the drift functions

$$f(t, x, z, \theta) := -az^3 - bz - dx + \gamma \cos(t), \qquad h(t, x, z, \theta) := x, \tag{4.6}$$

the diffusion matrix $\boldsymbol{G} := [\sigma_{\mathrm{D}}]$, the prior log-density

$$\ln \pi(x, z, \theta) := -\frac{1}{2}\left(\frac{x^2}{\sigma_{\mathrm{x}}^2} + \frac{z^2}{\sigma_{\mathrm{z}}^2} + \frac{a^2 + b^2 + d^2}{\sigma_{\theta}^2}\right) + (r-1)\ln \sigma_{\mathrm{y}} - \frac{1}{s}\sigma_{\mathrm{y}}, \tag{4.7}$$

where the constant terms are omitted as they do not influence the location of maxima, the measurement space $\mathcal{Y} := \mathbb{R}^{N+1}$, and the measurement log-likelihood

$$\ln \psi(y \,|\, x, z, \theta) := -\frac{1}{2}\sum_{k=0}^{N} \frac{\big(y_k - z(kt_{\mathrm{s}})\big)^2}{\sigma_{\mathrm{y}}^2} - (N+1)\ln \sigma_{\mathrm{y}},$$

where the constant terms have been omitted as well. Additionally, we note that the measure $\nu$ with respect to which the density $\psi$ is defined is the Lebesgue measure over $\mathcal{Y}$.

A total of 100 Monte Carlo simulations were performed for each value of the total experiment duration $t_{\mathrm{f}}$. The simulated and estimated state paths for one of the simulations are shown in Figure 4.2 for the whole experiment interval and in Figure 4.3 for a portion of it, so that finer details can be noticed. For this experiment, the approximations of nonlinear Kalman filters are reasonable and the PEM is an adequate ground truth for evaluation of the MAP and minimum energy estimates. It can be seen that all estimated state-paths are fairly close and that their overall error is comparable. This can also be seen in Figure 4.4, which shows boxplots of the state-path estimation error, quantified by the ISE metric defined in (4.1). It can be seen that the state estimation error of all three methods is statistically comparable. We note that, although not shown here, the same qualitative behaviour of Figure 4.4 is also observed with the integrated absolute error (IAE) metric.

Although the state-path estimates of all three methods are very similar, however, their parameter estimates are not. Boxplots of the parameter estimates obtained with all three methods are shown in Figure 4.5. We can see that the minimum energy $D$ parameter estimates are consistently lower than the MAP and PEM estimates. This can be understood by noting that the divergence of the noisy drift is given by

$$\operatorname{div} f(t, x, z, \theta) = -d.$$

Consequently, the Onsager–Machlup functional favors higher $d$ values, as they attenuate the fluctuations due to the process noise. The energy functional, on the other hand, favors lower $d$ values, as noise paths with less energy are needed to maintain the same state path. This is also evident from a physical interpretation of the parameters, as $d$ is the oscillators' damping constant and is proportional to the rate at which the unforced system loses energy (Kanamaru, 2008). The estimates for the remaining parameters drift parameters $A$ and $B$ with all methods were statistically very similar. Both the MEE and JMAPSPPE also have a small bias for the $\sigma_{\mathrm{y}}$ parameter, while the PEM estimator displays no bias.

Another result of this experiment was that, for the drift parameters, the MAP estimates obtained with the PEM were very close to the estimates of the JMAPSPPE, that is, the marginal and joint modes were close, especially for the longer experiment. We believe this qualitative behaviour should occur for most systems for which the Kalman filter is applicable, i.e., systems subject to Gaussian noise and with frequent measurements, for which the state posterior is approximately Gaussian. For these systems, the joint MAP state-path and parameter estimator could be used as a replacement for Kalman-filter based prediction error methods with similar results, similarly to what was proposed by Varziri et al. (2008b). The MAP and minimum energy estimators are applicable for systems with more general measurement distributions, however, which we illustrate with the next examples.

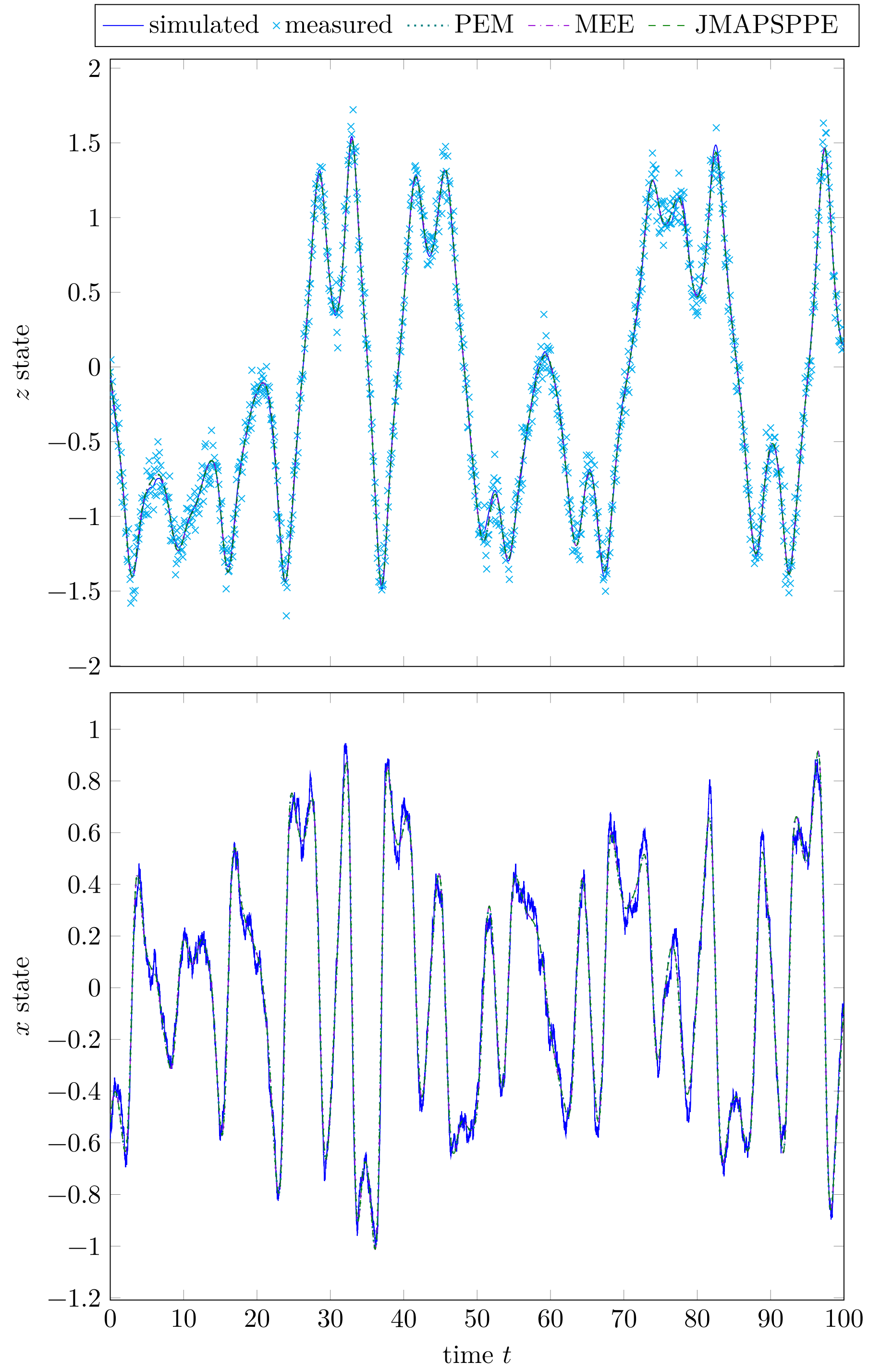

Figure 4.2: Simulated and estimated state paths for one of the Monte Carlo simulations of the Duffing oscillator, with $t_\mathrm{f} = 100$.

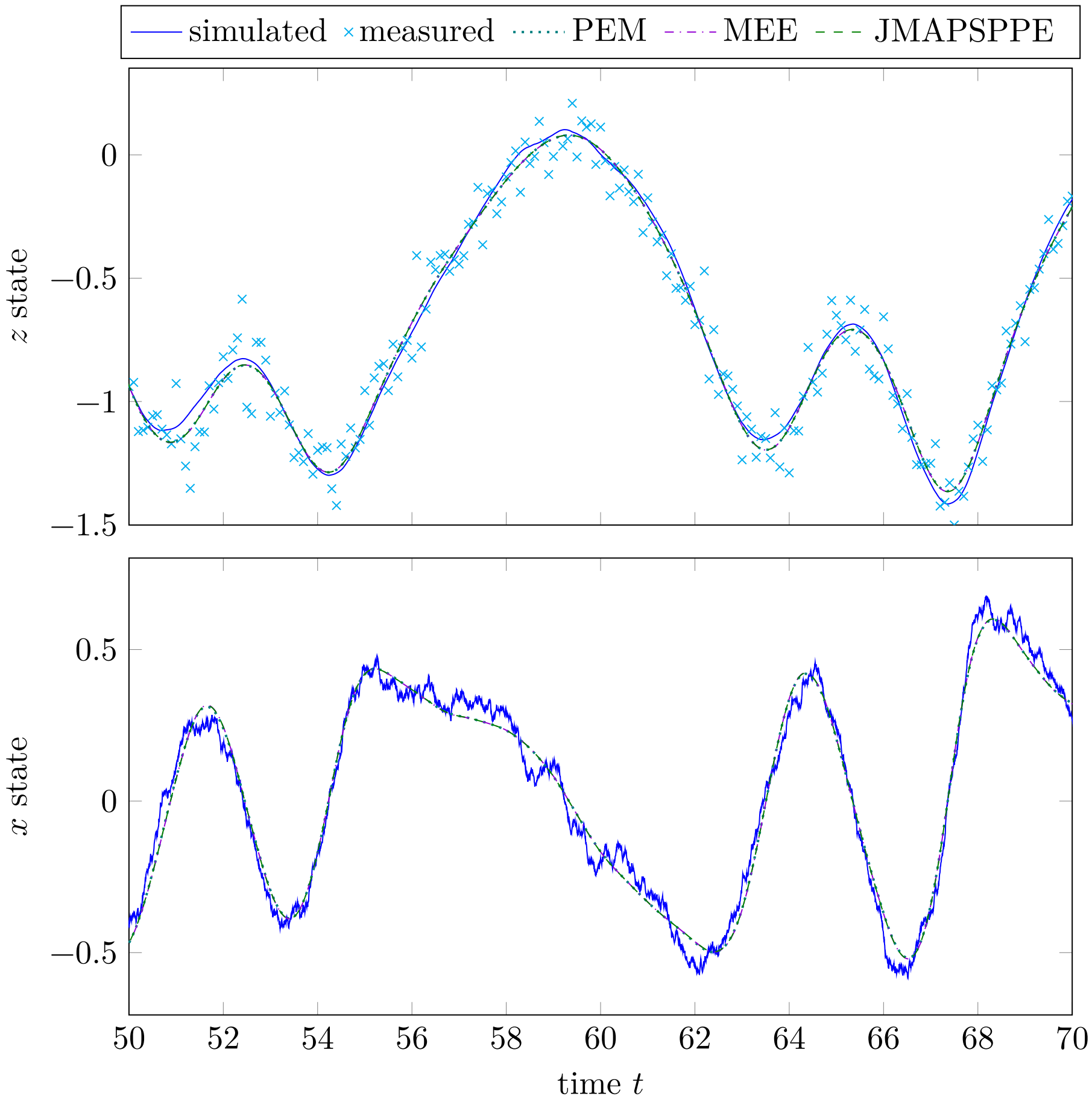

Figure 4.3: Detail of the simulated and estimated state paths for one of the Monte Carlo simulations of the Duffing oscillator, with $t_\mathrm{f} = 100$.

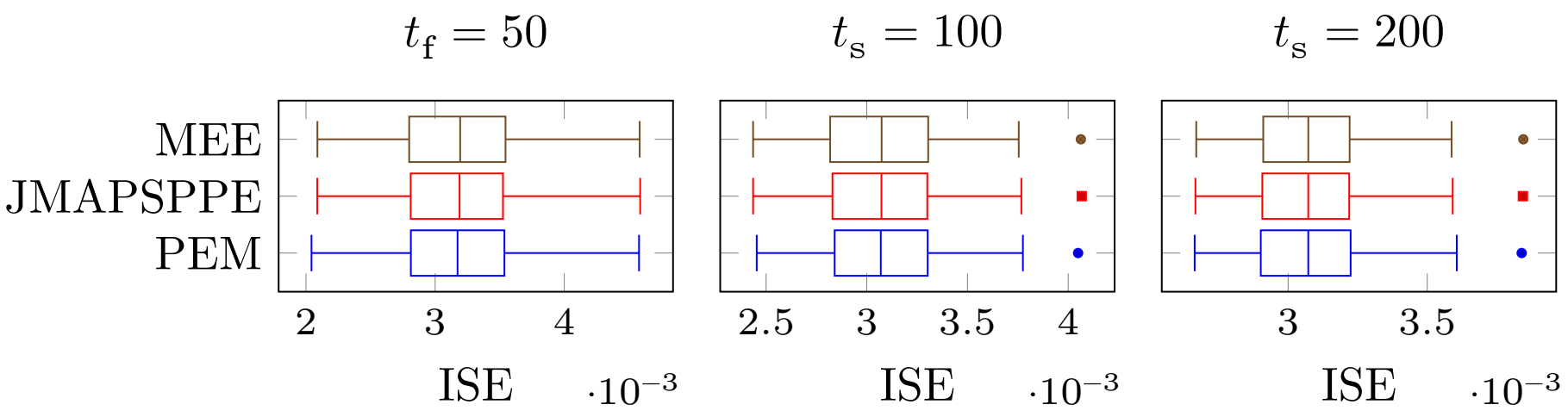

Figure 4.4: Boxplots of the integrated square error (ISE) of the Duffing oscillator estimated state paths over all Monte Carlo simulations.

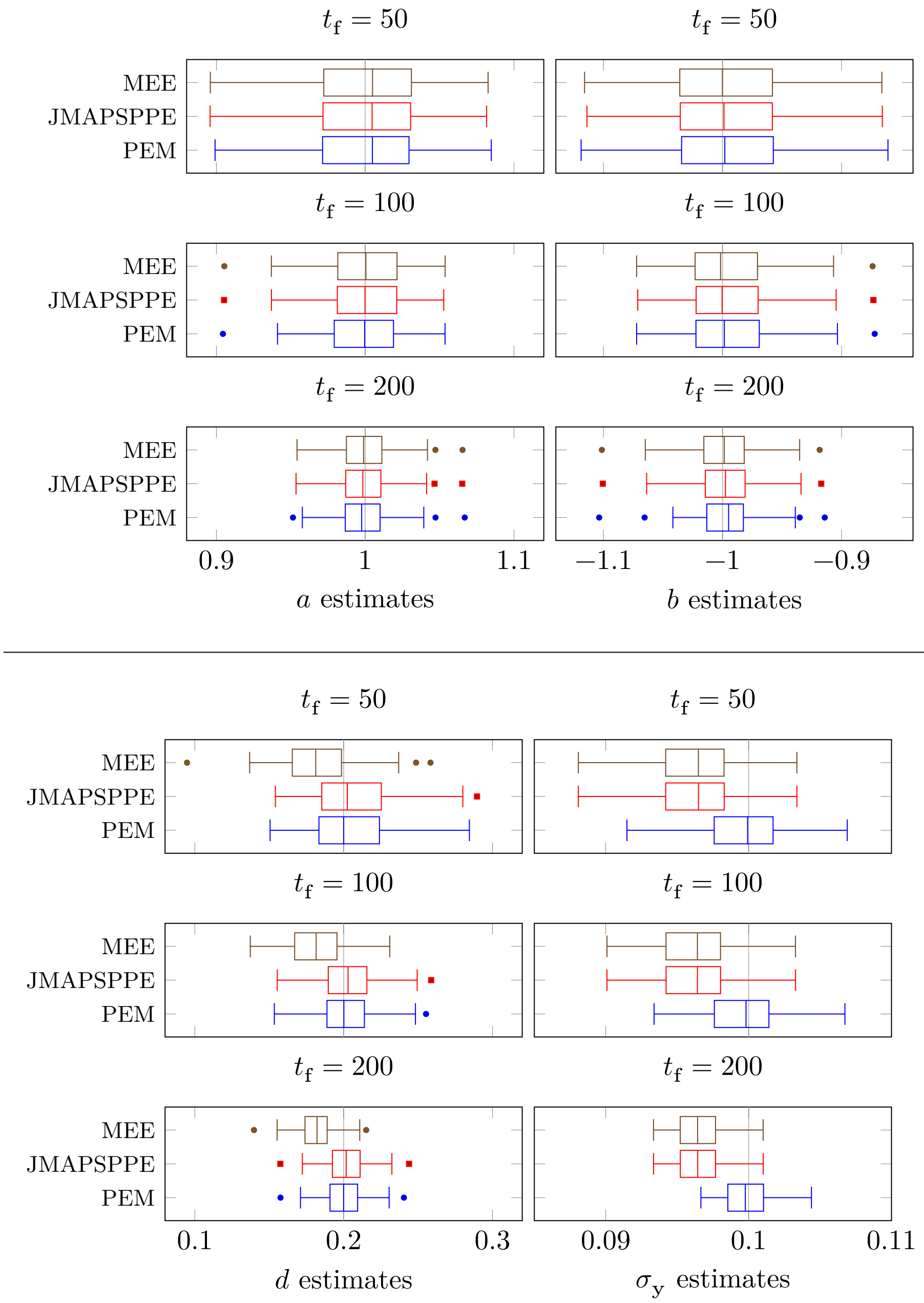

Figure 4.5: Boxplot of the Duffing oscillator parameter estimates. The nominal parameter values, $a = 1$, $b = -1$, $d = 0.2$, and $\sigma_y = 0.1$, are marked with gridlines in the plots.

### 4.1.2 Duffing oscillator with outlier measurements

For the second simulated example we used the Duffing oscillator, as in the previous example, but with measurements containing outliers, as in the examples of Aravkin et al. (2011, 2012c) and Dutra et al. (2014). The system's dynamics are given by the SDE (4.2) and its initial states sampled according to (4.3). The total simulation length was $t_\mathrm{f} = 100$.

Discrete-time measurements of the $Z$ state, corrupted by independent noise, were sampled with period $t_\mathrm{s}$. Given $Z$, each element $Y_k$ of the measurement vector $Y := [Y_0, \dots, Y_N]$ was drawn independentelly from the Gaussian mixture distribution

$$Y_k|Z \sim p_\mathrm{o}\mathcal{N}\left(Z_{kt_\mathrm{s}}, \sigma_\mathrm{o}^2\right) + (1 - p_\mathrm{o})\mathcal{N}\left(Z_{kt_\mathrm{s}}, \sigma_\mathrm{r}^2\right),$$

where $\sigma_\mathrm{r}, \sigma_\mathrm{o} \in \mathbb{R}_{>0}$ are the regular measurements' and outliers' standard deviation, respectively, and $p_\mathrm{o} \in \mathbb{R}_{>0}$ is the outlier probability. The outlier probability $p_\mathrm{o}$ was varied to investigate different outlier noise contamination scenarios.

As in Section 4.1.1, for the estimation the unknown parameter vector was given by $\theta := [a, b, d, \sigma_\mathrm{y}]$. The estimators used the same dynamic model that was used to generate the data, with the $f$ and $h$ functions given by (4.6) and the diffusion matrix $\boldsymbol{G} := [\sigma_\mathrm{D}]$. In addition, the same initial state and parameter priors of the previous section were used as well, with distributions given by (4.4)–(4.5) and log-density given by (4.7).

A different measurement model was used in the estimation, however, to investigate robustness of the estimator against outliers. Like in (Aravkin et al., 2012c; Dutra et al., 2014), Student's $t$-distribution with 4 degrees of freedom was used as the measurement distribution, with the following expression for its log-likelihood:

$$\ln \psi(y \,|\, x, z, \theta) := -\frac{5}{2}\sum_{k=0}^{N} \ln\left(1 + \frac{\big(y_k - z(kt_\mathrm{s})\big)^2}{4\sigma_\mathrm{y}^2}\right) - \ln \sigma_\mathrm{y},$$

where the constant terms, which do not influence the location of maxima, have been omitted and the unknown parameter $\sigma_\mathrm{y}$ is the measurement noise scale, to be estimated. Additionally, we note that the measure $\nu$ with respect to which the density $\psi$ is taken is the Lebesgue measure over $\mathcal{Y} := \mathbb{R}^{N+1}$.

The nominal values of the parameters used in the simulation are presented in Table 4.2. The simulated and estimated state paths for one of these simulations are shown in Figure 4.6 for the whole experiment interval and in Figure 4.7 for a portion of it. For each tested value of $p_\mathrm{o}$, a total of 100 Monte Carlo simulations were performed.

Boxplots of the state-path estimation error, quantified by the ISE metric defined in (4.1), are shown in Figure 4.8. We can see that the state-path estimation error of both the JMAPSPPE and the MEE is comparable on all tested values of $p_o$. The state-path estimation error of the PEM, however, is significantly larger. This is because the outliers make the measurement distribution heavy-tailed and, consequently, better represented by Student's $t$-distribution than by the normal distribution.

Table 4.2: Nominal parameter values for the Duffing oscillator experiment with outliers.

| symbol | $A$ | $B$ | $D$ | $\gamma$ | $\sigma_{\mathrm{D}}$ | $\sigma_{\mathrm{x}}$ | $\sigma_{\mathrm{z}}$ | $\sigma_{\theta}$ | $\sigma_{\mathrm{o}}$ | $\sigma_{\mathrm{r}}$ | $r$ | $s$ | $t_{\mathrm{s}}$ |
|---|---|---|---|---|---|---|---|---|---|---|---|---|---|
| value | 1 | -1 | 0.2 | 0.3 | 0.1 | 0.4 | 0.4 | 10 | 1 | 0.2 | 1.1 | 10 | 0.1 |

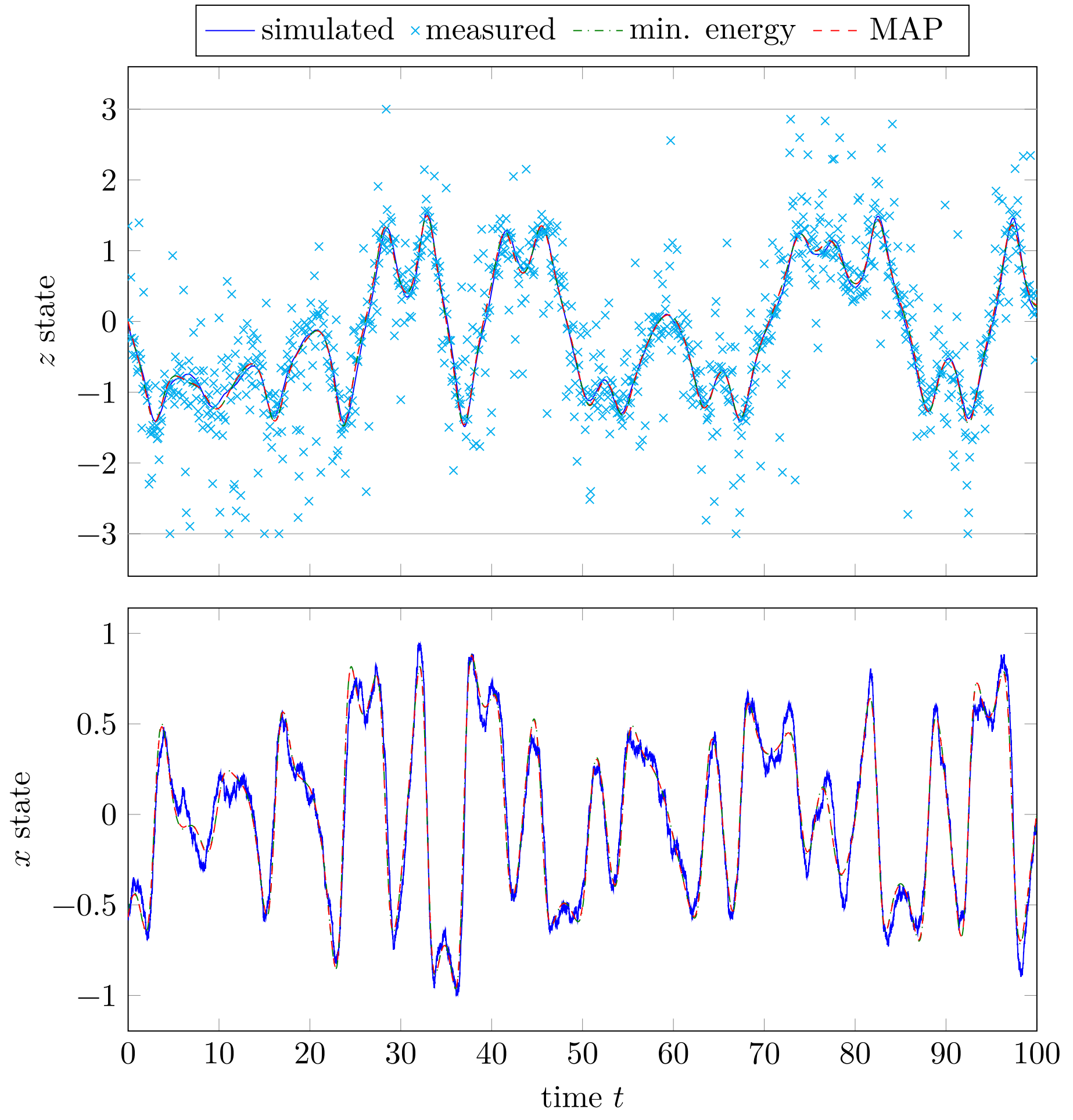

Figure 4.6: Simulated and estimated state paths for one of the Monte Carlo simulations of the Duffing oscillator with outliers, for $p_{\mathrm{o}} = 0.4$. The measurements outside the plot range are shown at $\pm 3$.

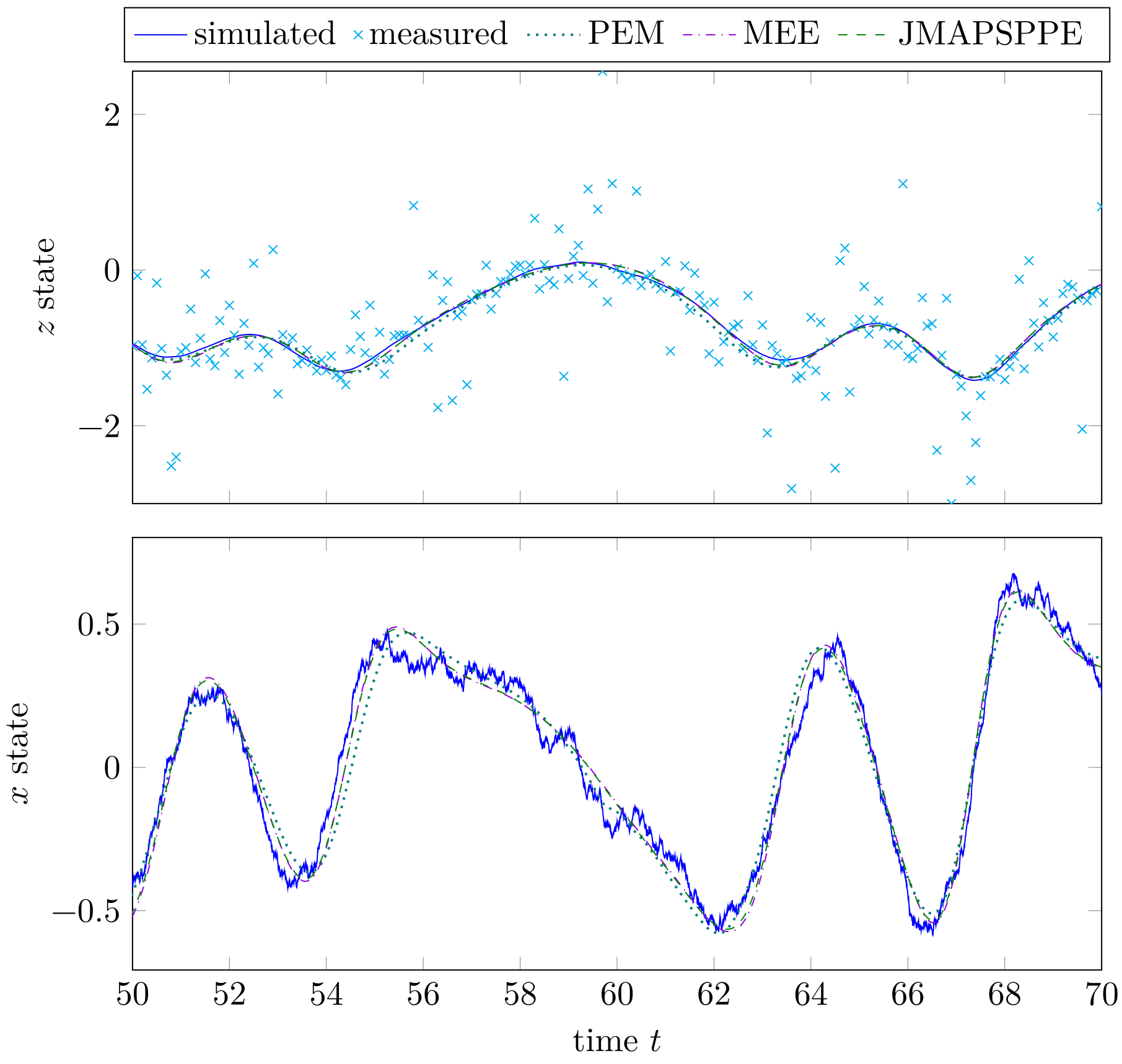

Figure 4.7: Detail of the simulated and estimated state paths for one of the Monte Carlo simulations of the Duffing oscillator experiment with outliers, for $p_o = 0.4$.

Boxplots of the parameter estimates of all three methods are shown in Figure 4.9. Like in the previous example, it can be seen that the MEE is clearly biased in the $d$ parameter, due to the fact that it ignores the amplification of the noise by the drift.

In this example, we have shown that the JMAPSPPE and MEE can be used for robust state-path and parameter estimation in systems with measurements under intense outlier contamination for which nonlinear Kalman smoothers and the prediction error method yields larger errors. Robustness against outliers was gained by modeling the measurements with heavy-tailed distributions. In the next example, we show how modeling the measurements with slender-tailed distributions can be used extract more information from the data by taking into account the analog to digital conversion.

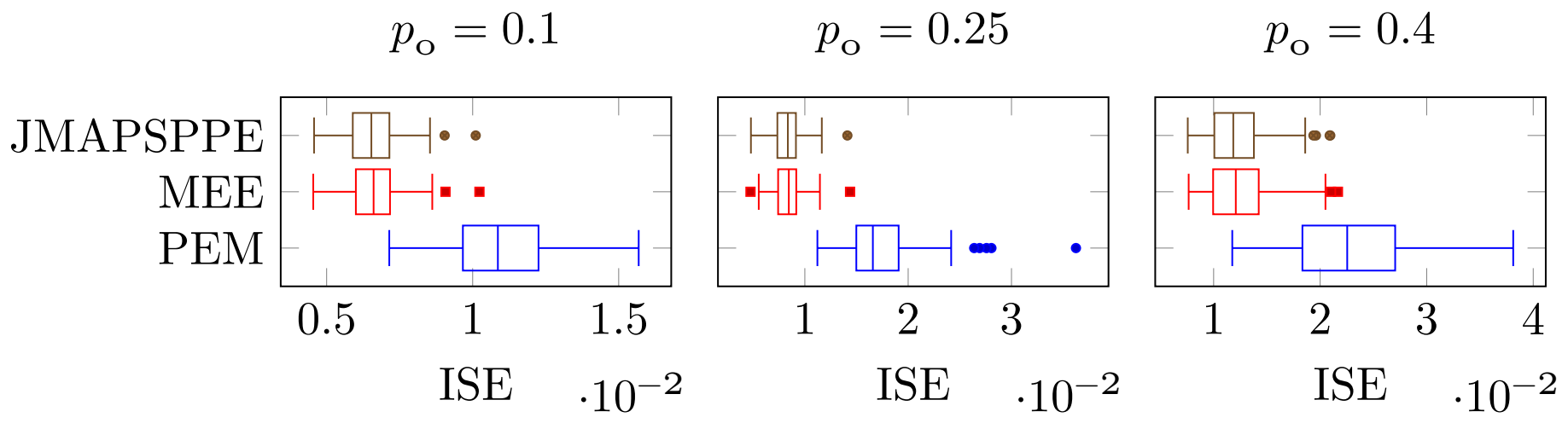

Figure 4.8: Boxplots of the state-path integrated square error (ISE) for the Duffing oscillator with outlier measurements.

### 4.1.3 Holmes–Rand oscillator with quantitized measurements

The final simulated application is made on the Holmes–Rand oscillator, a general nonlinear system which includes both the Duffing and Van der Pol oscillators as special cases (Holmes and Rand, 1980). The system has two states and its dynamics is given by the following SDEs:

$$\begin{aligned} \mathrm{d}X_t &= [-(A + \Gamma Z_t^2)X_t - BZ_t - DZ_t^3 + \phi\cos(t)]\,\mathrm{d}t + \sigma_\mathrm{D}\,\mathrm{d}W_t, \\ \mathrm{d}Z_t &= X_t\,\mathrm{d}t, \end{aligned}$$

where $A$, $\Gamma$, $B$, and $D$ are parameters considered unknown, to be estimated; $\phi$ and $\sigma_\mathrm{D}$ are parameters considered known; and $W$ is a Wiener process representing the process noise. The initial states $X_0$ and $Z_0$ were drawn from independent normal distributions with zero mean and standard deviations $\sigma_\mathrm{x}$ and $\sigma_\mathrm{z}$ respectively, i.e.,

$$X_0 \sim \mathcal{N}(0, \sigma_\mathrm{x}^2)\,, \qquad Z_0 \sim \mathcal{N}(0, \sigma_\mathrm{z}^2)\,.$$

The nominal parameter values used in the simulation are presented in Table 4.3. The system exhibits chaos with these parameters values and is characterized by a double-well potential. An example of the system's simulated state path is shown in Figure 4.10.

Discrete-time measurements of the $Z$ state corrupted by independent noise were sampled at regular time intervals. A sampling period $t_\mathrm{s} = 0.1$ was chosen and a total of $N+1 := 501$ measurements were taken in each realization of the simulation. To further emulate the effect of digital data aquisition, each element $Y_k$ of the measurement vector $Y := [Y_0, \ldots, Y_N]$ was drawn independentelly, given $Z$ and $\Theta$, from a the Gaussian mixture distribution with mean $Z_{kt_{\mathrm{s}}}$ and standard deviation $\Sigma_{\mathrm{y}}$ and then rounded towards the nearest multiple of the bit length $l_{\mathrm{b}}$. The nominal value of the unknown parameter $\Sigma_{\mathrm{y}}$ was then varied to evaluate different analog to digital conversion scenarios. The measurement probabilities of this model are illustrated graphically in Figure 4.11. The probability mass function and likelihood function of the measurements, for different $l_{\mathrm{b}}/\sigma_{\mathrm{y}}$ ratios, are also shown in Figures 4.12 and 4.13, respectively.

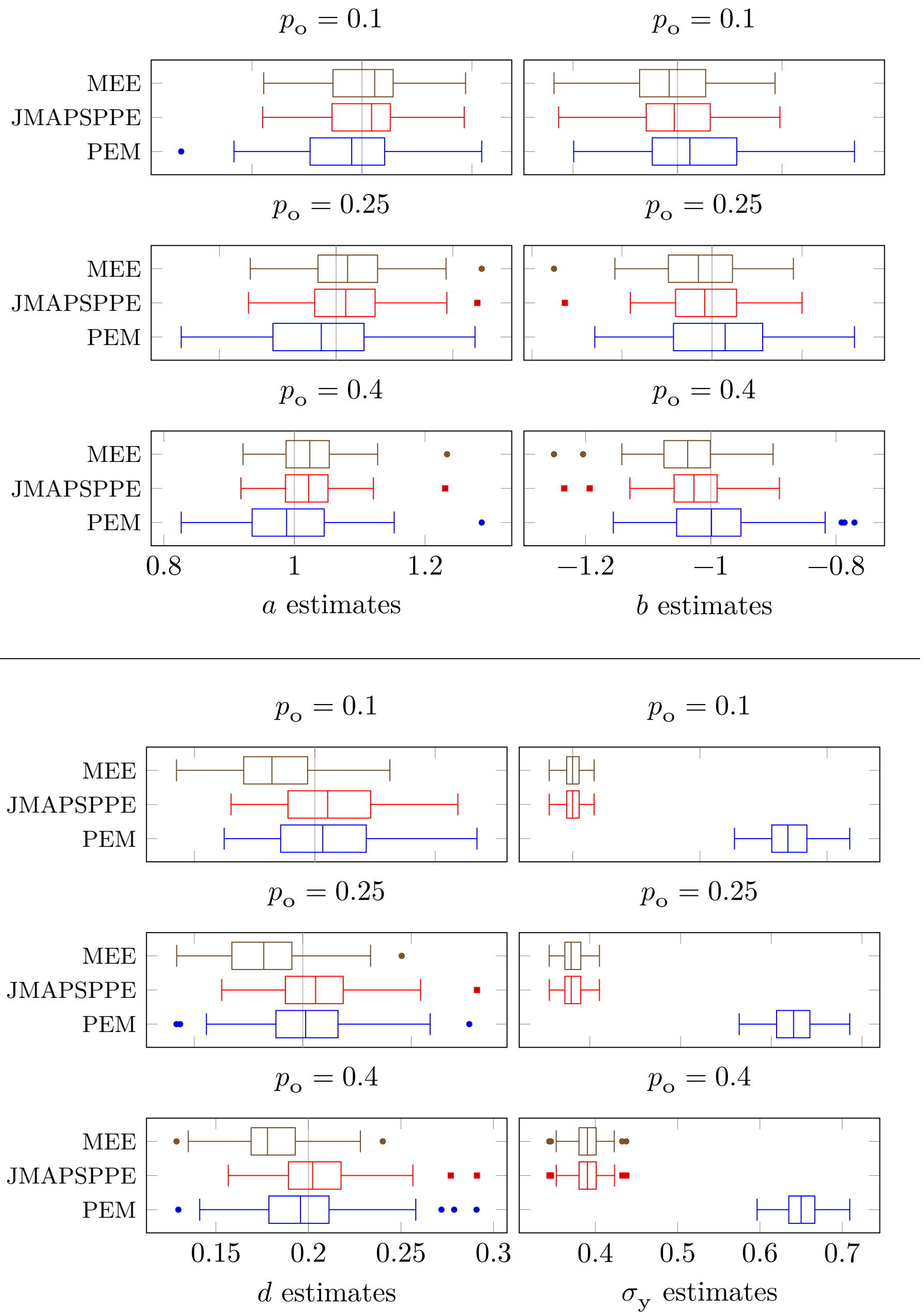

Figure 4.9: Boxplot of the Duffing oscillator parameter estimates with outlier measurements. The nominal parameter values, $a = 1$, $b = -1$ and $d = 0.2$, are marked with gridlines in the plots.

Table 4.3: Nominal parameter values for the Holmes–Rand oscillator experiment.

| symbol | $A$ | $B$ | $\Gamma$ | $D$ | $\phi$ | $\sigma_{\mathrm{D}}$ | $\sigma_{\mathrm{x}}$ | $\sigma_{\mathrm{z}}$ | $\sigma_{\theta}$ | $r$ | $l_{\mathrm{b}}$ | $t_{\mathrm{f}}$ |
|---|---|---|---|---|---|---|---|---|---|---|---|---|
| value | 0.2 | -1.0 | 0.2 | 1.0 | 0.4 | 0.1 | 0.1 | 0.1 | 10 | 4 | 0.05 | 50 |

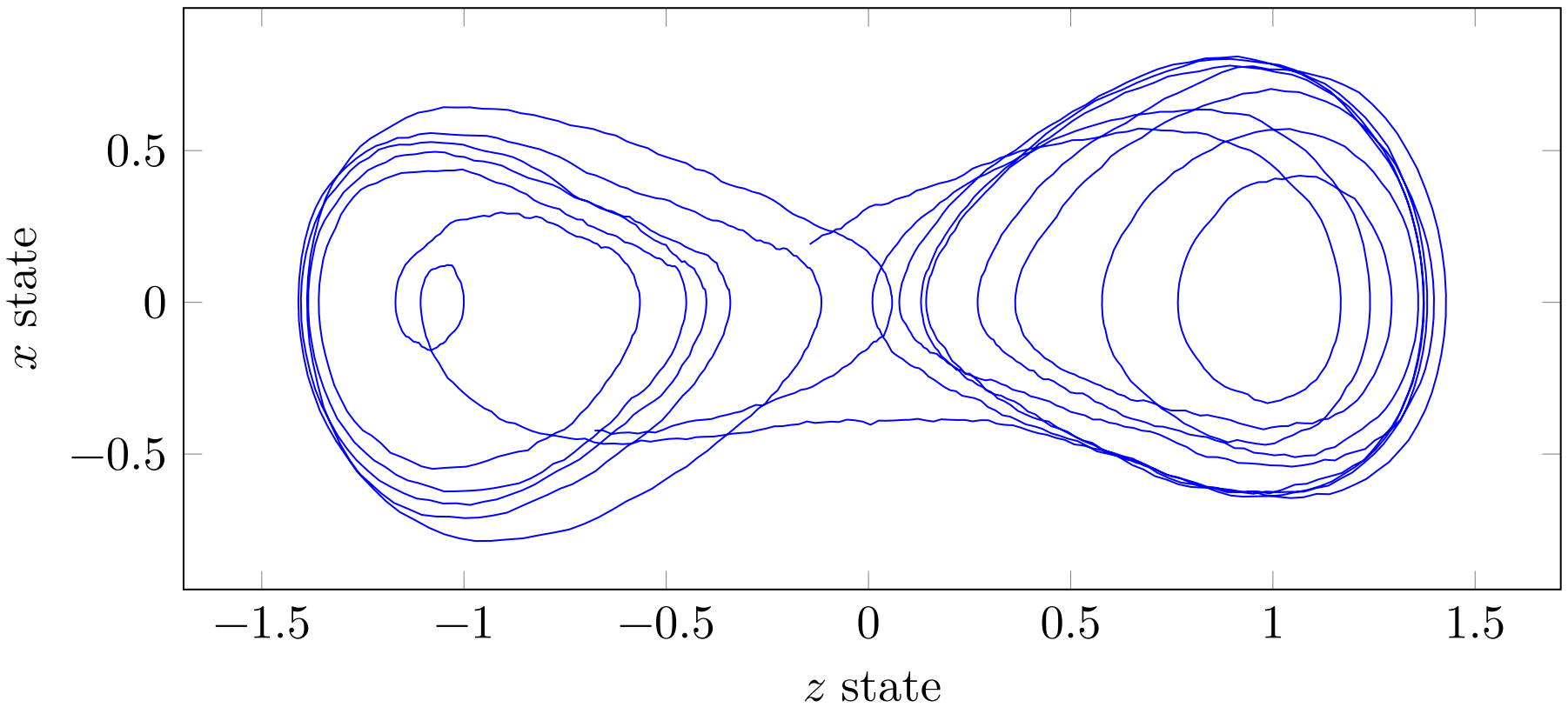

Figure 4.10: Example sample path of the Holmes–Rand oscillator, in the state space, with the double well visible.

For the estimation, we used the same dynamic and measurement model that was used for generating the data. Non-informative Gaussian priors were used for the estimation of the drift parameters, with zero mean and a large standard deviation $\sigma_{\theta}$:

$$A \sim \mathcal{N}(0, \sigma_{\theta}^2)\,, \qquad B \sim \mathcal{N}(0, \sigma_{\theta}^2)\,, \qquad \Gamma \sim \mathcal{N}(0, \sigma_{\theta}^2) \qquad D \sim \mathcal{N}(0, \sigma_{\theta}^2)\,.$$

For the measurement parameter $\Sigma_{\mathrm{y}}$, the gamma distribution with shape $r$ and scale $\frac{l_{\mathrm{b}}}{3}$ was used for the prior:

$$\Sigma_{\mathrm{y}} \sim \Gamma\!\left(r, \tfrac{l_{\mathrm{b}}}{3}\right).$$

Putting the estimation model in the format of Chapter 2, we have that it is characterized by the unknown parameter vector $\theta := [a, b, \gamma, d, \sigma_{\mathrm{y}}]$, the drift

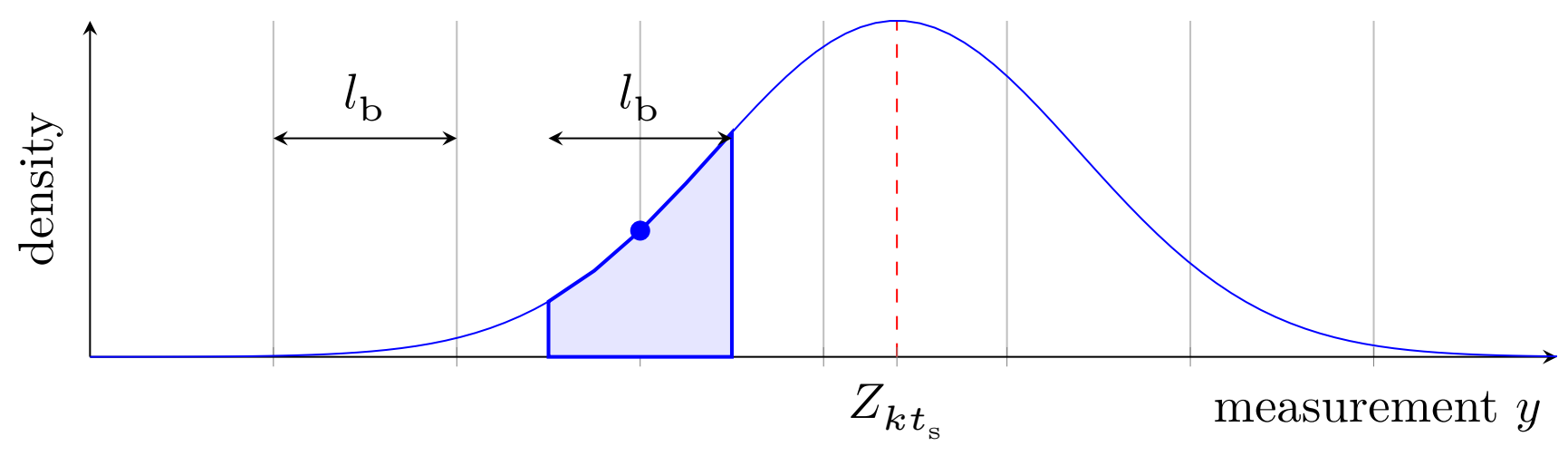

Figure 4.11: Graphical illustration of the measurement model used in the Holmes–Rand experiment. The dashed line indicates the simulated value of $Z$ at the corresponding measurement instant and the vertical gridlines indicate the integer multiples of the bit length $l_\mathrm{b}$. The probability of the outcome shown as a mark on the plot is the area under the curve shown as the shaded region.

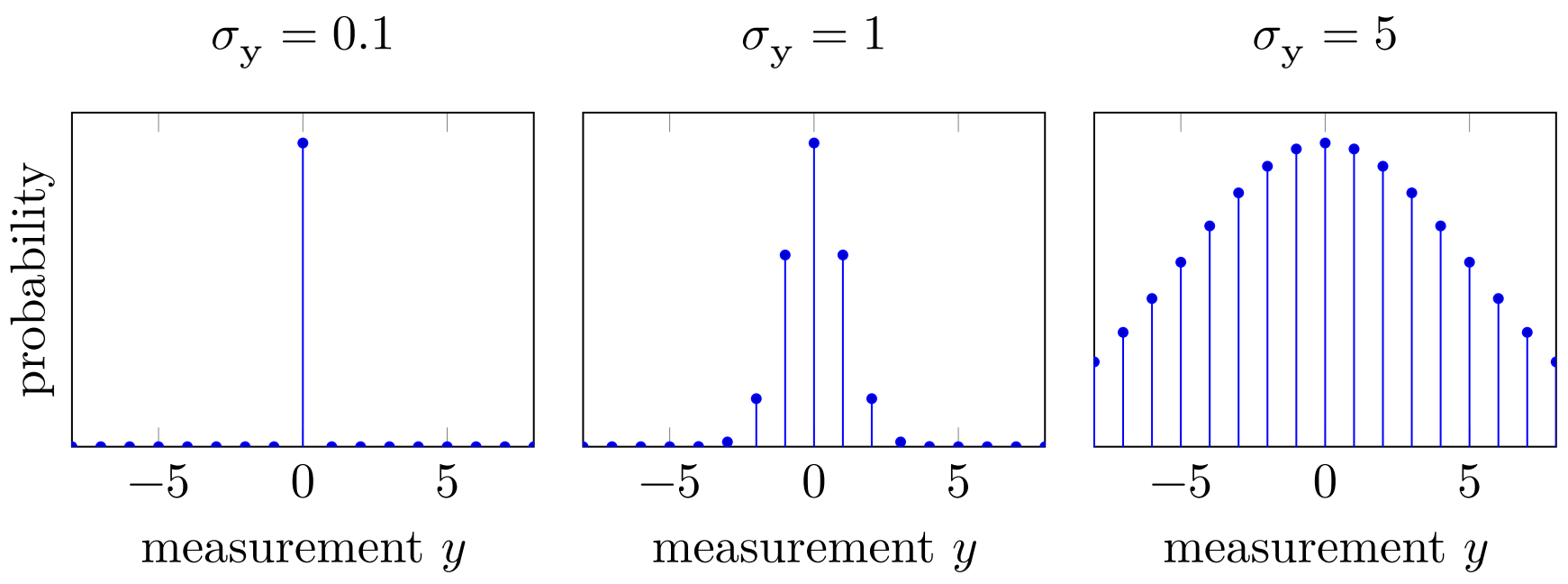

Figure 4.12: Conditional probability mass functions of the Holmes–Rand measurement $Y_k$, given $Z_{kt_\mathrm{s}} = 0$, $l_\mathrm{b} = 1$, and different values of $\sigma_\mathrm{y}$.

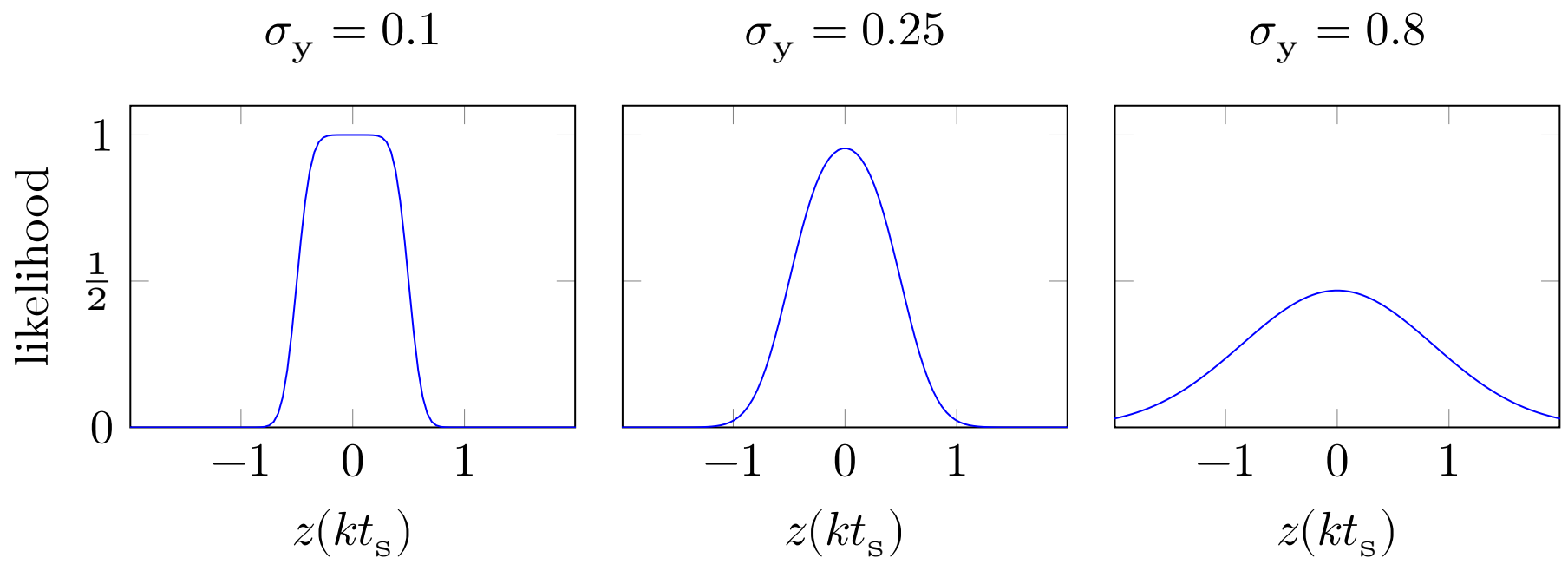

Figure 4.13: Likelihood functions of the Holmes–Rand measurement $Y_k = 0$, given $Z_{kt_\mathrm{s}}$, for $l_\mathrm{b} = 1$ and different values of $\sigma_\mathrm{y}$.

functions

$$f(t,x,z,\theta) := -(a+\gamma z^2)x - bz - dz^3 + \phi\cos(t), \tag{4.8a}$$
$$h(t,x,z,\theta) := x, \tag{4.8b}$$

the diffusion matrix $\boldsymbol{G} := [\sigma_{\mathrm{D}}]$, the prior log-density

$$\ln\pi(x,z,\theta) := -\frac{1}{2}\left(\frac{x^2}{\sigma_{\mathrm{x}}^2} + \frac{z^2}{\sigma_{\mathrm{z}}^2} + \frac{a^2+b^2+\gamma^2+d^2}{\sigma_{\theta}^2}\right) + (r-1)\ln\sigma_{\mathrm{y}} - \frac{3}{l_{\mathrm{b}}}\sigma_{\mathrm{y}},$$

where the constant terms are omitted as they do not influence the location of maxima, the measurement space $\mathcal{Y} := \mathbb{R}^{N+1}$, and the measurement log-likelihood

$$\ln\psi(y\,|\,x,z,\theta) := \sum_{k=0}^{N}\ln\left(\Phi\left(\tfrac{z(kt_{\mathrm{s}})-y_k+\frac{1}{2}l_{\mathrm{b}}}{\sigma_{\mathrm{y}}}\right) - \Phi\left(\tfrac{z(kt_{\mathrm{s}})-y_k-\frac{1}{2}l_{\mathrm{b}}}{\sigma_{\mathrm{y}}}\right)\right),$$

where $\Phi$ is the standard normal cumulative distribution function. Additionally, we note that the measure $\nu$ with respect to which the density $\psi$ is defined is the counting measure over $\mathcal{Y}$, as $Y$ is a discrete random variable.

For each scenario 100 Monte Carlo simulations were performed. As in the previous examples, the ISE metric defined in (4.1) was used to quantify the state estimation error. It can be seen that when $\sigma_{\mathrm{y}}$ is of the same order of magnitude as $l_{\mathrm{b}}$, then the state-path estimation error of all three methods is comparable. However, when the noise standard deviation $\sigma_{\mathrm{y}}$ is much smaller than the bit length $l_{\mathrm{b}}$, then the PEM yields larger errors. This can be understood by analysing Figure 4.13 once more. When $l_{\mathrm{b}} \leq \sigma_{\mathrm{y}}$, then the measurement likelihood is approximatelly Gaussian. When $\sigma_{\mathrm{y}}$ decreases, however, then the measurement likelihood approaches the uniform distribution, which is not well represeted by the Gaussian distribution. The simulated, measured, and estimated processes for one Monte Carlo simulation can be seen in Figure 4.15 for the whole experiment interval and in Figure 4.16 for a portion of it. It can be seen qualitatively that the measurements are indeed not Gaussian. Boxplots of the parameter estimates obtained with all three methods are shown in Figure 4.17. Unlike in the Duffing oscillator, a clear biasing of the MEE for the damping parameters is not observed.

In this example, it can be seen that the general form of the JMAPSPPE and MEE allows them to better model analog to digital conversion. It shows that, overall, the prediction error method attains a higher state-path estimation error when the measurement noise is much smaller than the bit length, as it does not represent well the analog-to-digital conversion.

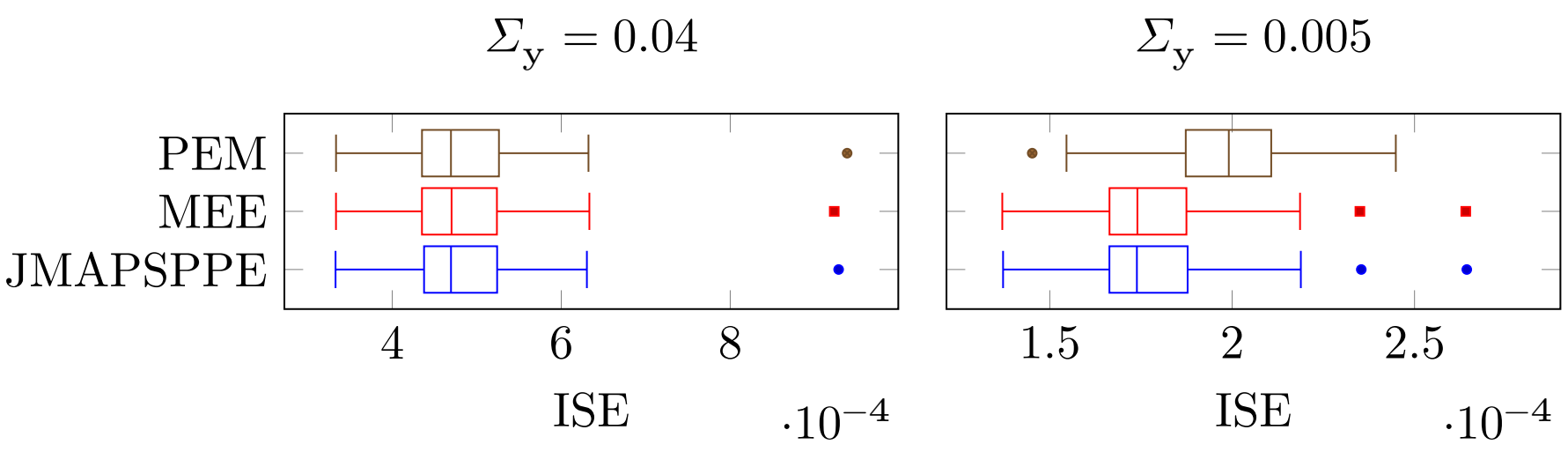

Figure 4.14: Boxplots of the integrated square error (ISE) of the Holmes–Rand oscillator estimated state paths over all Monte Carlo simulations.

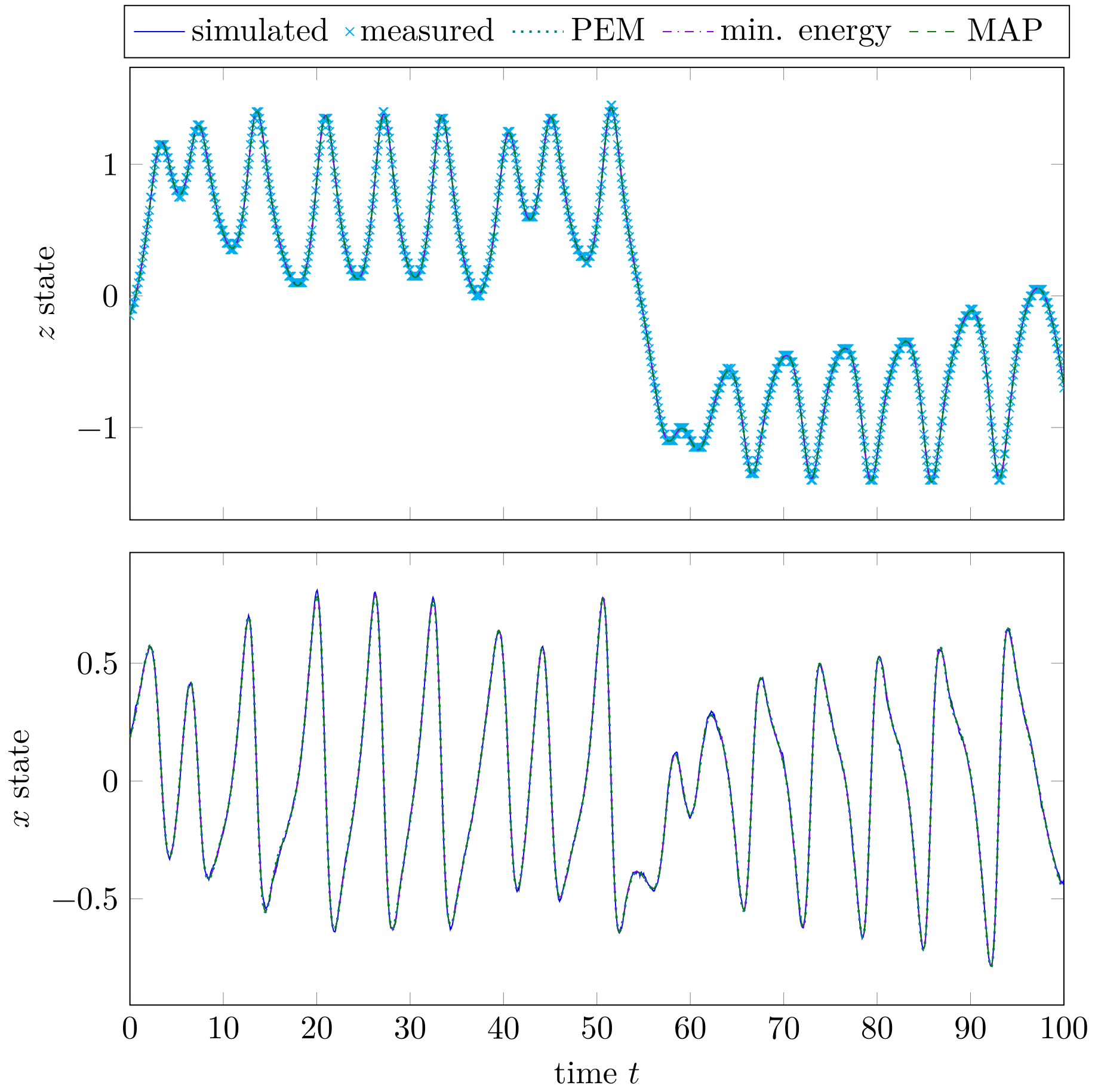

Figure 4.15: Simulated and estimated state paths for one of the Monte Carlo simulations of the Holmes–Rand oscillator, with $\Sigma_y = 0.005$.

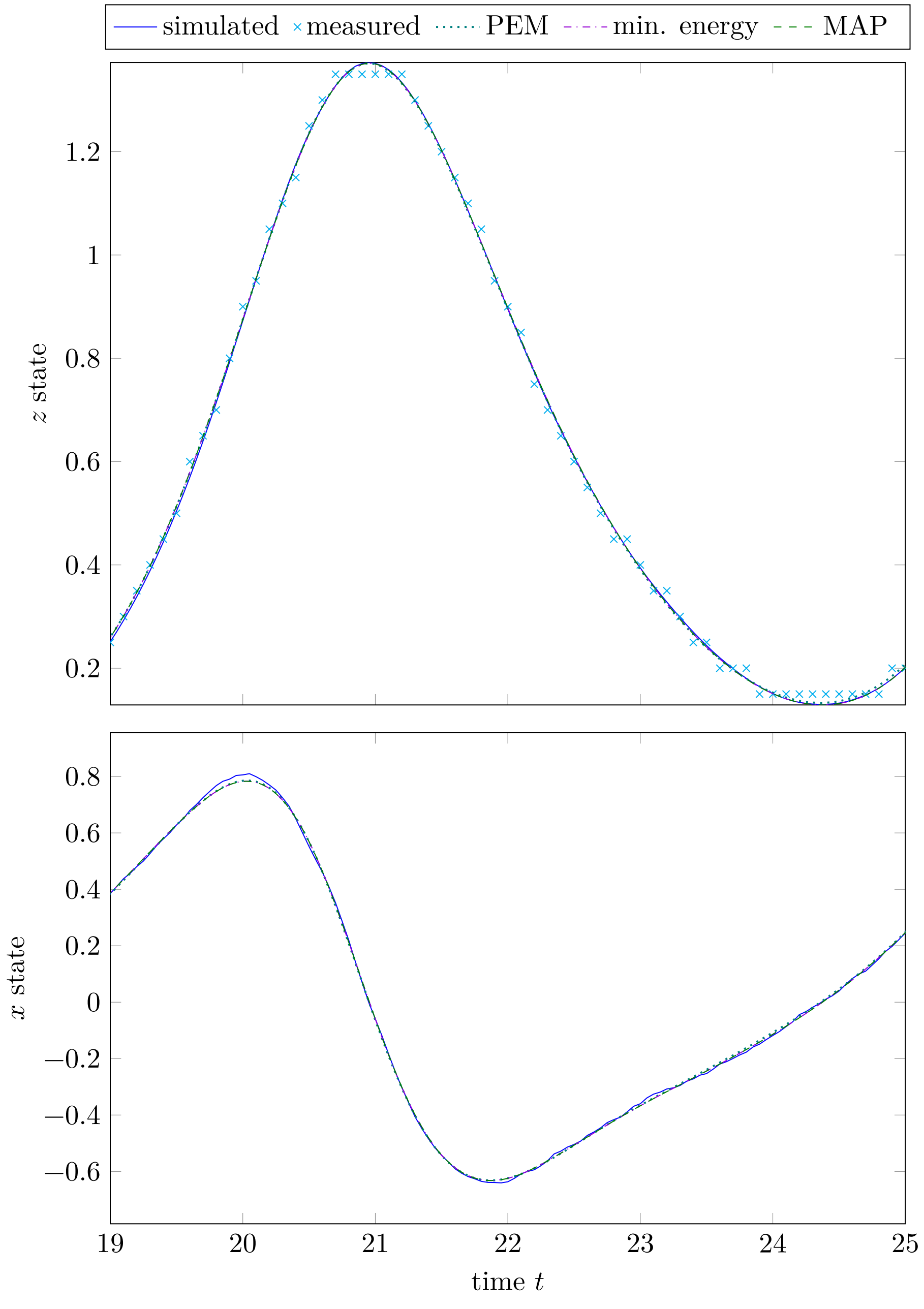

Figure 4.16: Detail of the simulated and estimated state paths for one of the Monte Carlo simulations of the Holmes–Rand oscillator, with $\Sigma_{\mathrm{y}} = 0.005$.

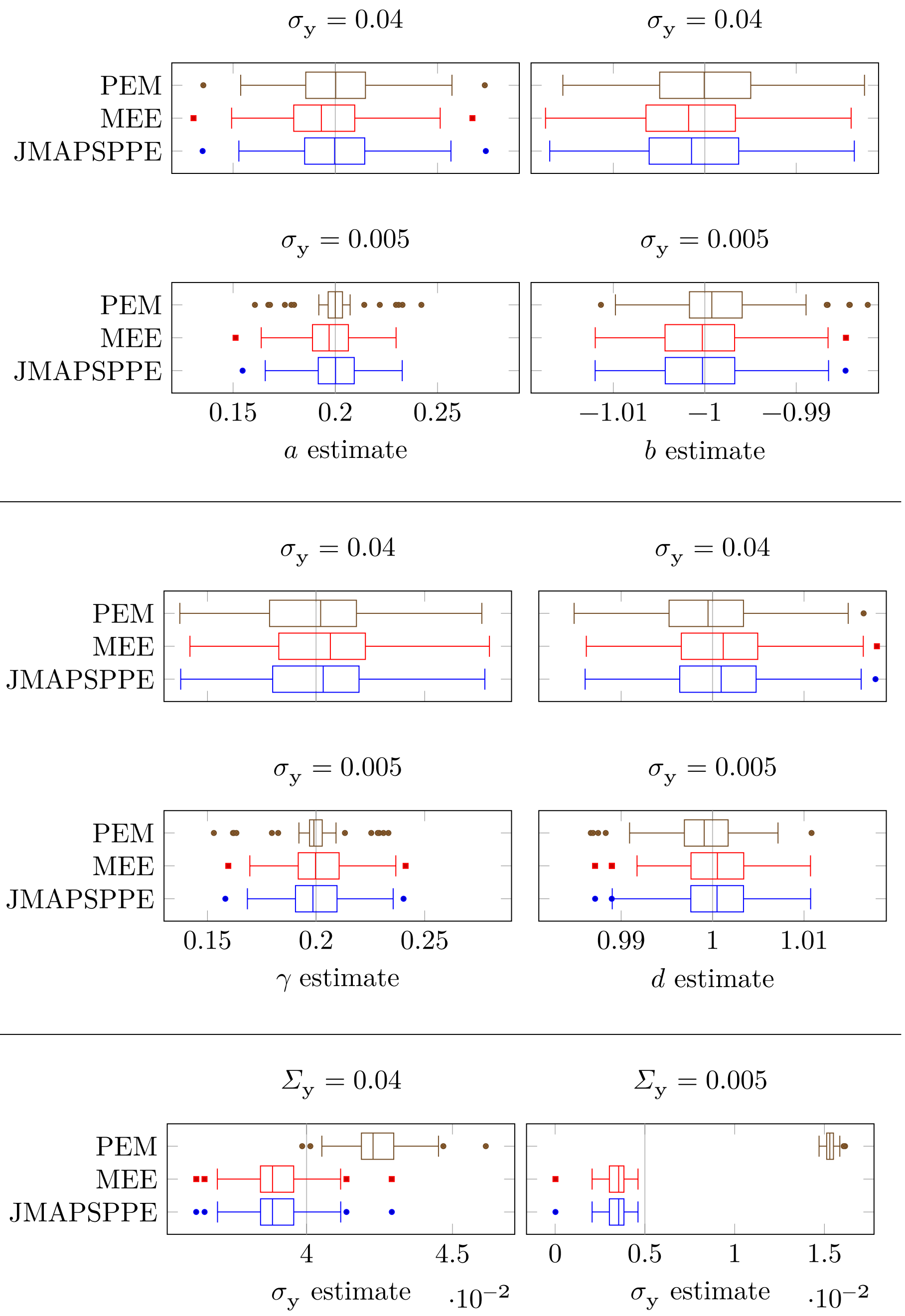

Figure 4.17: Boxplots of the Holmes–Rand oscillator parameter estimates. The nominal parameter values, $a = 0.2$, $b = -1$, $\gamma = 0.2$, and $d = 1$, are marked with gridlines in the plots.

## 4.2 Applications with experimental data

To illustrate a more practical application of the proposed estimators, we use them for flight-path reconstruction. In flight testing of aircraft, data from various sensors is collected on a wide range of physical quantities related to its movement. Common sensors are altimeters, airspeed sensors, accelerometers, global positioning system (GPS), gyroscopes, and magnetometers, to name a few. The physical quantities these sensors measure are all related through a well-known kinematic model. However, biases and noise in the sensor readings make the measured paths incoherent with the model. For example, the measured velocities might differ from the integral of the acceleration, the measured positions might differ from the integral of the velocity. Similar effects occur for the attitude, angular velocities, and angular accelerations. Flight-path reconstruction overcomes these limitations by using the kinematic model of the aircraft to obtain a coherent flight path and measurement model from the data (Mulder et al., 1999; Jategaonkar, 2006, Chap. 10; Klein and Morelli, 2006, Chap. 10). In the context of flight testing, flight path reconstruction is also known as data compatibility check.

For this example we used data collected from a VFW-Fokker 614 of the German Aerospace Center's (DLR) Advanced Technologies Testing Aircraft System (ATTAS) project, corresponding to a bank-to-bank roll manuever. This data accompanies the book "Flight Vehicle System Identification" (Jategaonkar, 2006) and is used with the author's permission. The measurements were collected at a rate of 25 Hz and correspond to the following quantities: airspeed $\tilde{v}_{\mathrm{m}}$, angle of attack at the noseboom $\alpha_{\mathrm{m}}$, angle of sideslip at the noseboom $\beta_{\mathrm{m}}$, roll angle $\phi_{\mathrm{m}}$, pitch angle $\theta_{\mathrm{m}}$, yaw angle $\psi_{\mathrm{m}}$, altitude $h_{\mathrm{m}}$, roll rate $p_{\mathrm{m}}$, pitch rate $q_{\mathrm{m}}$, yaw rate $r_{\mathrm{m}}$, longitudinal acceleration $a_{\mathrm{xm}}$, lateral acceleration $a_{\mathrm{ym}}$, and vertical acceleration $a_{\mathrm{zm}}$.

For the estimation, we considered the derivatives $D_{\mathrm{x}}$, $D_{\mathrm{y}}$ and $D_{\mathrm{z}}$ of the external accelerations the derivatives $D_{\mathrm{l}}$, $D_{\mathrm{m}}$, $D_{\mathrm{n}}$ and of the normalized moments of force to be linear diffusion processes evolving according to the following SDEs:

$$
\begin{aligned}
\mathrm{d}D_{\mathrm{x}}(t) &= -K_{\mathrm{a}}D_{\mathrm{x}}(t)\,\mathrm{d}t + \sigma_{\mathrm{a}}\,\mathrm{d}W_t^{(1)},\\
\mathrm{d}D_{\mathrm{y}}(t) &= -K_{\mathrm{a}}D_{\mathrm{y}}(t)\,\mathrm{d}t + \sigma_{\mathrm{a}}\,\mathrm{d}W_t^{(2)},\\
\mathrm{d}D_{\mathrm{z}}(t) &= -K_{\mathrm{a}}D_{\mathrm{z}}(t)\,\mathrm{d}t + \sigma_{\mathrm{a}}\,\mathrm{d}W_t^{(3)},\\
\mathrm{d}D_{\mathrm{l}}(t) &= -K_{\mathrm{m}}D_{\mathrm{l}}(t)\,\mathrm{d}t + \sigma_{\mathrm{m}}\,\mathrm{d}W_t^{(4)},\\
\mathrm{d}D_{\mathrm{m}}(t) &= -K_{\mathrm{m}}D_{\mathrm{m}}(t)\,\mathrm{d}t + \sigma_{\mathrm{m}}\,\mathrm{d}W_t^{(5)},\\
\mathrm{d}D_{\mathrm{n}}(t) &= -K_{\mathrm{m}}D_{\mathrm{n}}(t)\,\mathrm{d}t + \sigma_{\mathrm{m}}\,\mathrm{d}W_t^{(6)},
\end{aligned}
$$

where $K_{\mathrm{a}}$ and $K_{\mathrm{m}}$ are unknown parameters representing the damping coefficient of the diffusion, $\sigma_{\mathrm{a}}$ and $\sigma_{\mathrm{m}}$ are the diffusion coefficients of the accelerations and the moments, respectively, and $W^{(i)}$ are Wiener processes representing the process noise. The remaining states are not directly subject to noise. This representation is similar to the estimation before modeling (EBM) approach (Sri-Jayantha and Stengel, 1988; Jategaonkar, 2006, Sec. 10.4).

The gamma distribution was used as the prior of $K_\mathrm{a}$ and $K_\mathrm{m}$ parameters. For the remaining parameters and the initial conditions, non-informative Gaussian priors were used. For the measurement model, similarly, we considers that all sensors were subject to independent Gaussian noise, and that the accelerometers and gyrometers had, furthermore, unknown biases to be estimated. Putting the whole model in the format of Chapter 2, we have that it is characterized as follows. The clean state vector $z \in \mathbb{R}^{16}$ is composed of the following elements, in order: roll angle $\phi_\mathrm{a}$, pitch angle $\theta_\mathrm{a}$, yaw angle $\psi_\mathrm{a}$, roll rate $p$, pitch rate $q$, yaw rate $r$, normalized rolling moment $\ell$, normalized pitching moment $m$, normalized yawing moment $n$, altitude $h_\mathrm{a}$, longitudinal inertial velocity $u$, lateral inertial velocity $v$, vertical inertial velocity $w$, longitudinal inertial acceleration $a_\mathrm{x}$, lateral inertial acceleration $a_\mathrm{y}$, vertical inertial acceleration $a_\mathrm{z}$. The noisy state vector $x \in \mathbb{R}^6$, in turn, is composed of the following elements, in order: longitudinal, lateral and vertical jerk $d_\mathrm{x}$, $d_\mathrm{y}$, and $d_\mathrm{z}$, and rolling, pitching and yawing jerk $d_\mathrm{l}$, $d_\mathrm{m}$, and $d_\mathrm{n}$. The parameter vector $\theta \in \mathbb{R}^8$ is composed of the following elements, in order: the angular velocity measurement biases $b_\mathrm{p}$, $b_\mathrm{q}$, and $b_\mathrm{r}$, the acceleration measurement biases $b_\mathrm{ax}$, $b_\mathrm{ay}$, and $b_\mathrm{az}$, and the damping coefficients $k_\mathrm{a}$ and $k_\mathrm{m}$. The unnormalized pitching, rolling and yawing moments, respectively, can be obtained from the normalized moments as

$$\frac{J_\mathrm{x} J_\mathrm{z} - J_{\mathrm{x}z}^2}{J_\mathrm{z}}\ell, \qquad \frac{1}{J_\mathrm{y}}m, \qquad \frac{J_\mathrm{x} J_\mathrm{z} - J_{\mathrm{x}z}^2}{J_\mathrm{x}}n.$$

The noisy drift function is given by the following linear model:

$$f(t, x, z, \theta) := \begin{bmatrix} -k_\mathrm{a} d_\mathrm{x} \\ -k_\mathrm{a} d_\mathrm{y} \\ -k_\mathrm{a} d_\mathrm{z} \\ -k_\mathrm{m} d_\mathrm{l} \\ -k_\mathrm{m} d_\mathrm{m} \\ -k_\mathrm{m} d_\mathrm{n} \end{bmatrix}.$$

The clean drift, in turn, is given by the equations of rigid body kinematics

(Sri-Jayantha and Stengel, 1988; Jategaonkar, 2006, Eq. 10.26):

$$h(t,x,z,\theta) := \begin{bmatrix} p + q\sin(\phi_\mathrm{a})\tan(\theta_\mathrm{a}) + r\cos(\phi_\mathrm{a})\tan(\theta_\mathrm{a}) \\ q\cos(\phi_\mathrm{a}) - r\sin(\phi_\mathrm{a}) \\ q\sin(\phi_\mathrm{a})\sec(\theta_\mathrm{a}) + r\cos(\phi_\mathrm{a})\sec(\theta_\mathrm{a}) \\ pqC_{11} + qrC_{12} + qC_{13} + \ell + nC_{14} \\ prC_{21} + (r^2 - p^2)C_{22} - rC_{23} + m \\ pqC_{31} + qrC_{32} + qC_{33} + \ell C_{34} + n \\ d_\mathrm{l} \\ d_\mathrm{m} \\ d_\mathrm{n} \\ u\sin(\theta_\mathrm{a}) - v\cos(\theta_\mathrm{a})\sin(\phi_\mathrm{a}) - w\cos(\theta_\mathrm{a})\cos(\phi_\mathrm{a}) \\ -qw + rv - g\sin(\theta_\mathrm{a}) + a_\mathrm{x} \\ -ru + pw + g\cos(\theta_\mathrm{a})\sin(\phi_\mathrm{a}) + a_\mathrm{y} \\ qu - pv + g\cos(\theta_\mathrm{a})\cos(\phi_\mathrm{a}) + a_\mathrm{z} \\ d_\mathrm{x} \\ d_\mathrm{y} \\ d_\mathrm{z} \end{bmatrix}.$$

The prior log-density is given by

$$\begin{aligned} \ln \pi(x,z,\theta) :=& -\frac{\phi_\mathrm{a}^2 + \theta_\mathrm{a}^2 + \psi_\mathrm{a}^2 + p^2 + q^2 + r^2 + \ell^2 + m^2 + n^2}{2\sigma_0} \\ &- \frac{u^2 + v^2 + w^2 + a_\mathrm{z}^2 + a_\mathrm{y}^2 + a_\mathrm{z}^2 + d_\mathrm{x}^2 + d_\mathrm{y}^2 + d_\mathrm{z}^2 + d_\mathrm{l}^2 + d_\mathrm{m}^2 + d_\mathrm{n}^2}{2\sigma_0} \\ &- \frac{b_\mathrm{ax}^2 + b_\mathrm{ay}^2 + b_\mathrm{az}^2 + b_\mathrm{p}^2 + b_\mathrm{q}^2 + b_\mathrm{r}^2}{2\sigma_0} + 9\ln(k_\mathrm{a}k_\mathrm{m}) - 9(k_\mathrm{a} + k_\mathrm{m}), \end{aligned}$$

where a large value for the standard deviation $\sigma_0$ was chosen to make the prior non-informative. The measurement log-likelihood is given by

$$\ln \psi(y|x,z,\theta) := \sum_{k=0}^{N} \ln \psi_k \Big( y_k | x(kt_\mathrm{s}), z(kt_\mathrm{s}), \theta \Big),$$

where $\psi_k \colon \mathbb{R}^{13} \times \mathbb{R}^{16} \times \mathbb{R}^{6} \times \mathbb{R}^{12} \to \mathbb{R}_{\geq 0}$, the log-likelihood of each measurement, is given by

$$\begin{aligned} \ln \psi_k(y|x,z,\theta) :=& -\frac{(\phi_\mathrm{m} - \phi_\mathrm{a})^2}{2\sigma_\phi^2} - \frac{(\theta_\mathrm{m} - \theta_\mathrm{a})^2}{2\sigma_\theta^2} - \frac{(\psi_\mathrm{m} - \psi_\mathrm{a})^2}{2\sigma_\psi^2} \\ &- \frac{[\alpha_\mathrm{m} - \operatorname{atan2}(w - qx_\mathrm{nb} + py_\mathrm{nb}, u - ry_\mathrm{nb} + qz_\mathrm{nb})]^2}{2\sigma_\alpha^2} \\ &- \frac{[\beta_\mathrm{m} - \operatorname{atan2}(v - pz_\mathrm{nb} + rx_\mathrm{nb}, u - ry_\mathrm{nb} + qz_\mathrm{nb})]^2}{2\sigma_\alpha^2} \end{aligned}$$

$$- \frac{(h_\mathrm{m} - h_\mathrm{a})^2}{2\sigma_\mathrm{h}^2} - \frac{(p_\mathrm{m} - b_\mathrm{p} - p)^2}{2\sigma_\mathrm{p}^2} - \frac{(q_\mathrm{m} - b_\mathrm{q} - q)^2}{2\sigma_\mathrm{q}^2} - \frac{(r_\mathrm{m} - b_\mathrm{r} - r)^2}{2\sigma_\mathrm{r}^2}$$
$$- \frac{(b_\mathrm{ax} + a_\mathrm{x} + qz_\mathrm{as} - ry_\mathrm{as} - a_\mathrm{xm})^2}{2\sigma_\mathrm{ax}^2} - \frac{(b_\mathrm{ay} + a_\mathrm{y} + rx_\mathrm{as} - pz_\mathrm{as} - a_\mathrm{ym})^2}{2\sigma_\mathrm{ay}^2}$$
$$- \frac{(b_\mathrm{az} + a_\mathrm{z} + py_\mathrm{as} - qx_\mathrm{as} - a_\mathrm{zm})^2}{2\sigma_\mathrm{az}^2} - \frac{[\tilde{v}_\mathrm{m} - \sqrt{u^2 + v^2 + w^2}]^2}{2\sigma_\mathrm{v}^2},$$

where $x_\mathrm{nb}$, $y_\mathrm{nb}$ and $z_\mathrm{nb}$ are the coordinates, with respect to the center of gravity, of the tip of the nose boom where the aerodynamic probe is located, $x_\mathrm{as}$, $y_\mathrm{as}$ and $z_\mathrm{as}$ are the coordinates of the location of the accelerometer, $g$ is the acceleration of gravity, and $\sigma__$ are the measurement noise standard deviations. Constant terms which do not influence the location of maxima have been omitted from the expression above and the standard deviations were treated as known parameters, chosen empirically.

Both the JMAPSPPE and MEE estimates were compared with the output error method (OEM) described in (Jategaonkar, 2006, Sec. 10.5) and implemented in its accompanying materials. The outputs corresponding to the reconstructed paths can be seen in Figures 4.18 to 4.22

This example mostly demonstrates the feasibility of applying the proposed

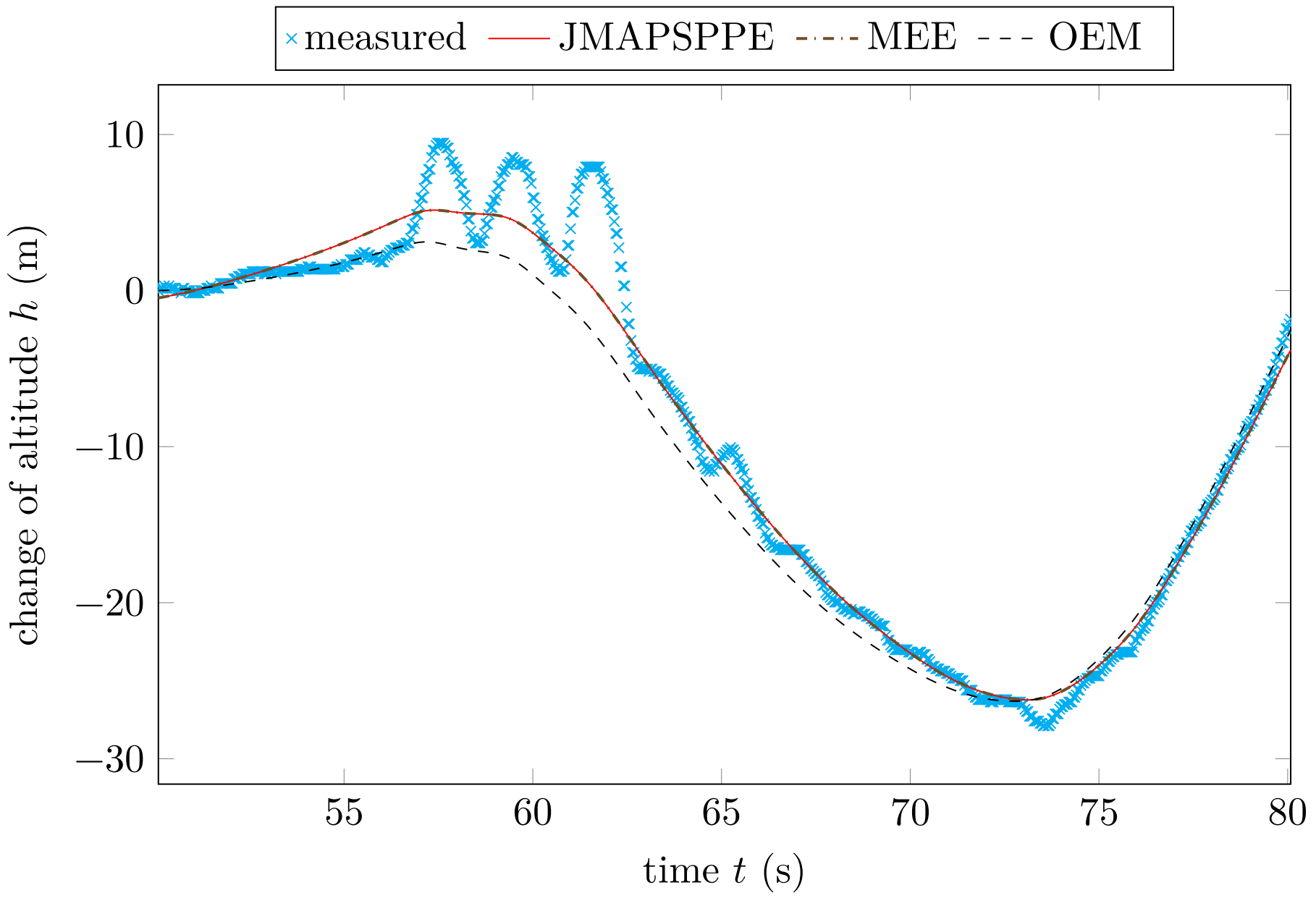

Figure 4.18: Altitude output corresponding to the reconstructed flight path using the various methods.

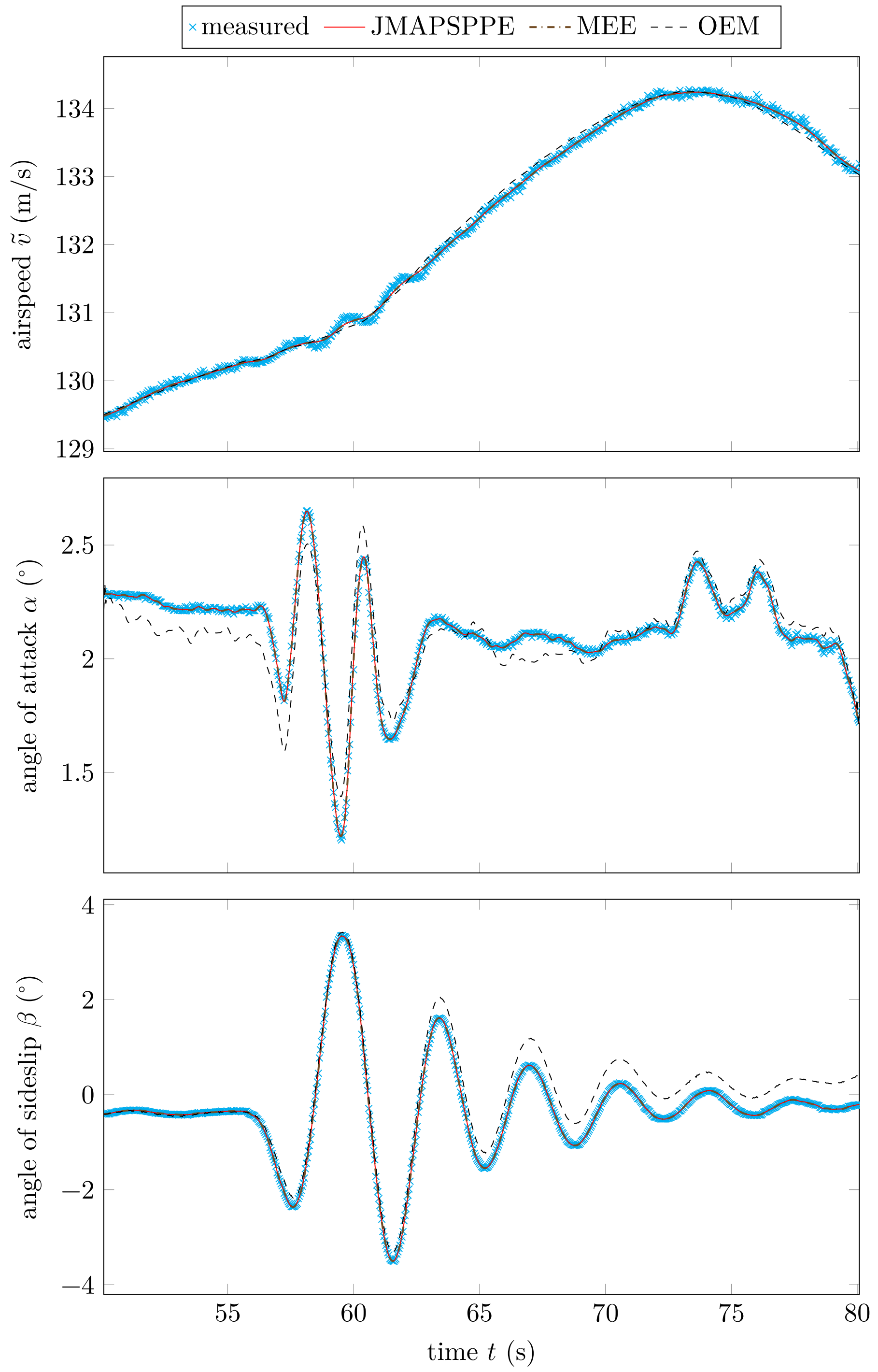

Figure 4.19: Velocity outputs corresponding to the reconstructed flight path using the various methods.

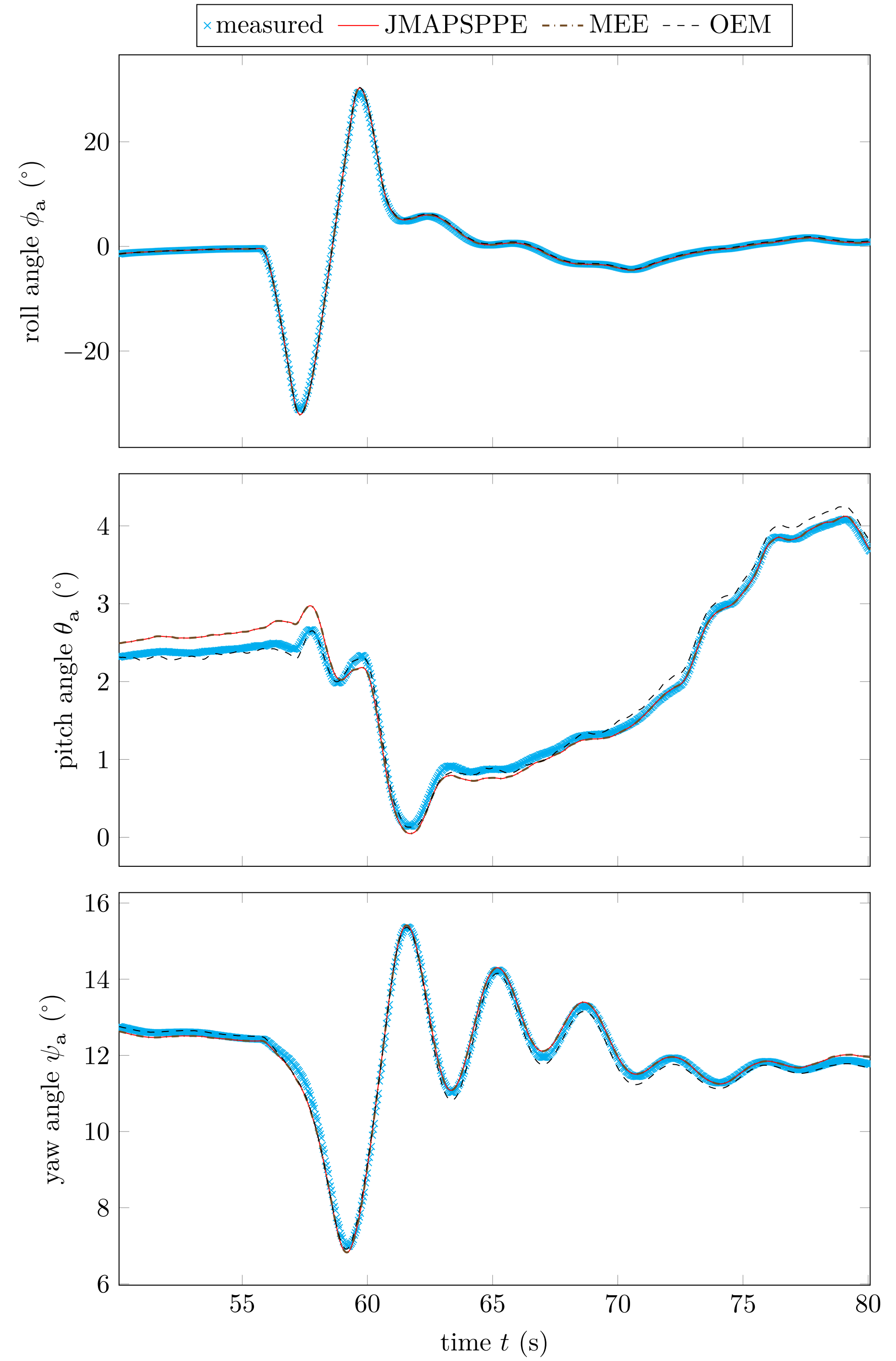

Figure 4.20: Attitude outputs corresponding to the reconstructed flight path using the various methods.

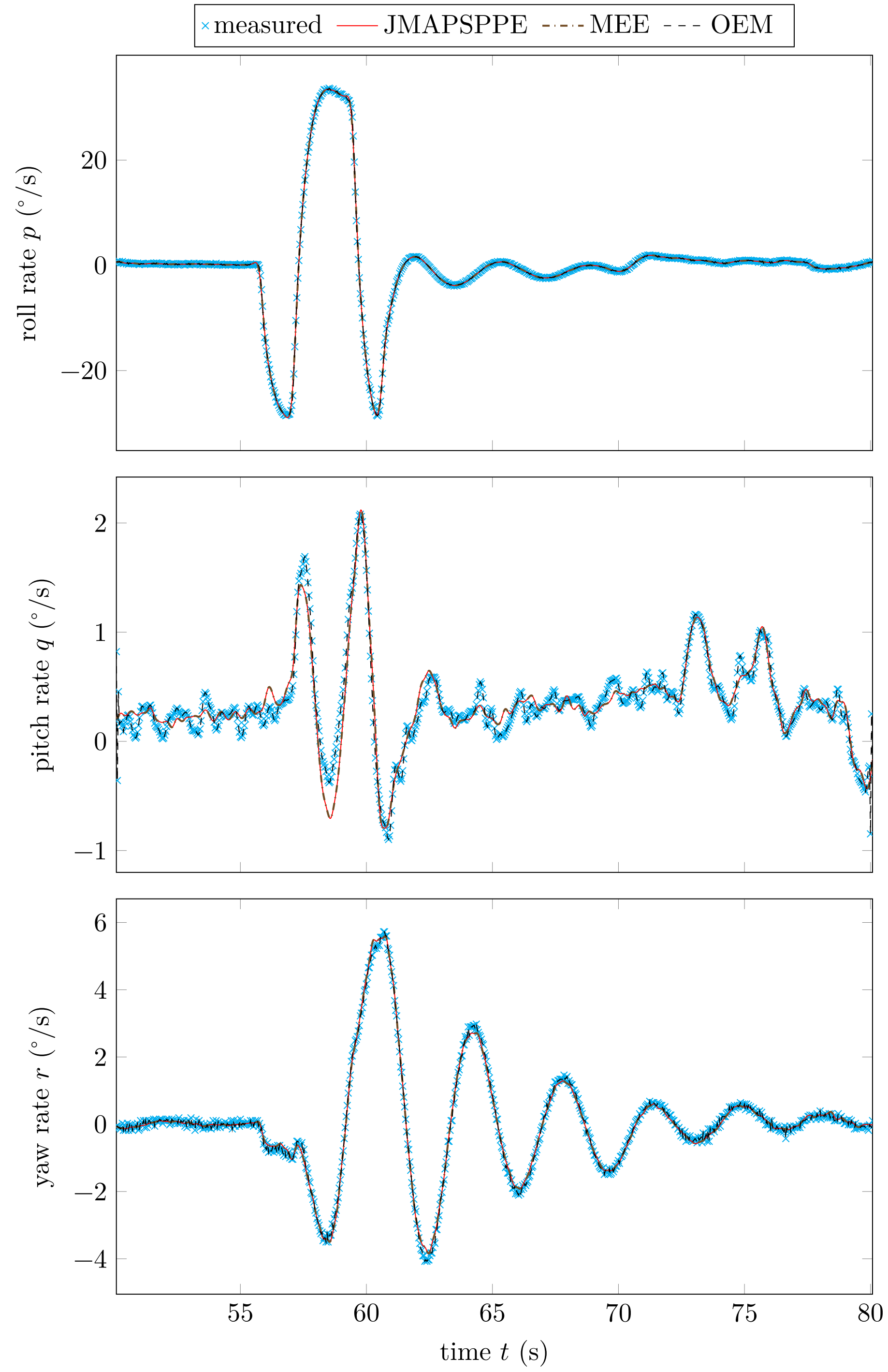

Figure 4.21: Angular velocity outputs corresponding to the reconstructed flight path using the various methods.

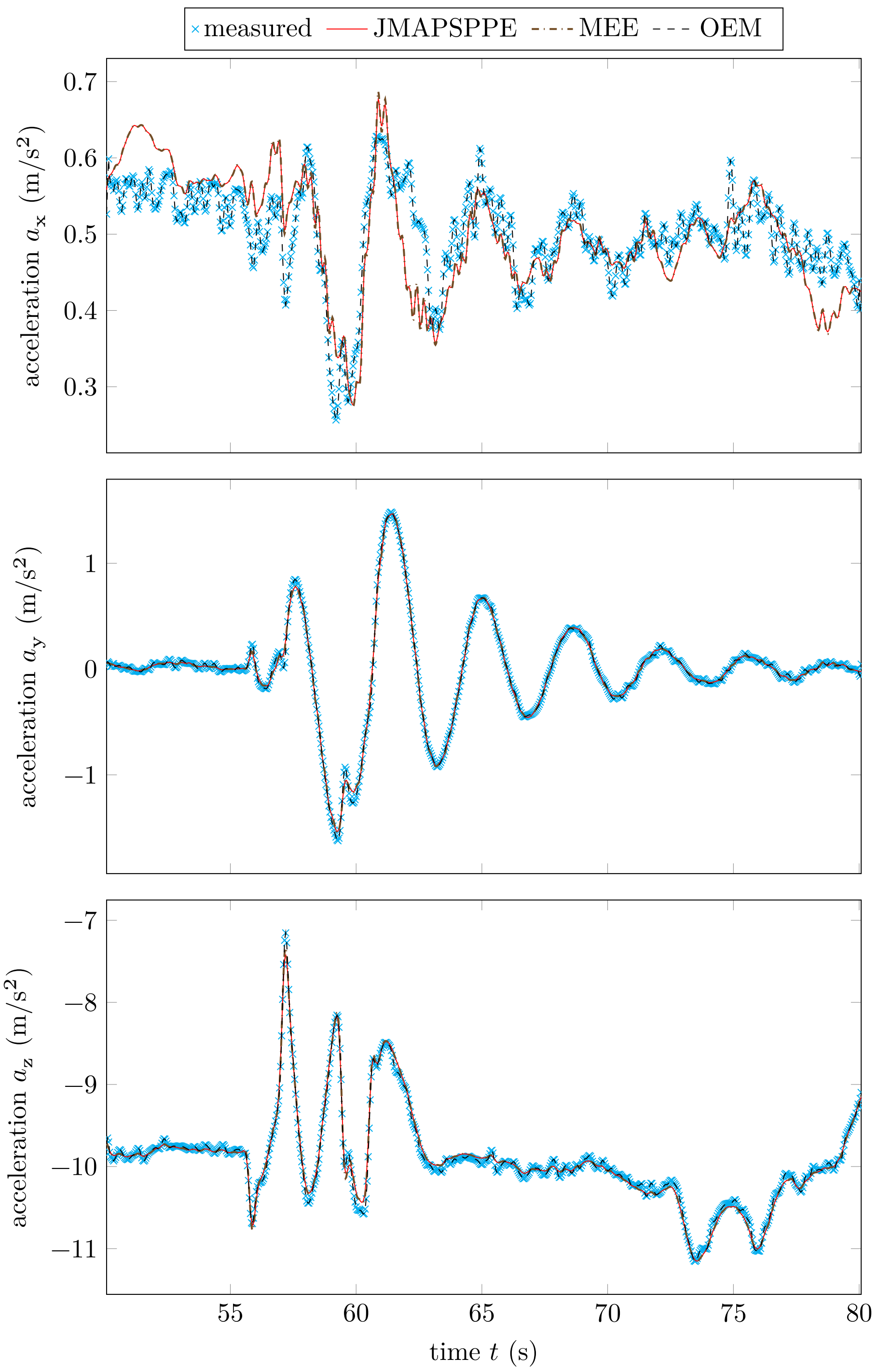

Figure 4.22: Acceleration outputs corresponding to the reconstructed flight path using the various methods.

JMAPSPPE and MEE in a practical problem of technological interest with a high dimensionality. We note that the use of flight-path reconstruction in industry and production requires some further tuning of the algorithm, dynamical model structure, and measurement model which is outside the scope of this thesis. Nevertheless, the results of this example suggest that the MAP estimator can be a good substitute for the output error method in flight path reconstruction of small unmanned aerial vehicles for which the small payload limits the quality of the inertial sensors. The OEM relies on high-quality inertial sensors and provides inadequate results when this is not the case (Mulder et al., 1999).

# Chapter 5

# Conclusions

> **Cueball:** I used to think correlation implied causation. Then I took a statistics class. Now I don't.
> **Megan:** Sounds like the statistics class helped.
> **Cueball:** Well, maybe.
>
> Randall Munroe, xkcd #552: *Correlation*

In this chapter we recall the conclusions and contributions of this thesis, and finalize with directions for future work.

## 5.1 Conclusions

The main contribution of this thesis is building a solid theoretical foundation for maximum *a posteriori* (MAP) estimation in stochastic differential equations (SDEs). The first one is laying out a rigorous definition of mode and MAP estimation for general random variables in possibly infinite-dimensional spaces, presented in Section 2.1. This definition is shown to be coincident with the traditional one for continous and discrete random variables, and its interpretations in the context of Bayesian decision theory and Bayesian estimation are provided.

Next, we showed how this definition of mode can be used to obtain the posterior and prior joint fictious state-path and parameter densities for systems described by SDEs where not necessarily all states are under direct influence of noise. In addition, we showed that the minimum energy estimates correspond to the state-paths associated with the joint MAP noise-paths and parameters.

We then related the popular approach of using the discretized SDE for joint MAP state-path and parameter estimation to the continuous-time MAP estimators by using variational analysis. We proved that the Euler-discretized MAP estimator converges hypographically to the minimum energy estimator. The trapezoidally-discretized MAP estimator, on the other hand, converges hypographically to the continuous-time MAP estimator. This implies that the discretized estimates might have different interpretations depending on the discretization method used.

Some example applications with both simulated and experimental were then presented. In the simulated examples, we saw that the state-path estimates of both methods were comparable. However, the minimum-energy estimates of parameters which appear in the drift divergence were biased. This is because the Onsager–Machlup functional penalizes parameter values which increase the amplification of the noise by the system, while the energy functional does not. In addition, the use of the proposed estimators was demonstrated in non-Gaussian applications in which nonlinear Kalman smoothers and prediction error methods yield poor results or are not applicable. The example with experimental data showed the viability of the proposed methods in an application of technological interest with high model order and system complexity.

The theoretical foundations built with these contributions provide a firm basis on which new work can be built upon.

## 5.2 Future work

The work herein can be expanded in a number of directions. Important questions that need to be answered are conditions for existence of the MAP estimates. Preliminary analysis indicates that when the negative of the log-prior is coercive and the measurement likelihood and drift divergence are bounded in the parameter support, then at least one MAP estimate exists.

Under some regularity conditions, the posterior fictious density of state-paths and parameters admit a second order Fréchet derivative. We believe that this derivative is analogous to the Hessian matrix of posterior probability density function in $\mathbb{R}^d$, itself the Bayesian equivalent of the Fisher information matrix. The second Fréchet derivative of the fictitious density seems, furthermore, related to reproducing kernels associated with Gaussian processes (Parzen, 1963). As such it could be used to build Gaussian approximations to the posterior state-path and parameter distribution, and to provide a measure of uncertainty associated with the modal estimates.

Preliminary analysis shows that the inverse of the optimization problem's Hessian matrix is very close to the covariance matrices obtained with the unscented Kalman smoother. These ideas underlie the approach of (Karimi and McAuley, 2013, 2014a; Varziri et al., 2008a) to marginalization of the states. Care should be taken to understand and develop the theoretical foundations of this approach, however. In particular, to use the Hessian matrix of the optimization problem its relationship to the original problem's second order derivative should be understood using tools such as second-order variational analysis and hypo-convergence (cf. Rockafellar and Wets, 1998, Chap. 13).

The hypo-convergence of the Euler and trapezoidal discretizations can be strengthened by proving $\tau/\sigma$-equi-semicontinuity (Attouch, 1984, Sec. 2.6.2). This would imply hypo-convergence in the weak topology, which can guarantee the existence of a convergent sequence discretized MAP estimates. In addition, the techniques of variational analysis used in this thesis can be used to prove the hypo-convergence of the Legendre–Gauss–Lobatto direct transcription methods used in Chapter 4 to discretize the optimization problems. A particular strong kind of hypo-convergence seems to hold, in which not only the decision variables but also the adoint variables (co-states) converge hypographically. We also believe that $\tau/\sigma$-equi-semicontinuity holds under some regularity conditions.

With respect to applications, we intend to use the MAP estimator for flight path reconstruction and system identification of the small-scale hand-launched unmanned aerial vehicles being developed at the Universidade Federal de Minas Gerais. Additionally, we believe that the MAP state-path estimates can be used together with particle smoothers (Klaas et al., 2006) to help design a proposal distribution that samples in regions of high posterior probability. similar to the unscented particle filter (van der Merwe et al., 2000).

# Appendix A

# Collected theorems and definitions

In this chapter, we collect some useful theorems and definitions which are used throughout the thesis.

## A.1 Simple identities

**Lemma A.1.** *For all* $a, b \in \mathbb{R}_{\geq 0}$*,*

$$\sqrt{a} + \sqrt{b} \leq \sqrt{2}\sqrt{a+b}. \tag{A.1}$$

*Proof.* The square root is a concave function, so from Jensen's inequality we have that

$$\frac{\sqrt{a} + \sqrt{b}}{2} \leq \sqrt{\frac{a+b}{2}}.$$

Rearranging, we obtain (A.1). □

**Lemma A.2.** *For all* $a, b \in \mathbb{R}$*,*

$$(a+b)^2 \leq 2(a^2 + b^2).$$

*Proof.* Since $(a-b)^2 \geq 0$,

$$(a+b)^2 \leq (a+b)^2 + (a-b)^2 = 2a^2 + 2b^2.$$ □

**Lemma A.3.** *For all* $x \in \mathbb{R}_{\geq 0}$*,* $\exp(x) \geq 1 + x$*.*

*Proof.* From the Taylor series expansion of the exponential we have that

$$\exp(x) - (1+x) = \sum_{k=2}^{\infty} \frac{x^k}{k!} \geq 0$$ □

## A.2 Linear algebra

**Lemma A.4** (matrix Mercator series, Higham and Al-Mohy, 2010, Sec. 5.2)**.** *For all matrices $\boldsymbol{A} \in \mathbb{R}^{d\times d}$ with spectral radius smaller than unity, we have that the matrix logarithm satisfies*

$$\ln(\boldsymbol{I}+\boldsymbol{A}) = \sum_{k=1}^{\infty}(-1)^{k+1}\frac{\boldsymbol{A}^k}{k}.$$

**Lemma A.5** (matrix exponential determinant, Bernstein, 2009, Cor. 11.2.4)**.** *For all matrices $\boldsymbol{A} \in \mathbb{R}^{d\times d}$, the determinant of a matrix exponential satisfies*

$$\det\exp(\boldsymbol{A}) = \exp(\operatorname{tr}\boldsymbol{A})\,.$$

## A.3 Analysis

**Lemma A.6** (Dominance of $\|\cdot\|_{L^1_d}$ by $\|\cdot\|_{L^2_d}$)**.** *Let $(\mathcal{X},\mathcal{F},\mu)$ be a finite measure space. Then, for all $f \in L^2_d(\mathcal{X},\mathcal{F},\mu)$,*

$$\|f\|_{L^1_d} \le \sqrt{\mu(\mathcal{X})}\,\|f\|_{L^2_d}$$

*Proof.* From the Cauchy–Schwarz inequality we have that

$$\|f\|_{L^1_d} = \int_{\mathcal{X}} |f(x)|\,\mathrm{d}\mu(x) \le \left(\int_{\mathcal{X}} |1|^2\,\mathrm{d}\mu(x)\right)^{1/2}\left(\int_{\mathcal{X}} |f(x)|^2\,\mathrm{d}\mu(x)\right)^{1/2} \qquad \square$$

**Lemma A.7** (Dominance of $|||\cdot|||$ by $\|\cdot\|_{\mathcal{W}^2_d}$)**.** *Let $\mathcal{T} := [0,t_\mathrm{f}]$. Then there some $c \in \mathbb{R}_{>0}$ such that, for all $x \in \mathcal{W}^2_d(\mathcal{T})$,*

$$|||x||| \le c\,\|x\|_{\mathcal{W}^2_d}\,, \tag{A.2}$$

*where $|||x||| := \sup_{t\in\mathcal{T}}|x(t)|$ and $\|x\|_{\mathcal{W}^2_d} := \left(|x(0)|^2 + \int_0^{t_\mathrm{f}}|\dot{x}(t)|^2\,\mathrm{d}t\right)^{1/2}$.*

*Proof.* Every $x \in \mathcal{W}^2_d$ is absolutely continuous, meaning that for all $t \in \mathcal{T}$

$$x(t) = x(0) + \int_0^t \dot{x}(t)\,\mathrm{d}t,$$

implying that

$$|||x||| = \sup_{t\in\mathcal{T}}\left|x(0) + \int_0^t \dot{x}(t)\,\mathrm{d}t\right| \le |x(0)| + \int_0^{t_\mathrm{f}} |\dot{x}(t)|\,\mathrm{d}t.$$

Next, using Lemma A.6 we have that

$$|||x||| \leq |x(0)| + \sqrt{t_\mathrm{f}} \left( \int_0^{t_\mathrm{f}} |\dot{x}(t)|^2 \, \mathrm{d}t \right)^{1/2}.$$

Lemma A.1, in turn, gives us

$$|||x||| \leq \sqrt{2} \left( |x(0)|^2 + t_\mathrm{f} \int_0^{t_\mathrm{f}} |\dot{x}(t)|^2 \, \mathrm{d}t \right)^{1/2}.$$

Consequently, (A.2) is satisfied with $c := \max\left(\sqrt{2t_\mathrm{f}}, \sqrt{2}\right)$. □

**Lemma A.8** (Grönwall–Bellman inequality, Lem. 5.6.4 of Polak, 1997)**.** *Let $c, K, t_\mathrm{f} \in \mathbb{R}_{\geq 0}$ and $\mathcal{T} := [0, t_\mathrm{f}]$. Then, if an integrable function $h\colon \mathcal{T} \to \mathbb{R}$ satisfies*

$$h(t) \leq c + K \int_0^t h(\tau) \, \mathrm{d}\tau \qquad \text{for all } t \in \mathcal{T},$$

*we have that*

$$h(t) \leq c \exp(K t_\mathrm{f}) \qquad \text{for all } t \in \mathcal{T}.$$

**Lemma A.9** (adapted from Lemma 5.6.3 of Polak, 1997)**.** *Let $t_\mathrm{f} \in \mathbb{R}_{>0}$, $\mathcal{T} := [0, t_\mathrm{f}]$, and $w \in \mathcal{C}(\mathcal{T}, \mathbb{R}^d)$. In addition, let the continuous function $f\colon \mathcal{T} \times \mathbb{R}^d \to \mathbb{R}^d$ be Lipschitz continuous with respect to its second argument, uniformly with respect to its first, i.e., there exists $L \in \mathbb{R}_{>0}$ such that for all $t \in \mathcal{T}$ and $x', x'' \in \mathbb{R}^d$*

$$\left| f(t, x') - f(t, x'') \right| \leq L \left| x' - x'' \right|. \tag{A.3}$$

*We then have that, for all $x_0 \in \mathbb{R}^d$ there is a unique solution $x\colon \mathcal{T} \to \mathbb{R}^d$ to the integral equation*

$$x(t) = x_0 + \int_0^t f(t, x(\tau)) \, \mathrm{d}\tau + w(t), \tag{A.4}$$

*which satisfies, for all $\tilde{x} \in \mathcal{W}_d^2(\mathcal{T})$, the inequality*

$$|||\tilde{x} - x||| \leq \exp(L t_\mathrm{f}) \left( |x_0 - \tilde{x}(0)| + |||w||| + \int_0^{t_\mathrm{f}} \left| \dot{\tilde{x}}(t) - f(t, \tilde{x}(t)) \right| \mathrm{d}t \right). \tag{A.5}$$

*Proof.* First, note that the triangle inequality implies that

$$|f(t, x)| \leq |f(t, x) - f(t, 0)| + |f(t, 0)| \qquad \text{for all } t \in \mathcal{T}, x \in \mathbb{R}^d.$$

The Lipschitz condition (A.3) and the extreme value theorem, in turn, imply that there exists $c \in \mathbb{R}_{>0}$ such that

$$|f(t, x)| \leq L |x| + c \qquad \text{for all } t \in \mathcal{T}, x \in \mathbb{R}^d. \tag{A.6}$$

Next, let $\{\tilde{x}_i\}_{i=0}^{\infty}$ be a sequence of continuous functions $\tilde{x}_i \in \mathcal{C}(\mathcal{T}, \mathbb{R}^d)$ with its first element $\tilde{x}_0 := \tilde{x}$. The following elements of the sequence are given by the recursion

$$\tilde{x}_{i+1}(t) := x_0 + \int_0^t f(\tau, \tilde{x}_i(\tau))\,\mathrm{d}\tau + w(t), \qquad \text{for all } t \in \mathcal{T}. \tag{A.7}$$

The growth bound (A.6) implies that $t \mapsto f(t, \tilde{x}_i(t))$ is integrable if $\tilde{x}_i$ is continuous. Furthermore, (A.7) and the continuity of $w$ imply that $\tilde{x}_{i+1}$ is continuous if $\tilde{x}_i$ is continuous. As $\tilde{x}_0$ is integrable, by recursion we then have that all $x_i$ are well-defined and continuous.

Since $\tilde{x} \in \mathcal{W}_d^2(\mathcal{T})$, we have that it is weakly differentiable and satisfies

$$\tilde{x}(t) = x(0) + \int_0^t \dot{\tilde{x}}(\tau)\,\mathrm{d}\tau, \qquad \text{for all } t \in \mathcal{T}.$$

Together with (A.7) and the fact that $\tilde{x}_0 := \tilde{x}$, this implies that

$$|\tilde{x}_1(t) - x_0(t)| \leq |x_0 - \tilde{x}(0)| + |||w||| + \int_0^{t_\mathrm{f}} \left|\dot{\tilde{x}}(t) - f(t, \tilde{x}(t))\right| \mathrm{d}t =: \epsilon.$$

For any $i > 0$ and $t \in \mathcal{T}$, from the recursion (A.7) and the Lipschitz condition (A.3) we have that

$$\begin{aligned} |\tilde{x}_{i+1}(t) - \tilde{x}_i(t)| &\leq \int_0^t |f(\tau, \tilde{x}_i(\tau)) - f(\tau, \tilde{x}_{i-1}(\tau))|\,\mathrm{d}\tau \\ &\leq \int_0^t L\,|\tilde{x}_i(\tau) - \tilde{x}_{i-1}(\tau)|\,\mathrm{d}\tau. \end{aligned}$$

By induction we then have that for all $i \in \{0, 1, \dots\}$

$$|\tilde{x}_{i+1}(t) - \tilde{x}_i(t)| \leq \frac{(Lt)^i}{i!}\epsilon,$$

which in turn implies that

$$|||\tilde{x}_{i+1} - \tilde{x}_i||| \leq \frac{(Lt_\mathrm{f})^i}{i!}\epsilon, \qquad |||\tilde{x}_\ell - \tilde{x}_k||| \leq \epsilon \sum_{i=k}^{\ell-1} \frac{(Lt_\mathrm{f})^i}{i!}. \tag{A.8}$$

Noting the similarity to the Taylor series expansion of the exponential, we can then see that $\{\tilde{x}_i\}_{i=1}^{\infty}$ is a Cauchy sequence, which converges due to completeness of $\mathcal{C}(\mathcal{T}, \mathbb{R}^d)$. Denoting its limit by $x$, we then have by (A.8) that for all $i \in \{0, 1, \dots\}$,

$$|||x - \tilde{x}_i||| \leq \epsilon \sum_{i=0}^{\infty} \frac{(Lt_\mathrm{f})^i}{i!} = \exp(Lt_\mathrm{f})\epsilon,$$

implying that (A.5) holds.

Finally, to see that (A.4) holds for the limit of the sequence, note that

$$\lim_{i\to\infty}\left|\int_0^t [f(\tau,\tilde{x}_i(\tau)) - f(\tau,x(\tau))]\,\mathrm{d}\tau\right| \le \lim_{i\to\infty} Lt_{\mathrm{f}}\,\|\!|\tilde{x}_i - x|\!\|=0.$$

Consequently, by letting $i\to\infty$ in (A.7) we obtain (A.4). □

**Corollary A.10.** *Let $\mathcal{T}$ and $f$ be as in Lemma A.9. Then for all $x_0\in\mathbb{R}^d$ there exists a unique solution $x\in\mathcal{W}_d^2([0,\tau])$ to the initial value problem*

$$\dot{x}(t) = f(t,x(t))\,,\qquad x(0)=x_0. \tag{A.9}$$

*Furthermore, for all $\tilde{x}\in\mathcal{W}_d^2([0,\tau])$,*

$$\|\!|\tilde{x}-x|\!\| \le \exp(L\tau)\left(|x(0)-\tilde{x}(0)| + \int_0^\tau \left|\dot{\tilde{x}}(t) - f(t,\tilde{x}(t))\right|\mathrm{d}t\right).$$

*Proof.* This corollary follows directly from Lemma A.9 by letting $w=0$. If $x\mathcal{C}(\mathcal{T},\mathbb{R}^n)$ is the solution to the integral equation

$$x(t) = x_0 + \int_0^t f(t,x(\tau))\,\mathrm{d}\tau + w(t),$$

then by Lebesgue's differentiation theorem we have that it is also a solution to the initial value problem (A.9). □

**Remark A.11** (regularization to ensure Lipschitz continuity)**.** *Let $r\colon\mathbb{R}\to\mathbb{R}$ and $f\colon\mathbb{R}^m\to\mathbb{R}^n$ be continuously differentiable and let $r$ satisfy, for some $M\in\mathbb{R}_{>0}$,*

$$r(\epsilon)=\begin{cases}1 & \text{if } \epsilon\le M\\ 0 & \text{if } \epsilon\ge 2M.\end{cases}$$

*Then the function $\tilde{f}(x) := f(x)r(|x|)$ is bounded, Lipschitz continuous, and satisfies $\tilde{f}(x)=f(x)$ for all $|x|\le M$.*

*Proof.* From the definition of $r$ and $\tilde{f}$, we have that

$$\sup_{x\in\mathbb{R}^m}\left|\tilde{f}(x)\right| = \sup_{|x|\le 2M}\left|\tilde{f}(x)\right|$$
$$\sup_{x\in\mathbb{R}^m}\left|\nabla\tilde{f}(x)\right| = \sup_{|x|\le 2M}\left|\nabla\tilde{f}(x)\right|.$$

From the extreme value theorem we then have that both $\tilde{f}$ and its derivatives are bounded. To conclude, we have that differentiable functions with bounded derivatives are Lipschitz continuous. □

## A.4 Probability theory and stochastic processes

The following two lemmas are adapted from Ikeda and Watanabe (1981, p. 449) and are frequently used to obtain the Onsager–Machlup functional.

**Lemma A.12** (limit of conditional exponential moments). *Let $\{\mathbb{A}_\epsilon\}_{\epsilon\in\mathbb{R}_{>0}}$ be a family of non-null events in $\mathcal{E}$. Then, if $A$ is a $\mathbb{R}$-valued random variable such that*

$$\limsup_{\epsilon\downarrow 0} \mathrm{E}[\exp(cA) \,|\, \mathbb{A}_\epsilon] \le 1 \qquad \text{for all } c \in \mathbb{R} \tag{A.10}$$

*then*

$$\lim_{\epsilon\downarrow 0} \mathrm{E}[\exp(cA) \,|\, \mathbb{A}_\epsilon] = 1 \qquad \text{for all } c \in \mathbb{R}. \tag{A.11}$$

*Proof.* Applying the Cauchy–Schwarz inequality we have that, for all $c \in \mathbb{R}$,

$$1 = \left(\mathrm{E}\left[\exp\left(\tfrac{cA}{2}\right)\exp\left(\tfrac{-cA}{2}\right) \,\middle|\, \mathbb{A}_\epsilon\right]\right)^2 \le \mathrm{E}[\exp(cA) \,|\, \mathbb{A}_\epsilon]\,\mathrm{E}[\exp(-cA) \,|\, \mathbb{A}_\epsilon]\,.$$

This in turn implies that

$$\mathrm{E}[\exp(cA) \,|\, \mathbb{A}_\epsilon] \ge \frac{1}{\mathrm{E}[\exp(-cA) \,|\, \mathbb{A}_\epsilon]}.$$

Taking the limit inferior and applying (A.10) we obtain the following bound:

$$\liminf_{\epsilon\downarrow 0} \mathrm{E}[\exp(cA) \,|\, \mathbb{A}_\epsilon] \ge \frac{1}{\limsup_{\epsilon\downarrow 0} \mathrm{E}[\exp(-cA) \,|\, \mathbb{A}_\epsilon]} \ge 1. \tag{A.12}$$

As the limit superior is always greater than or equal to the limit inferior, we conclude that the limits of (A.10) and (A.12) coincide and equal to one. Consequently, (A.11) holds. □

**Lemma A.13** (conditional exponential moments of a sum). *Let $\{\mathbb{A}_\epsilon\}_{\epsilon\in\mathbb{R}_{>0}}$ be a family of non-null events in $\mathcal{E}$. Then, if $A_1, \dots, A_n$ are $\mathbb{R}$-valued random variables such that*

$$\limsup_{\epsilon\downarrow 0} \mathrm{E}[\exp(cA_i) \,|\, \mathbb{A}_\epsilon] \le 1 \qquad \text{for all } c \in \mathbb{R} \text{ and } i = 1, \dots, n, \tag{A.13}$$

*then*

$$\lim_{\epsilon\downarrow 0} \mathrm{E}[\exp(cA_1 + \dots + cA_n) \,|\, \mathbb{A}_\epsilon] = 1 \qquad \text{for all } c \in \mathbb{R}. \tag{A.14}$$

*Proof.* Define the $\mathbb{R}$-valued random variables $B_i$, $i = 1, \dots, n$ as

$$B_i := \sum_{j=1}^{i} A_i.$$

We then have that, for $i = 1$,

$$\limsup_{\epsilon \downarrow 0} \mathrm{E}[\exp(cB_i) \,|\, \mathbb{A}_\epsilon] \leq 1 \qquad \text{for all } c \in \mathbb{R}. \tag{A.15}$$

Next, assume that (A.15) holds for some $i < n$. Then, applying the Cauchy–Scharwz inequality,

$$\mathrm{E}[\exp(cB_i)\exp(cA_{i+1}) \,|\, \mathbb{A}_\epsilon] \leq \sqrt{\mathrm{E}[\exp(2cB_i) \,|\, \mathbb{A}_\epsilon]\, \mathrm{E}[\exp(2cA_{i+1}) \,|\, \mathbb{A}_\epsilon]}.$$

Taking the limit superior and applying (A.13) and (A.15) we then have that (A.15) holds for $i + 1$ and, by induction, up to $i = n$. Applying Lemma A.12 we then have that (A.14) holds. □

**Lemma A.14** (Dembo and Zeitouni, 1998, Lemma 5.2.1)**.** *For all $\tau, \delta \in \mathbb{R}_{>0}$, if $W$ is an $n$-dimensional Wiener process then*

$$P\Big(\sup_{t \in [0,\tau]} \|W_t\| \geq \delta\Big) \leq 4n \exp\left(-\frac{\delta^2}{2n\tau}\right).$$

The following three lemmas are used to calculate what is referred to in the literature as the first exit or sojourn probability of the Wiener process (cf. Fujita and Kotani, 1982; Gasanenko, 1999). It corresponds to the probability that the Wiener process starting at $w$ sojourns in a set $\mathcal{W}$ for longer than $t$. Theorem A.15 is one of the many connections between partial differential equations and diffusion processes.

**Theorem A.15** (sojourn probability of the Wiener process, Patie and Winter, 2008, Thm. 7)**.** *Let $W$ be an $n$-dimensinal Wiener process over $\mathbb{R}_{\geq 0}$, $\mathcal{W} \subset \mathbb{R}^n$ be an open, bounded and connected domain with a smooth boundary $\partial\mathcal{W}$, the $\overline{\mathbb{R}}$-valued random variable $T_w$ be the time of first exit of the Wiener process starting at $w$ from $\mathcal{W}$ and $u\colon \mathbb{R}_{\geq 0} \times \bar{\mathcal{W}} \to [0, 1]$ be the probability that $T_w$ occurs after a given time interval, i.e.,*

$$T_w(\omega) := \inf\left\{t \in \mathbb{R}_{\geq 0} \,\middle|\, w + W_t \notin \mathcal{W}\right\}, \qquad u(w, t) := P(T_w > t).$$

*Then $u$ satisfies the boundary value problem*

$$\frac{\partial u}{\partial t}(t, w) = \frac{1}{2} \operatorname{tr} \nabla^2_{\mathrm{w}} u(t, w) \qquad \text{for all } w \in \mathcal{W} \text{ and } t \in \mathbb{R}_{\geq 0} \tag{A.16a}$$

$$u(t, w) = 0 \qquad \text{for all } w \in \partial\mathcal{W} \text{ and } t \in \mathbb{R}_{>0} \tag{A.16b}$$

$$u(0, w) = 1 \qquad \text{for all } w \in \mathcal{W}, \tag{A.16c}$$

*where $\operatorname{tr} \nabla^2_{\mathrm{w}} u$ denotes the Laplacian operator applied to $u$, i.e., the trace of its Hessian matrix with respect to its second argument $w$.*

The boundary value problem (A.16) in Theorem A.15 happens to have an analytical solution, which is useful in calculating assymptotic sojourn probabilities for more general diffusion processes, as done in Sec. 2.2.1. To obtain these solutions, we first state the following lemma.

**Lemma A.16** (Dirichlet eigenfunction basis). *Let $\mathcal{W} \subset \mathbb{R}^n$ be an open, bounded and connected domain with a smooth boundary $\partial\mathcal{W}$. Then there exists a sequence of not identically zero $\mathcal{C}(\bar{\mathcal{W}}, \mathbb{R})$ functions $\{\phi_i\}_{i=1}^{\infty}$ and a corresponding nondecreasing sequence of $\mathbb{R}_{>0}$ numbers $\{\lambda_i\}_{i=1}^{\infty}$ that solve the Dirichlet eigenvalue problem*

$$\lambda_i \phi_i(w) = -\frac{1}{2} \operatorname{tr} \nabla^2 \phi_i(w) \qquad \text{for all } w \in \mathcal{W} \tag{A.17a}$$

$$\phi_i(w) = 0 \qquad \text{for all } w \in \partial\mathcal{W}. \tag{A.17b}$$

*Furthermore, $\{\phi_i\}_{i=1}^{\infty}$ is an orthonormal basis for $L^2(\bar{\mathcal{W}})$ and $0 < \lambda_1 < \lambda_2$.*

For a proof of Lemma A.16 and some other interesting properties of the eigenvalue problem, refer to Larsson and Thomée (2003, Chap. 6). In particular, the fact that the eigenfunctions form an orthonormal basis for $L^2(\bar{\mathcal{W}})$ is proved in Theorem 6.4; the fact that the eigenvalues are nondecreasing and $\lambda_1 < \lambda_2$ is proved with Theorem 6.3. The eigenvalues and eigenfunctions of Theorem A.16 are used to construct the solution to the boundary value problem (A.16) of Theorem A.15, as proved in the proposion below.

**Proposition A.17** (solution to the sojourn probability BVP). *Let $\mathcal{W} \subset \mathbb{R}^n$ be an open, bounded and connected domain with a smooth boundary $\partial\mathcal{W}$. Then the solution $u\colon \mathbb{R}_{\geq 0} \times \bar{\mathcal{W}} \to [0, 1]$ to the boundary value problem* (A.16) *over $\mathcal{W}$ is given by*

$$u(t, w) = \sum_{i=1}^{\infty} \exp(-\lambda_i t) \phi_i(w) \int_{\mathcal{W}} \phi_i(v) \, \mathrm{d}v. \tag{A.18}$$

*Proof.* To begin, note that (A.16a) is satisfied since

$$\frac{\partial u}{\partial t}(t, w) = -\sum_{i=1}^{\infty} \exp(-\lambda_i t) \lambda_i \phi_i(w) \int_{\mathcal{W}} \phi_i(v) \, \mathrm{d}v \tag{A.19a}$$

$$= \sum_{i=1}^{\infty} \exp(-\lambda_i t) \frac{1}{2} \operatorname{tr} \nabla^2 \phi_i(w) \int_{\mathcal{W}} \phi_i(v) \, \mathrm{d}v \tag{A.19b}$$

$$= \frac{1}{2} \operatorname{tr} \nabla^2_{\mathrm{w}} u(t, w), \tag{A.19c}$$

where to obtain (A.19b) the eigenfunction partial differential equation (A.17a) was used. In addition, the boundary condition (A.16b) is trivially satisfied since $\phi_i(w) = 0$ for all $w \in \partial\mathcal{W}$ due to the eigenvalue problem boundary condition (A.17b).

Finally, to see that the boundary condition (A.16c) is also satisfied by the expression (A.18), recall that $\{\phi_i\}_{i=1}^{\infty}$ is an orthonormal basis for $L^2(\bar{\mathcal{W}})$, according to Theorem A.16. This means that any function $g \in L^2(\bar{\mathcal{W}})$ can be written as an orthogonal projection onto the basis:

$$g(w) = \sum_{i=1}^{\infty} \phi_i(w) \int_{\mathcal{W}} g(v)\phi_i(v)\,\mathrm{d}v. \tag{A.20}$$

If we take $g(w) = 1$ for all $w \in \mathcal{W}$ and substitute (A.20) into (A.18), we get that

$$u(0, w) = \sum_{i=1}^{\infty} \phi_i(w) \int_{\mathcal{W}} \phi_i(v)\,\mathrm{d}v = 1 \qquad \text{for all } w \in \mathcal{W}. \qquad \square$$

The next theorem is a stochastic version of Stokes' theorem. It was first proved by Takahashi and Watanabe (1981, Lem. 2.3) for use in obtaining the Onsager–Machlup functional for diffusions in manifolds (see also Hara and Takahashi, 1996, Lem. 4.3; Capitaine, 2000, Lem. 2). We generalize the original theorems for random differential forms and integrators starting outside the origin. This generalizations are necessary for use in obtaining the Onsager–Machlup functional for diffusions with unknown parameters and initial conditions.

**Lemma A.18** (stochastic Stokes' theorem)**.** *Let $X$ be an $\mathbb{R}^n$-valued semimartingale over $\mathcal{T} := [0, 1]$, adapted to the filtration $\{\mathcal{E}_t\}_{t\geq 0}$; and the function $f\colon \mathcal{T} \times \mathbb{R}^n \times \Omega \to \mathbb{R}^n$ be almost surely differentiable with respect to its first argument and twice continuously differentiable with respect to its second argument. In addition, let the process $F$, defined as*

$$F_t := f(t, X_t, \omega),$$

*be adapted to the filtration $\{\mathcal{E}_t\}_{t\geq 0}$ as well. Then,*

$$\begin{aligned}\int_0^1 F_t^{\mathsf{T}} \circ \mathrm{d}X_t + \bar{f}(0, X_0, \omega)^{\mathsf{T}} X_0 - \bar{f}(1, X_1, \omega)^{\mathsf{T}} X_1 \\ = -\int_0^1 \frac{\partial \bar{f}}{\partial t}(t, X_t, \omega)^{\mathsf{T}} X_t\,\mathrm{d}t - \sum_{i,j=1}^{n} \int_0^1 \frac{\partial \bar{f}^{(i)}}{\partial x^{(j)}}(t, X_t, \omega) \circ \mathrm{d}S_t^{(ij)}\end{aligned} \tag{A.21}$$

*where $S\colon \mathcal{T} \times \Omega \to \mathbb{R}^{n\times n}$ is Lévy's area process.*

$$S_t^{(ij)} := \int_0^t X_\tau^{(i)} \circ \mathrm{d}X_\tau^{(j)} - \int_0^t X_\tau^{(j)} \circ \mathrm{d}X_\tau^{(i)}$$

*and the function $\bar{f}\colon \mathcal{T} \times \mathbb{R}^n \times \Omega \to \mathbb{R}^n$ is defined as*

$$\bar{f}(t, x, \omega) := \int_0^1 f(t, \tau x, \omega)\,\mathrm{d}\tau. \tag{A.22}$$

*Proof.* To begin, define the $X^N$ and $F^N$ processes, for all $N \in \mathbb{N}$, as the piecewise linear interpolations of $X$ and $F$ over the sets $\mathcal{T}_N := \{0, 1/N, \dots, 1\}$. In addition, define the polygonal space-time 2-chains $\mathbb{c}_N$ as

$$\mathbb{c}_N := \{(t, uX_t^N), 0 \le u \le 1, 0 \le t \le 1\} \tag{A.23}$$

and the almost-surely differentiable space-time 1-form $\beta$ as

$$\beta = \sum_{i=1}^n f^{(i)}(t, x, \omega)\,\mathrm{d}x^{(i)}. \tag{A.24}$$

Applying the classical Stokes' theorem, we have that

$$\int_{\partial \mathbb{c}_N} \beta = \int_{\mathbb{c}_N} \mathrm{d}\beta, \tag{A.25}$$

where the exterior derivative of $\beta$ is the differential 2-form given by

$$\mathrm{d}\beta = \sum_{i=1}^n \frac{\partial f^{(i)}}{\partial t}(t, x, \omega)\,\mathrm{d}t \wedge \mathrm{d}x^{(i)} + \sum_{i,j=1}^n \frac{\partial f^{(i)}}{\partial x^{(j)}}(t, x, \omega)\,\mathrm{d}x^{(j)} \wedge \mathrm{d}x^{(i)}, \tag{A.26}$$

where $\wedge$ is the wedge product.

From the definition of the chain and the differential 1-form in (A.23) and (A.24), we have that the left-hand side of (A.25) is given by

$$\int_{\partial \mathbb{c}_N} \beta = \int_0^1 f(0, uX_0, \omega)^\mathsf{T} X_0\,\mathrm{d}u + \int_0^1 f(t, X_t^N, \omega)^\mathsf{T}\,\mathrm{d}X_t^N + \int_1^0 f(1, uX_1, \omega)^\mathsf{T} X_1\,\mathrm{d}u, \tag{A.27}$$

where the second integral on the right-hand side of (A.27) can be taken as either a Riemann–Stieltjes, an Itō or a Stratonovich integral since $X^N$ is a bounded variation process. The integrals of the right-hand side of (A.27) correspond to a closed cycle at the boundary of $\mathbb{c}_N$ obtained by, starting at $t$ and $u$ at the origin, increasing $u$ up to one, then increasing $t$ up to one while keeping $u = 1$, then decreasing $u$ down to zero while keeping $t = 1$, then drecreasing $t$ down to zero while $u = 0$. We note that integral corresponding to the last edge of the path is zero since $uX_t^N$ does not change when $u = 0$. Using the function $\bar{f}$ defined in (A.22), (A.27) can be simplified to

$$\int_{\partial \mathbb{c}_N} \beta = \bar{f}(0, X_0, \omega)^\mathsf{T} X_0 + \int_0^1 f(t, X_t^N, \omega)^\mathsf{T}\,\mathrm{d}X_t^N - \bar{f}(1, X_1, \omega)^\mathsf{T} X_1. \tag{A.28}$$

Similarly, from the definition of the chain in (A.23) and the formula of the exterior derivative of $\beta$ in (A.26), we have that the right-hand side of (A.25)

is given by

$$\int_{\mathbb{c}_N} \mathrm{d}\beta = -\int_0^1\!\!\int_0^1 \frac{\partial f}{\partial t}(t, uX_t, \omega)^{\mathsf{T}} X_t^N \,\mathrm{d}u\,\mathrm{d}t \\ - \sum_{i,j=1}^{n} \int_0^1\!\!\int_0^1 \frac{\partial f^{(i)}}{\partial x^{(j)}}(t, uX_t^N, \omega) \left( X_t^{N(i)} \dot{X}_t^{N(j)} - X_t^{N(j)} \dot{X}_t^{N(i)} \right) u \,\mathrm{d}u\,\mathrm{d}t. \quad \text{(A.29)}$$

Using the function $\bar{f}$ defined in (A.22), (A.29) can be simplified to

$$\int_{\mathbb{c}_N} \mathrm{d}\beta = -\int_0^1 \frac{\partial \bar{f}}{\partial t}(t, X_t^N, \omega)^{\mathsf{T}} X_t^N \,\mathrm{d}t \\ - \sum_{i,j=1}^{n} \int_0^1 \frac{\partial \bar{f}^{(i)}}{\partial x^{(j)}}(t, X_t^N, \omega) \left( X_t^{N(i)} \,\mathrm{d}X_t^{N(j)} - X_t^{N(j)} \,\mathrm{d}X_t^{N(i)} \right), \quad \text{(A.30)}$$

where the last two integral can be interpreted as either a Riemann–Stieltjes, an Itō or a Stratonovich integral.

Next, by letting $N \to \infty$ we have that the integrals with respect to $X^N$ on the right-hand side of both (A.28) and (A.30) converge almost surely to Stratonovich integrals with respect to $X$ (Sussmann, 1978). Substituting (A.28) and (A.30) into (A.25) we then obtain (A.21). □

**Lemma A.19** (Capitaine 2000, Lem. 3). *Let the processes $A\colon \mathcal{T} \times \Omega \to \mathbb{R}$ and $B\colon \mathcal{T} \times \Omega \to \mathbb{R}^n$ be two continuous square-integrable martingales over $\mathcal{T} := [0,1]$ with respect to the filtration $\{\mathcal{E}_t\}_{t\geq 0}$, such that $B$ has the predictable representation property and $[\![A, B]\!] = 0$, i.e., $A$ and $B$ are orthogonal martingales. Additionally, let $\mathbb{B} \in \mathcal{E}$ be an event from the $\sigma$-algebra generated by $B$. Then, for all $\mathbb{R}$-valued process $F$ adapted to $\{\mathcal{E}_t\}_{t\geq 0}$ for which there exists $c \in \mathbb{R}$ such that*

$$\sup_{\omega\in\mathbb{B}} \sup_{t\in\mathcal{T}} [\![A]\!]_1(\omega) + |F_t(\omega)| \leq c,$$

*we have that*

$$\mathrm{E}\left[\exp\left(\int_0^1 F_t \,\mathrm{d}A_s\right) \middle| \mathbb{B}\right] \leq \left(\mathrm{E}\left[\exp\left(2\int_0^1 F_t^2 \,\mathrm{d}\,[\![A]\!]_s\right) \middle| \mathbb{B}\right]\right)^{1/2}.$$

**Lemma A.20** (Georgii, 2008, Prop. 9.1). *Let $A$ be an $\mathbb{R}^d$-valued random variable admitting a probability density. If $q\colon \mathbb{R}^d \to \mathbb{R}^d$ is a $C^1$ diffeomorphism and the $\mathbb{R}^d$-valued random variable $B$ is defined as $B := q(A)$, then $B$ admits the probability density given by*

$$p_B(b) = \left|\det \nabla q^{-1}(b)\right| p_A(q^{-1}(b)),$$

*where $\det \nabla q^{-1}(b)$ is the Jacobian determinant of the inverse of $q$, evaluated at $b$.*

# Index